\documentclass[10pt,a4paper,twoside]{scrartcl}
\usepackage[utf8]{inputenc}
\usepackage[english]{babel}
\usepackage{lmodern}
\usepackage{amsmath,amssymb,mathtools,amsfonts}
\usepackage{graphicx}
\usepackage{xcolor}
\usepackage{subcaption}
\usepackage{paralist}
\usepackage{hyperref}
\usepackage[capitalize, nameinlink, noabbrev]{cleveref}
\crefname{equation}{}{}
\usepackage{autonum}
\usepackage{wrapfig}
\usepackage{lscape}

\usepackage{scrlayer-scrpage}
\ohead{N. Schneider, V. Michel, K. Sigloch and E. Totten: \\ A matching pursuit approach to the geophysical inverse problem of seismic travel time tomography under the ray theory approximation}

\newcommand{\T}{\mathcal{T}}

\newcommand{\ball}{\mathbb{B}}
\newcommand{\real}{\mathbb{R}}
\newcommand{\nat}{\mathbb{N}}
\newcommand{\intg}{\mathbb{Z}}
\newcommand{\dic}{\mathcal{D}}

\newcommand{\RFMP}{\mathrm{RFMP}}
\newcommand{\ROFMP}{\mathrm{ROFMP}}

\newcommand{\Hs}[1]{\mathcal{H}_#1\left(\ball\right)}
\newcommand{\Lp}[1]{\mathrm{L}^#1\left(\ball,\real\right)}
\newcommand{\lp}[1]{\mathrm{L}^#1\left(\ball,\real^3\right)}
\newcommand{\lon}{\varphi}
\newcommand{\lat}{\theta}
\newcommand{\eps}{\epsilon}
\newcommand{\era}{\varepsilon^r}
\newcommand{\ephi}{\varepsilon^\lon}
\newcommand{\ete}{\varepsilon^t}
\newcommand{\trans}{\mathrm{T}}
\newcommand{\proj}[1]{\mathcal{P}_{#1}}
\newcommand{\RAD}{\mathbf{R}}
\newcommand{\intd}{\mathrm{d}}
\newcommand{\pdervr}{\frac{\partial}{\partial r}}
\newcommand{\pdervlon}{\frac{\partial}{\partial \lon}}

\newcommand{\pdervt}{\frac{\partial}{\partial t}}
\DeclareMathOperator*{\sgn}{sgn}
\DeclareMathOperator*{\atanh}{atanh}

\title{A matching pursuit approach \\ to the geophysical inverse problem \\ of seismic travel time tomography \\ under the ray theory approximation}
\author{N. Schneider$^\ast$, V. Michel\footnote{Geomathematics Group Siegen, University of Siegen, michel@mathematik.uni-siegen.de, naomi.schneider@mathematik.uni-siegen.de} ,\ K. Sigloch\footnote{Université Côte d'Azur, CNRS, Observatoire de la Côte d'Azur, IRD, Géoazur, Sophia Antipolis, France}\phantom{x} and E. J. Totten\footnote{Earth Sciences Department, University of Oxford, UK, now at Dublin Institute for Advanced Studies}}
\date{}

\usepackage{natbib}
\begin{document}
\maketitle

\begin{abstract}
Seismic travel time tomography is a geophysical imaging method to infer the 3-D interior structure of the solid Earth. Most commonly formulated as a linear(ized) inverse problem, it maps differences between observed and expected wave travel times to interior regions where waves propagate faster or slower than the expected average. The Earth’s interior is typically parametrized by a single kind of localized basis function. Here we present an alternative approach that uses matching pursuits on large dictionaries of basis functions.\\
Within the past decade the (Learning) Inverse Problem Matching Pursuits ((L)IPMPs) have been developed. They combine global and local trial functions. An approximation is built in a so-called best basis, chosen iteratively from an intentionally overcomplete set or dictionary. In each iteration, the choice for the next best basis element reduces the Tikhonov--Phillips functional. This is in contrast to classical methods that use either global or local basis functions. The LIPMPs have proven its applicability in inverse problems like the downward continuation of the gravitational potential as well as the MEG-/EEG-problem from medical imaging.  \\ 
Here, we remodel the Learning Regularized Functional Matching Pursuit (LRFMP), which is one of the LIPMPs, for travel time tomography in a ray theoretical setting. In particular, we introduce the operator, some possible trial functions and the regularization. We show a numerical proof of concept for artificial travel time delays obtained from a contrived model for velocity differences. The corresponding code is available at \textcolor{blue}{\url{https://doi.org/10.5281/zenodo.8227888}} under the licence CC-BY-NC-SA 3.0 DE.
\end{abstract}



\paragraph{Acknowledgments}
The authors gratefully acknowledge the financial support by the German Research Foundation (DFG; Deutsche Forschungsgemeinschaft), project MI 655/14-1. K. Sigloch was supported by the French government through the UCAJEDI Investments in the Future project, reference number ANR-15-IDEX-01. We thank Maria Tsekhmistrenko, PhD, who owns the relevant DETOX code on GitHub and Afsaneh Mohammadzaheri who assisted E.J. Totten with some data in the time of the Corona pandemic. Last but not least, we are grateful for using the HPC Clusters Horus and Omni maintained by the ZIMT of the University of Siegen for our numerical results.

\section{Introduction}
\label{sect:intro}

Seismic travel time tomography serves to infer the 3-D interior structure of the solid Earth. For this, it measures the characteristics of seismic waves that propagate through the Earth’s deep interior, between sources (earthquakes) and receivers (seismometers) that are located at or near the surface. The most widely practised approach solves a linearized inverse problem, where measured travel times of waves propagating through the real, heterogeneous interior deviate moderately from forward-modelled travel times through a spherically symmetric reference Earth model.This reference velocity model is updated with moderate 3-D variations, in a linear inversion step that attributes travel time anomalies (observed minus predicted) to discrete regions where wave velocities must be moderately faster or slower than in the reference model. The quantitative connection is made by means of efficiently computed sensitivity kernels, see \cite{AkiRichards2009,Abel2007,BenMenahemSingh2000,DahlenTromp1998,Dahlenetal2000,Nolet2008}. 

In particular, a \textbf{s}pherically symmetric, \textbf{n}on-\textbf{r}otating, \textbf{e}lastic and \textbf{i}sotropic spherical shell (SNREI) model like the IASP91, see \cite{Kennettetal1990}, is usually used. When comparing modelled travel times with true measurements, we expect a travel time
difference caused by anomalies within the Earth. Specifically, for total body-wave travel times of hundreds to 1000 s, the unexplained travel time delay is typically fractions of a second to several seconds. Travel time tomography aims to approximate these deviations in the velocity field from the corresponding delays. Thus, travel time tomography is an inverse problem. Unfortunately, there are not many theoretical results known about it. In particular, for practical purposes, we are lacking a singular value decomposition. Moreover, it is known that, in practice, the inverse problem is ill-posed because, e.\,g.\, the solution is non-unique. For details, see e.\,g.\, \cite{Abel2007}. Nonetheless, it is approached with an infinite or more accurate, but also more demanding finite frequency strategy, \cite{DahlenTromp1998,Hosseinietal2020,Marqueringetal1998,Marqueringetal1999,Nolet2008,Sigloch2008,Tian2007,Tian2007-2,Yomogida1992}.
As can be seen at \textcolor{blue}{\url{https://www.earth.ox.ac.uk/~smachine/cgi/index.php}}, currently there exists a range of detailed interior Earth velocity deviation models instead of one standard model. 

Taking a closer look at these models raises a few questions: tetrahedra, voxels, spherical harmonics are chosen mostly for computational convenience and not due to a physical motivation. Moreover, mantle anomalies are a-priori expected to be wider in horizontal dimensions than in the vertical dimension. This puts not only the type of chosen basis into question but also their specific geometry. Further, the characteristics of  the structure of the Earth are still unknown. Hence, it might be potentially better to restrict basis functions as less as possible but let the data itself choose suitable functions. This is especially interesting since some areas of the mantle are very well illuminated, others very poorly. However, we do not expect fundamental changes in convection style in different places due to the general viscosity of the mantle. Hence, unsampled blanks might be filled more physically plausibly if the approximation method is free to choose types of basis functions taking into account the whole data distribution at once. This motivated us to apply a different type of approximation method for inverse problems to travel time tomography which appears to be suitable to tackle these challenges of parameterization, underdeterminedness and regularization.

Here, we propose to use the Learning Regularized Functional Matching Pursuit (LRFMP) which is one of the (Learning) Inverse Problem Matching Pursuit ((L)IPMP) algorithms. The LRFMP is the realization of the RFMP with a learning add-on. The RFMP algorithm iteratively approximates the solution of an inverse problem by minimizing the corresponding Tikhonov--Phillips functional. The unique characteristic of the (L)IPMPs as well as the RFMP, is that the obtained approximation will be given in a so-called best basis of dictionary elements. The dictionary is an intentionally redundant set of diverse trial functions. Here, we consider orthogonal polynomials and linear tesseroid-based finite element hat functions (FEHFs). The best basis is usually also made of all types of trial functions given in the dictionary. The learning add-on enables the method to choose the best basis out of infinitely many trial functions. More details on the methods can be found in 
\cite{Berkeletal2011,Fischer2011,Fischeretal2012,Fischeretal2013-1,Fischeretal2013-2,Guttingetal2017,Kontak2018,Kontaketal2018-2,Kontaketal2018-1,Michel2015-2,Micheletal2017-1,Micheletal2018-1,Micheletal2014,Micheletal2016-1,Leweke2018,Prakashetal2020,Schneider2020,Schneideretal2022,Telschow2014,Telschowetal2018}. 
These publications also show that these methods have already been successfully applied to inverse problems in geodesy and medical imaging as well as for normal mode tomography. Here, we will shortly introduce them as well in \cref{sect:ipmps}, in particular with an emphasis on adjustments made for travel time tomography. Regarding the theory, a focus is set on the derivation of certain gradients and inner products needed for this adjustment given in \cref{ssect:tfcs:reg}. These computations are provided as supplementary material in the appendix of this paper.

We demonstrate this new method on inversions for synthetic (i.e., invented) whole-mantle test structures, using realistic sampling geometries, i.e. actual earthquake-to-receiver paths that have yielded high-quality travel time data for tomography in the ISC-EHB catalogue. For the latter, see e.\,g.\, \cite{ISC-EHB,Westonetal2018}. 


In \cref{sect:problem}, we introduce the mathematical inverse problem of the travel time tomography. Afterwards, we present our choices of dictionary elements in \cref{sect:tfcs}. Further, we compute regularization terms related to a Sobolev space of these trial functions. Next, we introduce the LRFMP in \cref{sect:ipmps} and show a first proof of concept in \cref{sect:numerics}. In \cref{sect:app}, the interested reader will find more details on certain mathematical derivations regarding the trial functions and the methods. The corresponding code is available at \textcolor{blue}{\url{https://doi.org/10.5281/zenodo.8227888}} under licence the CC-BY-NC-SA 3.0 DE.

\subsection{Notation}
\label{ssect:intro:nota}
The set of positive integers is denoted by $\nat$. If we include 0, we use $\nat_0$. For integers, we write $\intg$. For the set of real numbers, we use $\real$. The $d$-dimensional Euclidean space is called $\real^d$. Only positive real numbers are collected in $\real_+$. For $x\in\real^3$, we can use the common spherical coordinate transformation
\begin{align}
x(r,\lon,\lat) = \left( \begin{matrix}
r\sqrt{1-t^2}\cos(\lon)\\
r\sqrt{1-t^2}\sin(\lon)\\
rt
\end{matrix} \right)
\end{align}
for the radius $r\in\real_+$, the longitude $\lon \in [0,2\pi[$ and the latitude $\lat \in [0,\pi]$ which yields the polar distance $\cos(\lat) = t \in [-1,1]$. A radius is denoted by $\RAD \in \real_+$ and the corresponding ball with radius $\RAD$ is defined by $\ball_\RAD \coloneqq \{x \in \real^3\ \colon |x|\leq \RAD\}$. If $\RAD = 1$, we also abbreviate $\ball \coloneqq \ball_1$. Thus, a point $\xi(\lon,t) \in \partial \ball \subset \real^3$ is given by
\begin{align}
\xi(\lon,t) = \left( \begin{matrix}
\sqrt{1-t^2}\cos(\lon)\\
\sqrt{1-t^2}\sin(\lon)\\
t
\end{matrix} \right).
\end{align}
Moreover, a well-known local orthonormal basis in $\real^3$ is defined by the vectors
\begin{align}
\era(\lon,t) = \left( \begin{matrix}
\sqrt{1-t^2}\cos(\lon)\\
\sqrt{1-t^2}\sin(\lon)\\
t
\end{matrix} \right),
\ 
\ephi(\lon,t) = \left( \begin{matrix}
-\sin(\lon)\\
\phantom{-}\cos(\lon)\\
0
\end{matrix} \right),
\ 
\ete(\lon,t) = \left( \begin{matrix}
-t\cos(\lon)\\
-t\sin(\lon)\\
\sqrt{1-t^2}
\end{matrix} \right),
\end{align}
see, for instance, \cite{Michel2013}. 

\section{Seismic travel time tomography}
\label{sect:problem}

Seismic travel time tomography estimates the propagation velocity of seismic (P-)waves as a function of spatial location inside the solid Earth. The sensitivity of travel time measurements to Earth structures along the wave propagation path is volumetrically extended in practice, but can often be reasonably approximated as an (infinitely narrow) ray, in analogy to optical ray theory. Here we consider only (seismic) ray theory for P-waves, although our approach could be extended to less approximative sensitivity modelling, such as Born/finite-frequency approximations, see e.\,g.\,\cite{Dahlenetal2000}, and to S- or other wave types. As we aim here for a proof of concept, we consider the accuracy of the classical ray theoretical infinite-frequency approach as adequate. For one seismic ray $\widetilde{X}$ between seismic source and receiver, the mathematical relationship between its travel time $\psi$ and the (P-)wave velocity $c$ (and its slowness $S$, respectively), is given by the Eikonal equation, see e.g. \cite{DahlenTromp1998,Michel2020,Nolet2008}, 
\begin{align}
\frac{1}{c(\widetilde{X}(s))} = S(\widetilde{X}(s)) = |\nabla_{\widetilde{X}} \psi(\widetilde{X}(s))|,
\label{eq:Eikonal}
\end{align}
where $\widetilde{X} \colon [0,R]\to\real^3$ and the arc length $s$ is chosen to parametrize $\widetilde{X}$. We repara\-meterize the ray here by a parameter $t\in[0,1]$ such that $X\colon [0,1]\to\real^3$ has the same curve, that is $X([0,1])=\widetilde{X}([0,R])$. The specific choice of the parameter transformation $s\leftrightarrow t$ depends on the implementation of the numerically calculated ray. The Eikonal equation yields the non-linear inverse problem 
\begin{align}
\int_{\widetilde{X}} S(\widetilde{X}(s))\ \intd \sigma(\widetilde{X}(s)) = \psi(\widetilde{X}(R)) \
&\Leftrightarrow  \
\int_0^1 S(X(t))\left|X'(t)\right| \intd t = \psi(X(1))
\label{eq:NLIPTT}
\intertext{with}
s(t) &= \int_0^t|X'(\tau)|d\tau.
\label{eq:arclen}
\end{align} 
Due to the difficulty of non-linear inverse problems, \cref{eq:NLIPTT} is usually linearized: instead of approximating the slowness $S$ itself, we consider the deviation $\delta S$ between reality and a reference model, see e.g. \cite{Abel2007,Nolet2008}. This yields the seismic ray perturbation operator (SPO) given by
\begin{align}
(\T\ \delta S)\left(\widetilde{X}_\mathrm{ref}\right) 
\coloneqq  \int_0^R \delta S\left(\widetilde{X}_\mathrm{ref}(s)\right)\ \intd s
= \delta \psi\left(\widetilde{X}_\mathrm{ref}(R)\right)
\label{eq:operatorgeneral}
\end{align}
where $\widetilde{X}_\mathrm{ref}$ stands for the ray obtained from a reference model and parameterized with the arc length $s$ between 0 and its total arc length $R$. 
Moreover, $\delta \psi = \psi_{\mathrm{obs}} - \psi_{\mathrm{ref}}$ is the delay of the observed and (due to the reference model) expected travel time. The more rays we consider, the better our approximation of $\delta S$ can be. Thus, in practice, we use a family of rays $\{\widetilde{X}_{\mathrm{ref},i}\}_{i=1,...,\ell}$ with respect to a set of source-receiver pairs 
\begin{align}
\left\{\left(\widetilde{X}_{\mathrm{ref},i}(0),\widetilde{X}_{\mathrm{ref},i}(R)\right)\right\}_{i=1,...,\ell}
\label{eq:sourcereceiverpairs}
\end{align} 
and consider the Discretized SPO (DSPO)
\begin{align}
\left(\T_\daleth^i\ \delta S\right)\left(X_{\mathrm{ref},i}\right) 
&\coloneqq \int_0^{R} \delta S\left(\widetilde{X}_{\mathrm{ref},i}(s)\right) \intd s = \delta \psi_i\left(\widetilde{X}_{\mathrm{ref},i}(R)\right)\\
&= \int_0^{1} \delta S\left(X_{\mathrm{ref},i}(t)\right)\left|X'_{\mathrm{ref},i}(t)\right| \intd t = \delta \psi_i\left(X_{\mathrm{ref},i}(1)\right),\\
\T_\daleth\ \delta S &= \delta \psi_\daleth
\label{eq:operatordiscret}
\end{align}
with $s$ as in \cref{eq:arclen}, $s(1)=R$, $\T_\daleth \coloneqq (\T_\daleth^i)_{i=1,...,\ell}$ and $\delta \psi_\daleth \coloneqq (\delta \psi_i(X_{\mathrm{ref},i}(1)) )_{i=1,...,\ell}.$ Our aim is then to construct an approximation 
\begin{align}
\delta S \approx \sum_{i=1}^I \alpha_id_i
\end{align}
for $I\in\nat,\ \alpha_i\in\real$ and some trial functions $d_i$. 
Note that, in seismology, instead of considering $\delta S$, the velocity anomaly (deviation from the reference velocity) is often considered. In our approach, this means, we would have to reformulate our approximation in the following way: if we obtain $\delta S$ as the linear combination, we approximate
\begin{align}
\delta S = \frac{1}{c} - \frac{1}{c_{\mathrm{ref}}} = - \frac{c-c_{\mathrm{ref}}}{cc_{\mathrm{ref}}}.
\label{eq:deltaS}
\end{align}
This reformulates to 
\begin{align}
\delta S + \frac{1}{c_{\mathrm{ref}}} = \frac{1}{c}.
\label{eq:c_reciproc}
\end{align}
Thus, after approximating $\delta S$, we would transform it via
\begin{align}
-\frac{\delta S}{\delta S + c_{\mathrm{ref}}^{-1}}
= -\left(- \frac{c-c_{\mathrm{ref}}}{cc_{\mathrm{ref}}} \right)c 
= \frac{c-c_{\mathrm{ref}}}{c_{\mathrm{ref}}}
\eqqcolon \frac{\mathrm{d}c}{c}
\end{align}
for a better comparability in the geophysical community. Since we are using here a non-real, artificial deviation model (see \cref{ssect:numerics:genset}), we abstain from this explicit reformulation but consider the numerator naturally in our resolution test.

\section{Trial functions}

\label{sect:tfcs}
We utilize two types of trial functions in our methods: global orthogonal polynomials and compactly supported linear finite element hat functions. We first present their definition and give examples of them.

\subsection{Definition and examples}
\label{ssect:tfcs:def}

\subsubsection{Tesseroid-based finite element hat functions}
\label{sssect:tfcs:def:FE}

\begin{figure}[htbp]
\centering
\includegraphics[width=.32\textwidth]{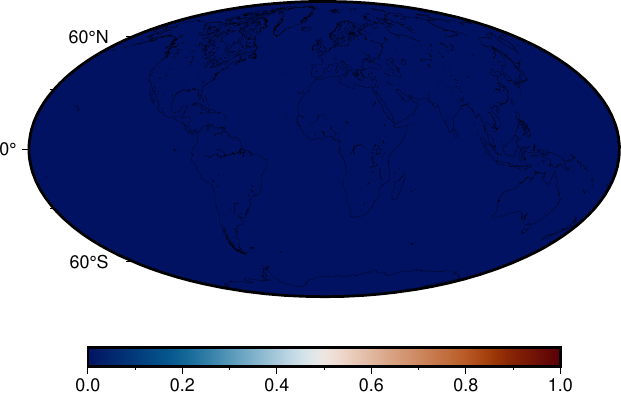}
\includegraphics[width=.32\textwidth]{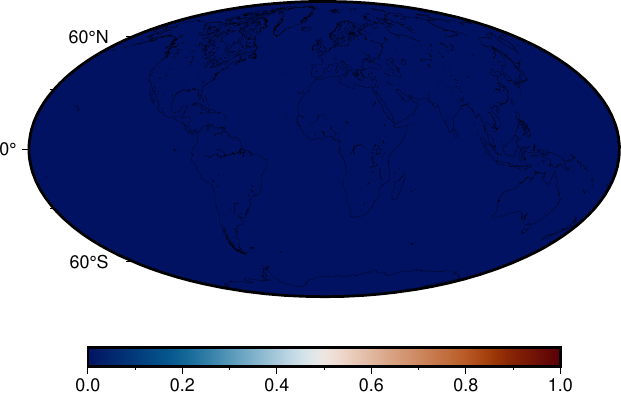}
\includegraphics[width=.32\textwidth]{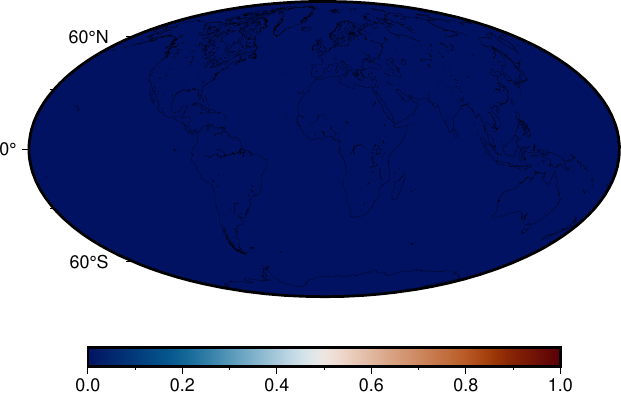}
\includegraphics[width=.32\textwidth]{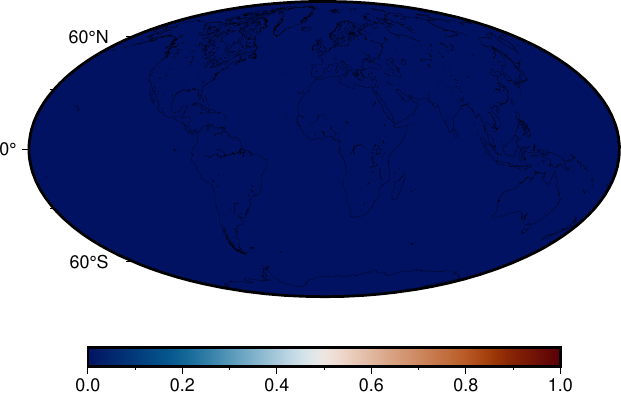}
\includegraphics[width=.32\textwidth]{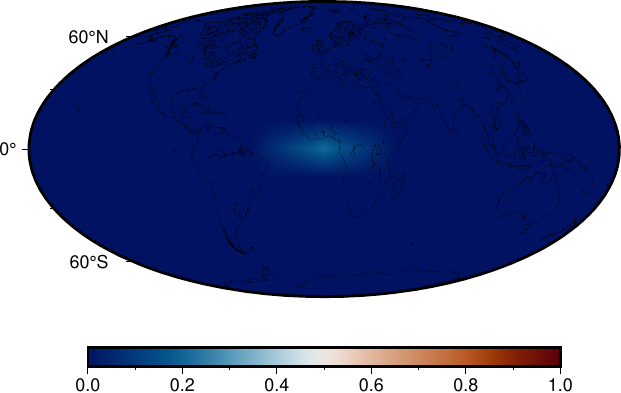}
\includegraphics[width=.32\textwidth]{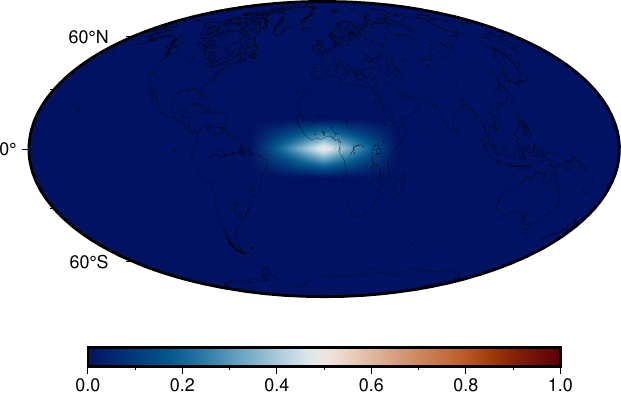}
\includegraphics[width=.32\textwidth]{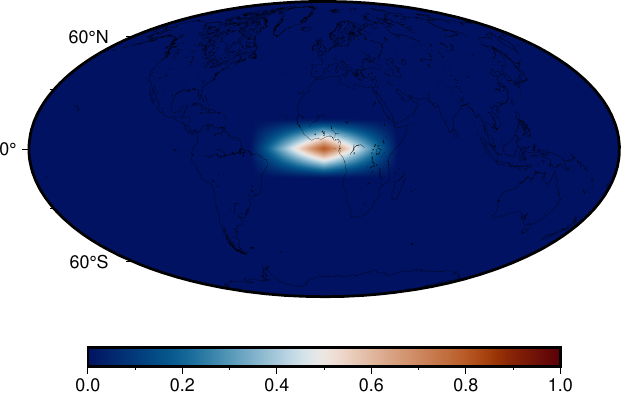}
\includegraphics[width=.32\textwidth]{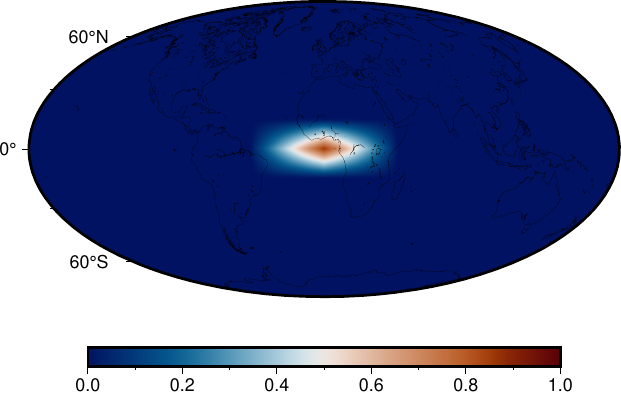}
\includegraphics[width=.32\textwidth]{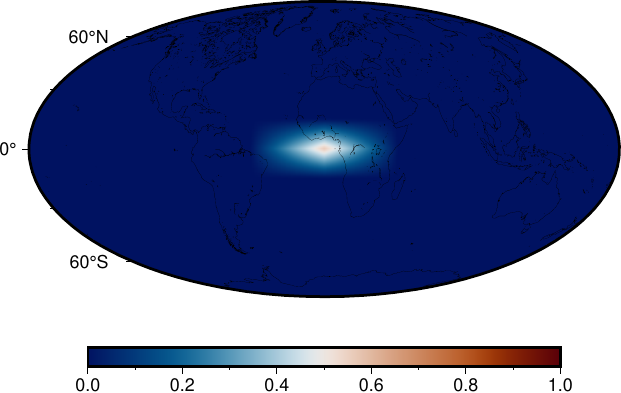}
\includegraphics[width=.32\textwidth]{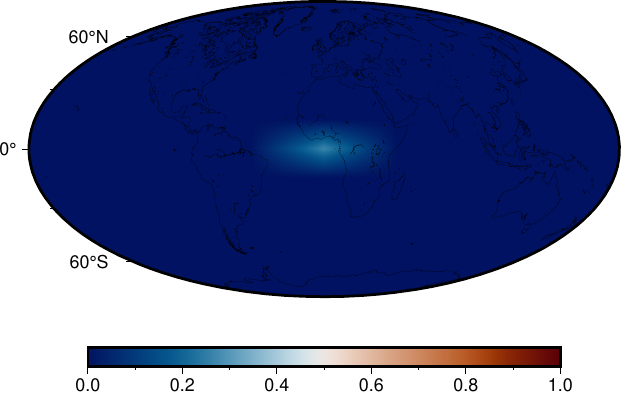}
\includegraphics[width=.32\textwidth]{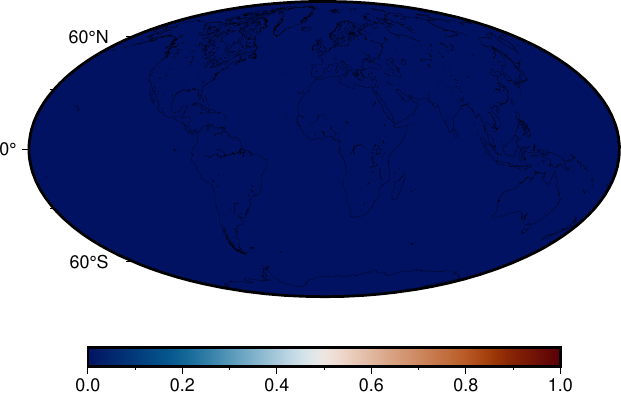}
\includegraphics[width=.32\textwidth]{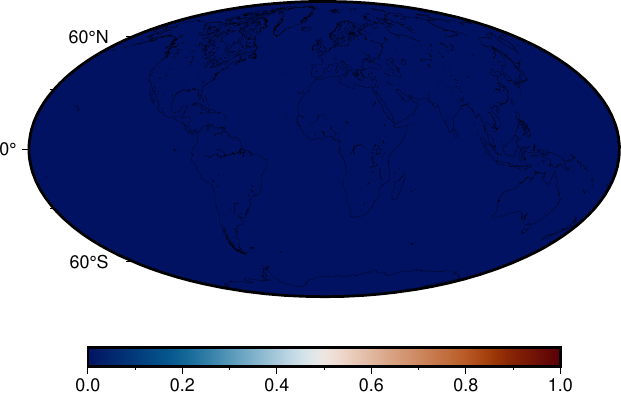}
\caption{Example of an FEHF. We show $N_{(5096.8,0,0),(955.65,\pi/4,0.25)}$ for the depth slices (left, middle and right) at the following radial distances to the centre of the Earth in km: 3193.1, 3482.0 and 3770.9 (first row), 4059.8, 4348.7 and 4637.6 (second row), 4926.5, 5215.4 and 5504.3 (third row), 5793.2, 6082.1 and 6371.0 (last row). The colour scales are adjusted for a better comparability.}
\label{fig:FEHF}
\end{figure}

For a general introduction to finite elements, see e.g. \cite{Braess2007,Grossmannetal2007,Johnson2009,Schwarzetal2011}. Regarding trial functions, we are interested in composite hat functions on tesseroids as the finite elements. Thus, in the sequel, we speak of finite element hat functions (FEHFs).

In seismology, using tetrahedra as finite elements are common, see e.g. \cite{Hosseinietal2020,Sigloch2008,Tian2007,Tian2007-2}. We change the underlying finite element structure here because, from a seismological perspective, there is no physical need to use tetrahedra. Tesseroids, see \cite{Fukushima2018}, are rectangular parallelepipeds defined by upper and lower bounds with respect to the radius, the longitude and the latitude. Thus, a tesseroid may be even more reasonable than tetrahedra: such an underlying structure would be a realization of the interior of the Earth being structured or layered by gravitation. Moreover, the depth of the mantle ($\approx$ 2 900 km) is much smaller than the circumference of the Earth ($\approx$ 40 000 km). Hence, the size of the heterogeneities in the mantle is smaller in the depth than it is latitudinally and longitudinally.

For the definition of a tesseroid, we obviously have six degrees of freedom. Due to the composite nature of the FEHFs, we move from the boundary-based view to a centre-based one. Then, the degrees of freedom are its radial centre $R$ and its distance to each side $\Delta R$ as well as the longitudinal and latitudinal centres $\Phi$ and $T$ and their respective distances to each side $\Delta \Phi$ and $\Delta T$. The natural constraints of a tesseroid are
\begin{align}
R &\in [\rho\RAD,\RAD] = [R_{\mathrm{min}},R_{\mathrm{max}}],& 0<\eps_R&\leq \Delta R\leq \RAD/2,\\
\Phi &\in [0,2\pi]  = [\Phi_{\mathrm{min}},\Phi_{\mathrm{max}}],& 0<\eps_\Phi &\leq \Delta \Phi\leq \pi, \label{def:N-A-DA-constraints}\\
T &\in [-1+\eps_T,1-\eps_T] = [T_{\mathrm{min}},T_{\mathrm{max}}],& 0<\eps_T&\leq \Delta T \leq 0.5.
\end{align}
Note that, with $\rho\in [0,1]$, we control the maximal possible depth of the tesseroids. For instance, if we set $\rho = 3482/6371 = 0.54654$, we allow the FEHF as deep as the core-mantle boundary. Further note that due to singularities at the poles, see \cref{ssect:app:H1IPs}, we need to stay away from the theoretically possible bounds $\pm 1$ in the latitudinal case. We can identify (in polar coordinates) the ball with the domain
\begin{align}
D &\coloneqq [0,\RAD] \times [P,P+2\pi] \times [-1,1]\quad \mathrm{with}  \quad
P\coloneqq\left\lfloor \frac{\Phi-\Delta \Phi}{\pi} \right\rfloor \pi,
\end{align}
the difference domain as
\begin{align}
\Delta D &\coloneqq [\eps_R,\RAD/2] \times [\eps_\Phi,\pi] \times [\eps_T,0.5]\\
&=  [\eps_R,R_{\mathrm{max}}/2] \times [\eps_\Phi,\Phi_{\mathrm{max}}/2] \times [\eps_T, (T_{\mathrm{max}}+\eps_T)/2]
\end{align}
and the tesseroid as
\begin{align}
E &\coloneqq [\max(R_{\mathrm{min}},R-\Delta R),\min(R_{\mathrm{max}},R+\Delta R)] \\ &\qquad \times [\Phi-\Delta \Phi,\Phi+\Delta \Phi] \\ &\qquad \times [\max(T_{\mathrm{min}},T-\Delta T),\min(T_{\mathrm{max}},T+\Delta T)].
\label{def:tesseroid}
\end{align}
For the sake of brevity, let $a\coloneqq(a_j)_{j=1,2,3} \coloneqq (r,\lon,t),\
A\coloneqq(A_j)_{j=1,2,3} \coloneqq (R,\Phi,T),\
\Delta A\coloneqq(\Delta A_j)_{j=1,2,3} \coloneqq (\Delta R,\Delta \Phi,\Delta T),$
\begin{align}
A_{\mathrm{min}} &\coloneqq (\rho\RAD,P,-1+\eps_T) = (R_{\mathrm{min}},\Phi_{\mathrm{min}},T_{\mathrm{min}}),\\
A_{\mathrm{max}} &\coloneqq (\RAD,P+2\pi,1-\eps_T) =  (R_{\mathrm{max}},\Phi_{\mathrm{max}},T_{\mathrm{max}}) 
\intertext{and the Cartesian product}
\mathrm{supp}_{A,\Delta A} (a) &\coloneqq \prod_{j=1}^3 [\max(A_{\mathrm{min},j},A_j-(\Delta A)_j), \min(A_{\mathrm{max},j},A_j+(\Delta A)_j)] = E.
\end{align}
With these definitions, we can consider the FEHFs. Commonly used in seismology are Lagrange finite element basis function, see \cite{Tian2007,Tian2007-2}. Finite element basis functions for cuboids are given, for instance, in \cite{Mazdziarz2010}. We take the linear examples from there. Via translation to $\mathrm{supp}_{A,\Delta A} $ and the general notation just introduced, we then obtain the dictionary elements
\begin{align}
&N_{A,\Delta A}(x(a)) \coloneqq
N_{(R,\Phi,T),(\Delta R,\Delta \Phi, \Delta T)}(r\xi(\lon,t))
\coloneqq 
\chi_{\mathrm{supp}_{A,\Delta A}} (a)\prod_{j=1}^3 \frac{\Delta A_j-|a_j-A_j|}{\Delta A_j} ,
\label{def:N-A-DA}
\end{align}
where $\chi$ denotes the characteristic function and $x(a)=r\xi(\lon,t)\in\ball_\RAD$. The FEHF $N_{A,\Delta A}$ attains its maximum in $A$, which is the centre of the respective volume element. It holds $N_{A,\Delta A}(A)=1.$ The function linearly decreases towards zero when moving towards $A\pm\Delta A$. It is constant zero outside of the volume element. Thus, it is only piecewise constant. An example is given in \cref{fig:FEHF}. Note that, in this example, we see that the hat is clearly visible. However, this shows that, though the theoretical support of an FEHF is a tesseroid, its visible support appears smaller.

\subsubsection{Polynomials}
\label{sssect:tfcs:def:POLY}
For polynomials on a ball with radius $\mathbf{R}$, we consider the system
\begin{align}
&G^{\mathrm{I}}_{m,n,j} (r\xi(\lon,t)) \\
&\coloneqq p_{m,n}P_m^{(0,n+1/2)}\left(I(r)\right) \left(\frac{r}{\mathbf{R}}\right)^{n} Y_{n,j} (\xi(\lon,t))
\label{def:GI_2}\\
&= p_{m,n}P_m^{(0,n+1/2)}\left(2\left(\frac{r}{\mathbf{R}}\right)^2-1\right) \left(\frac{r}{\mathbf{R}}\right)^{n}  q_{n,j}P_{n,|j|}(t) \left\{ \begin{matrix} \sqrt{2}\cos(j\varphi),&j<0\\ 1,&j=0\\ \sqrt{2}\sin(j\varphi),&j>0 \end{matrix} \right.
\label{def:GI_3}\\
&\eqqcolon p_{m,n}P_m^{(0,n+1/2)}\left(I(r)\right) \left(\frac{r}{\mathbf{R}}\right)^{n} q_{n,j}P_{n,|j|}(t) \mathrm{Trig}(j\lon)
\label{def:GI_4}
\end{align}

\begin{figure}[htbp]
\centering
\includegraphics[width=.32\textwidth]{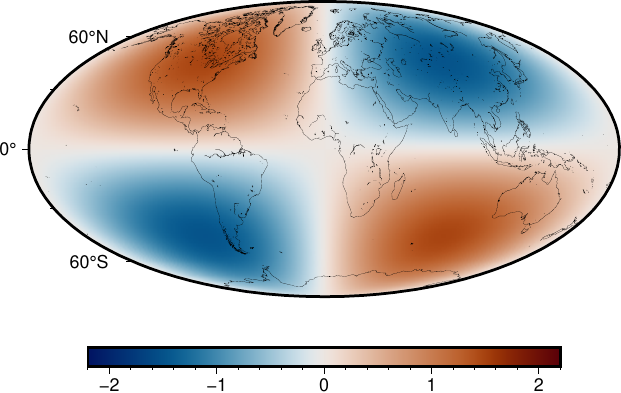}
\includegraphics[width=.32\textwidth]{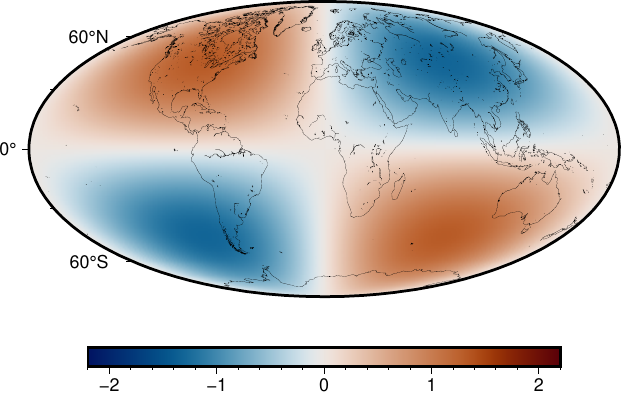}
\includegraphics[width=.32\textwidth]{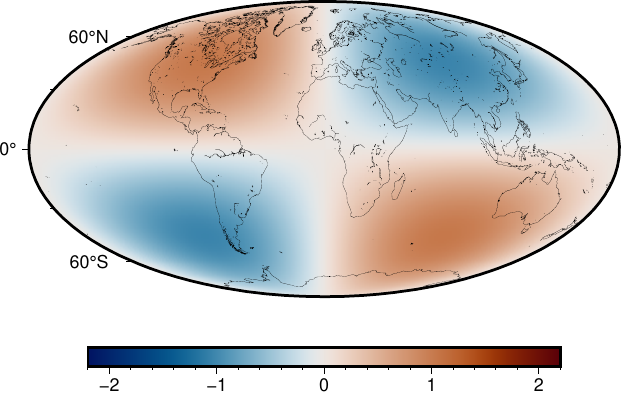}
\includegraphics[width=.32\textwidth]{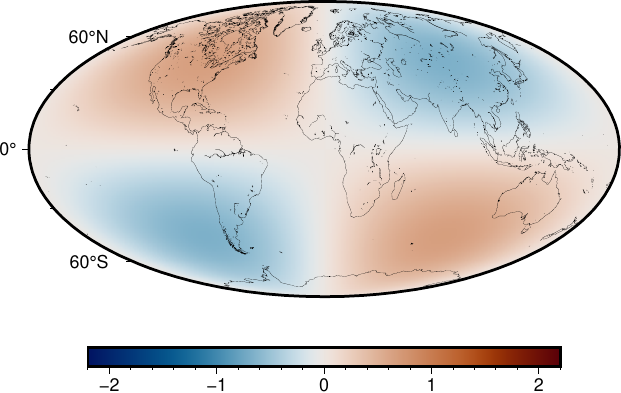}
\includegraphics[width=.32\textwidth]{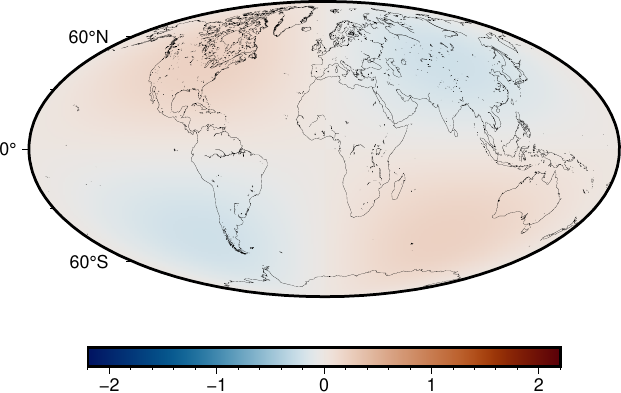}
\includegraphics[width=.32\textwidth]{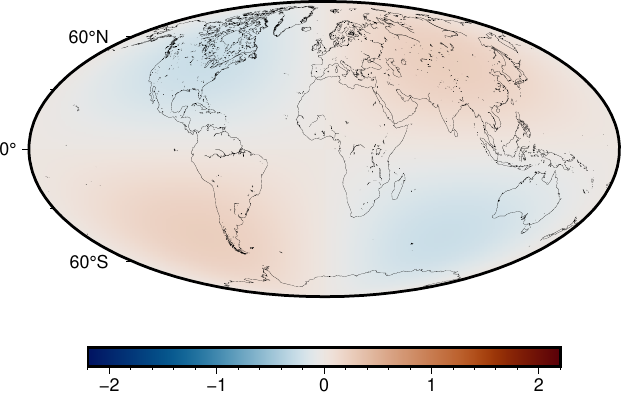}
\includegraphics[width=.32\textwidth]{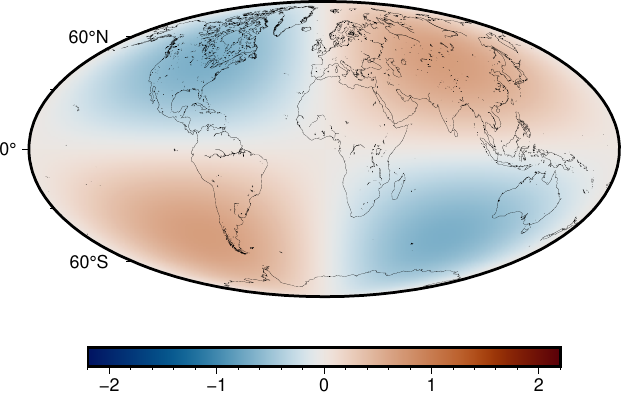}
\includegraphics[width=.32\textwidth]{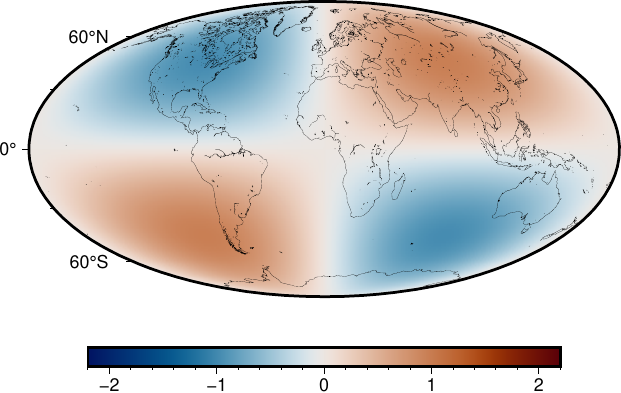}
\includegraphics[width=.32\textwidth]{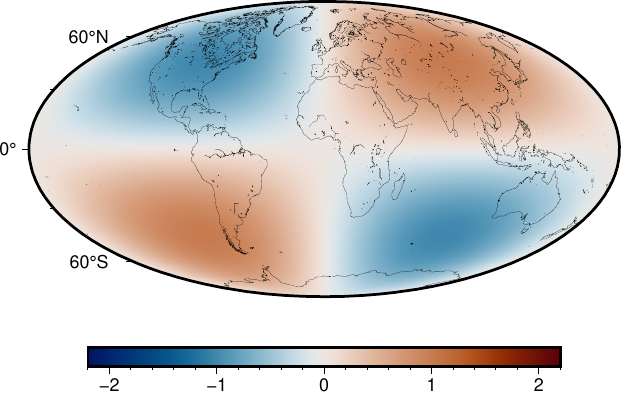}
\includegraphics[width=.32\textwidth]{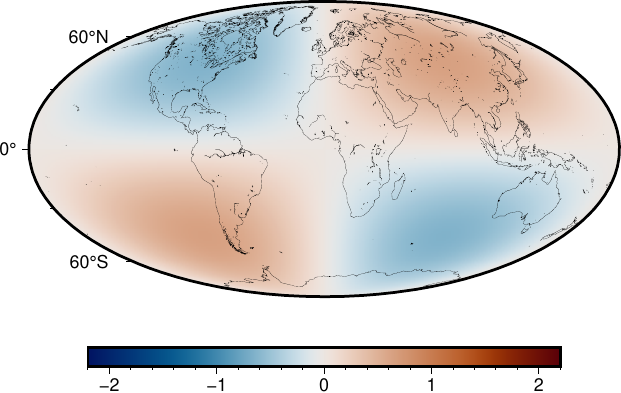}
\includegraphics[width=.32\textwidth]{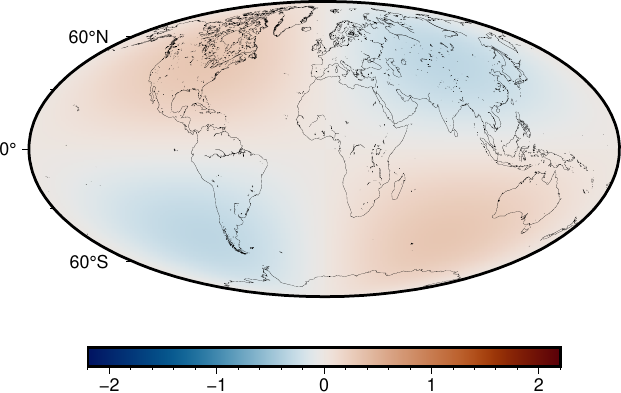}
\includegraphics[width=.32\textwidth]{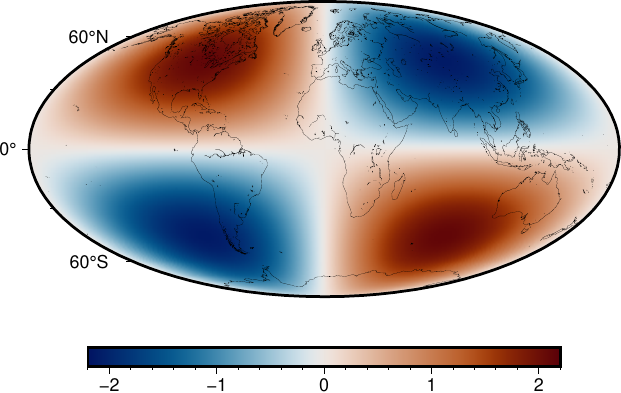}
\caption{Example of a polynomial. We show $G_{2,2,1}^{\mathrm{I}}$ for the depth slices (left, middle and right) at the following radial distances to the centre of the Earth in km: 3193.1, 3482.0 and 3770.9 (first row), 4059.8, 4348.7 and 4637.6 (second row), 4926.5, 5215.4 and 5504.3 (third row), 5793.2, 6082.1 and 6371.0 (last row). The colour scales are adjusted for a better comparability.}
\label{fig:GI}
\end{figure}

with the normalization factors
\begin{align}
&p_{m,n} \coloneqq \sqrt{\frac{4m+2n+3}{\mathbf{R}^3}}
\qquad \mathrm{and} \qquad 
q_{n,j} \coloneqq q_{n,|j|} \coloneqq \sqrt{\frac{2n+1}{4\pi}\frac{(n-|j|)!}{(n+|j|)!}}
\end{align}
for $m,\ n \in \nat_0,\ j \in \intg,\ |j|\leq n$
and where $P_m^{(\alpha,\beta)}$ denotes a Jacobi polynomial, $Y_{n,j}$ a spherical harmonic and $P_{n,|j|}$ an associated Legendre function. As can be seen here, we use fully normalized spherical harmonics in our implementation. For details on those composite functions, see e.g. \cite{AbramowitzStegun1965,Freedenetal1998,Freedenetal2013-1,Freedenetal2009,MagnusOberhettinger1966,Michel2020,Mueller1966,Szegoe1975}. The functions $G^{\mathrm{I}}$ were investigated first by \cite{Ballanietal1993,Dufour1977,Michel2013}. For generalizations of this system, see \cite{DunkelXu2014} as well as \cite{Micheletal2016}. Note that every $G^{\mathrm{I}}$ is an algebraic polynomial in $\real^3$ and is, therefore, well-defined on the whole ball. An example is given in \cref{fig:GI}.

\subsection{A dictionary}
\label{ssect:tfcs:dic}
With these two types of trial functions, we now build a so-called dictionary $\dic$ which is an intentionally redundant set of functions. In the sequel, we work with the following dictionary
\begin{align}
[\mathcal{G}_{M,N}]_{\mathrm{GI}} &\coloneqq \left\{G^{\mathrm{I}}_{m,n,j}\ \colon\ m,\ n \in \nat_0,\ m\leq M,\ n\leq N,\ j\in\intg,\ |j| \leq n\right\},\\
[\mathcal{A}]_{\mathrm{FEHF}} &\coloneqq \left\{N_{A,\Delta A}\ \colon\ A \in D,\ \Delta A \in \Delta D\right\},\\
\dic &\coloneqq \dic^{\mathrm{Inf}} \coloneqq [\mathcal{G}_{M,N}]_{\mathrm{GI}} \cup [\mathcal{A}]_{\mathrm{FEHF}}.
\end{align}
We see that the trial function classes $[\mathcal{G}_{M,N}]_{\mathrm{GI}}$ and $[\mathcal{A}]_{\mathrm{FEHF}}$ are defined via the characteristic parameters of the trial functions. Note that the parameters of the FEHFs are continuous while those of the polynomials are discrete. We allow polynomials up to a maximum radial and angular degree. Theoretically, we could also allow $m,\ n \in \nat_0$ and, thus, let $[\mathcal{G}_{M,N}]_{\mathrm{GI}}$ be infinite. In practice, however, this is not sensible as we will see later. Note, however, that $[\mathcal{A}]_{\mathrm{FEHF}}$ is truly infinite. We discuss the dictionary again in a larger context below when we introduce our approximation methods.

\subsection{Regularization terms}
\label{ssect:tfcs:reg}
As we are considering an ill-posed inverse problem, from a mathematical point of view, it is inevitable to include a regularization. For our approach, we use the Tikhonov--Phillips regularization. Thus, we need to determine a suitable function space for the penalty term. We decided to use the $\Hs{1}$-Sobolev space, see \cite{Adamsetal2003,Bhattacharyya2012,Braess2007,Heuser2006,Schwarzetal2011,Werner2018,Yosida1995}. The regularization terms are then obtained via 
\begin{align}
\langle d_i, d_j \rangle_{\mathcal{H}_1} 
&= \left\langle \mathrm{D}^{(0)}d_i, \mathrm{D}^{(0)} d_j \right\rangle_{\Lp{2}}
+ \left\langle \mathrm{D}^{(1)}d_i, \mathrm{D}^{(1)} d_j \right\rangle_{\lp{2}}\\
&= \left\langle d_i, d_j \right\rangle_{\Lp{2}}
+ \left\langle  \nabla_x d_i, \nabla_x d_j\right\rangle_{\lp{2}}
\label{def:H1IP}
\end{align}
for two dictionary elements $d_i,\ d_j \in \dic$. We have to determine these inner products for the different cases of dictionary elements in $\dic$. We start with the determination of the gradients. For the FEHFs, we obtain 
\begin{align}
&\nabla_{r\xi(\lon,t)} N_{(R,\Phi,T),(\Delta R,\Delta \Phi,\Delta T)}(r,\lon,t) \\&=
\chi_{\mathrm{supp}_{[(R,\Phi,T)-(\Delta R,\Delta \Phi,\Delta T),(R,\Phi,T)+(\Delta R,\Delta \Phi,\Delta T)]}}(r,\lon,t) \\
&\qquad\times \left(
\era \frac{[-\sgn(r-R)]}{\Delta R} \frac{\Delta \Phi-|\lon-\Phi|}{\Delta \Phi}\frac{\Delta T-|t-T|}{\Delta T}\right.\\ 
&\qquad\qquad  + \frac{1}{r}\ephi \frac{1}{\sqrt{1-t^2}} \frac{\Delta R-|r-R|}{\Delta R} \frac{[-\sgn(\lon-\Phi)]}{\Delta \Phi} \frac{\Delta T-|t-T|}{\Delta T}\\
&\qquad\qquad  \left. + \frac{1}{r}\ete \sqrt{1-t^2}\frac{\Delta R-|r-R|}{\Delta R}\frac{\Delta \Phi-|\lon-\Phi|}{\Delta \Phi}\frac{[-\sgn(t-T)]}{\Delta T}
\right).
\end{align}
The derivation of this formula is done straightforwardly and can be found in \cref{ssect:app:nablaFEHF}. For the polynomials, we have
\begin{align}
\nabla_{r\xi(\lon,t)} &G^{\mathrm{I}}_{m,n,j}(r\xi(\lon,t)) \\
&= \left(\frac{\partial}{\partial x_k}G^{\mathrm{I}}_{m,n,j} (r\xi(\lon,t)) \right)_{k=1,2,3} \\
&= p_{m,n}q_{n,j}\left(\left(P_m^{(0,n+1/2)}\left(I(r)\right)\right)'I'(r) \left(\frac{r}{\mathbf{R}}\right)^{n} P_{n,|j|}(t)\mathrm{Trig}(j\lon) \xi(\lon,t)\right.\\
&\qquad\qquad + \frac{n}{\mathbf{R}} P_m^{(0,n+1/2)}\left(I(r)\right) \left(\frac{r}{\mathbf{R}}\right)^{n-1}  P_{n,|j|}(t)\mathrm{Trig}(j\lon) \xi(\lon,t)\\
&\qquad\qquad + \frac{j}{\mathbf{R}} P_m^{(0,n+1/2)}\left(I(r)\right) \left(\frac{r}{\mathbf{R}}\right)^{n-1} \frac{1}{\sqrt{1-t^2}} P_{n,|j|}(t) \mathrm{Trig}(-j\lon) \ephi(\lon,t)\\
&\qquad\qquad \left.+ \frac{1}{\mathbf{R}} P_m^{(0,n+1/2)}\left(I(r)\right) \left(\frac{r}{\mathbf{R}}\right)^{n-1}  \sqrt{1-t^2} P'_{n,|j|}(t)\mathrm{Trig}(j\lon)\ete(\lon,t)\right)\\
&\eqqcolon p_{m,n}q_{n,j}\sum_{p=1}^4
G^{\mathrm{I}}_{m,n,j;p} (r\xi(\lon,t)),
\label{eq:nablaGI}
\end{align}
with 
\begin{align}
I'(r) = \frac{4r}{\RAD^2}.
\label{eq:gradIr}
\end{align}
The computation of this gradient is shown in \cref{ssect:app:nablaGI} in detail. In both cases, note that $\xi=\ \era,\ \ephi$ and $\ete$ are vectors.

Then, we obtain the following values for the inner products. The detailed derivation is given in \cref{ssect:app:H1IPs}. For two polynomials, we have
\begin{align}
&\left\langle G^{\mathrm{I}}_{m,n,j}, G^{\mathrm{I}}_{m',n',j'} \right\rangle_{\mathcal{H}_1(\mathbb{B})}\\
&= \delta_{m,m'}\delta_{n,n'}\delta_{j,j'} + \delta_{n,n'}\delta_{j,j'} p_{m,n}p_{m',n}\\
&\qquad\times\left[\frac{\mathbf{R}\sqrt{2}}{2^{n}}\int_{-1}^{1} \left(P_m^{(0,n+1/2)}(u)\right)'\left(P_{m'}^{(0,n+1/2)}(u)\right)' \left(1+u\right)^{n+3/2} \mathrm{d} u \right. \\
&\qquad\qquad + \frac{\mathbf{R}n}{2^{n}\sqrt{2}}\int_{-1}^{1} \left[ \left(P_m^{(0,n+1/2)}(u)\right)'P_{m'}^{(0,n+1/2)}(u) \right. \\ &\qquad\qquad\qquad\qquad\qquad\qquad \left. + P_m^{(0,n+1/2)}(u)\left(P_{m'}^{(0,n+1/2)}(u)\right)'\right] \left(1+u\right)^{n+1/2} \mathrm{d} u \\
&\qquad\qquad  +\left. \frac{\mathbf{R}n(2n+1)}{2^{n+1}\sqrt{2}}\int_{-1}^{1} P_m^{(0,n+1/2)}(u)P_{m'}^{(0,n+1/2)}(u) \left(1+u\right)^{n-1/2} \mathrm{d} u \right],
\end{align}
where $\delta_{a,b}$ denotes the Kronecker Delta.
Note that the remaining integrals must be computed numerically due to the exponents of $(1+u)$. For two FEHFs, we obtain 
\begin{align}
&\langle N_{A,\Delta A}, N_{A',(\Delta A)'} \rangle_{\mathcal{H}_1(\mathbb{B})}\\
&=\prod_{j=1}^{3} \int_{lb_{a_j}}^{ub_{a_j}} \frac{[\Delta A_j-|a_j-A_j|][(\Delta A_j)'-|a_j-A'_j|]}{\Delta A_j (\Delta A_j)'} \left\{\begin{matrix} a_1^2,j=1,\\ 1,\ j=2,3\end{matrix} \right\} \intd a_j\\
&\qquad + \sum_{k=1}^3\int_{lb_{a_k}}^{ub_{a_k}} \frac{\sgn(a_k-A_k)}{\Delta A_k} \frac{\sgn(a_k-A_k')}{(\Delta A_k)'} \left\{\begin{matrix} a_k^2, &k=1,\\ 1,& k=2,\\ 1-a_k^2, & k=3 \end{matrix} \right\} \intd a_k\\ &\qquad \times \prod_{j=1,\ j\not=k}^3 \int_{lb_{a_j}}^{ub_{a_j}} \frac{\Delta A_j - |a_j-A_j|}{\Delta A_j} \frac{(\Delta A_j)' - |a_j-A_j'|}{(\Delta A_j)'} \left\{\begin{matrix} \frac{1}{1-a_j^2}, &j=3,k=2,\\ 1,& \text{else} \end{matrix} \right\} \intd a_j,
\end{align}
where $lb_{a_i}$ and $ub_{a_i}$ are the lower and upper bound with respect to the dimension $a_i$ of the intersection of the respective domains of the FEHFs $N_{A,\Delta A}$ and $N_{A',(\Delta A)'}$. Note that we need to determine these boundaries as well as the critical points in between them.

At last, we consider the mixed case of a FEHF and a polynomial. This yields
\begin{align}
&\left\langle N_{A,\Delta A}, G^{\mathrm{I}}_{m,n,j} \right\rangle_{\mathcal{H}_1(\mathbb{B})}\\
&= 
p_{m,n}q_{n,j}\int_{lb_r}^{ub_r} \frac{\Delta R-|r-R|}{\Delta R} P_m^{(0,n+1/2)}(I(r))\left(\frac{r}{\mathbf{R}}\right)^nr^2 \mathrm{d} r\\
&\quad\times \int_{lb_\lon}^{ub_\lon} \frac{\Delta \Phi-|\varphi-\Phi|}{\Delta \Phi} \mathrm{Trig}(j\varphi) \mathrm{d} \varphi \int_{lb_t}^{ub_t} \frac{\Delta T-|t-T|}{\Delta T} P_{n,|j|}(t) \mathrm{d} t\\
& - p_{m,n}q_{n,j}
\sum_{k=1}^3 
\int_{lb_{a_k}}^{ub_{a_k}}
\frac{\sgn(a_k-A_k)}{\Delta A_k}
\left\{\begin{matrix} 
\left(P_m^{(0,n+1/2)}\left(I(a_k)\right)\right)'I'(a_k) \left(\frac{a_k}{\mathbf{R}}\right)^{n} a_k^2 &\\ + 
P_m^{(0,n+1/2)}\left(I(a_k)\right) \left(\frac{a_k}{\mathbf{R}}\right)^{n} n a_k,  &k=1\\ 
j\mathrm{Trig}(-j a_k),&k=2\\ 
(1-a_k^2)P_{n,|j|}'(a_k),&k=3 
\end{matrix}\right\}
\intd a_k\\
&\quad \times
\prod_{i=1,i\not=k}^3
\int_{lb_{a_i}}^{ub_{a_i}}
\frac{\Delta A_i-|a_i-A_i|}{\Delta A_i}
\left\{\begin{matrix} 
P_m^{(0,n+1/2)}\left(I(a_i)\right) \left(\frac{a_i}{\mathbf{R}}\right)^{n}, &i=1\\ 
\mathrm{Trig}(j a_i), &i=2\\ 
P_{n,|j|}(a_i),&i=3,k=1\\
\frac{1}{1-a_i^2} P_{n,|j|}(a_i),&i=3,k=2
\end{matrix}\right\}
\intd a_i
\end{align}
where $lb_r = lb_{a_1} = \max(R_{\mathrm{min}},R-\Delta R),\ ub_r = ub_{a_1} = \min(R_{\mathrm{max}},R+\Delta R),\ lb_\lon = lb_{a_2} = \Phi-\Delta \Phi,\ ub_\lon = ub_{a_2} = \Phi+\Delta \Phi,\ lb_t = lb_{a_3} = \max(T_{\mathrm{min}},T-\Delta T)$ and $ub_t = ub_{a_3} = \min(T_{\mathrm{max}},T+\Delta T)$.
Note that we use the same notation as in \cref{eq:nablaGI}. The longitudinal integrals can be derived analytically (see \cref{ssect:app:H1IPs}). The radial and the latitudinal integrals must be computed numerically. Note that, in our experience, it is difficult but critical for our approach to determine the most efficient implementation for the latitudinal integrals.

\section{The (Learning) Inverse Problem Matching Pursuits}
\label{sect:ipmps}
Next, we introduce our suggested algorithms for the approximation of ill-posed inverse problems: the (Learning) Inverse Problem Matching Pursuits ((L)IPMPs). The IPMPs include: the Regularized Functional Matching Pursuit (RFMP) and the Regularized Orthogonal Functional Matching Pursuit (ROFMP). The LIPMPs are the respective counterparts that include a learning add-on: the LRFMP and the LROFMP. Note that there also exists the Regularized Weak Functional Matching Pursuit (RWFMP) whose idea can be included in the LRFMP and the LROFMP quite naturally as we explain below. To the best of our knowledge, the Geomathematics Group Siegen is among the first to adapt these algorithms for inverse problems. For more details and applications, see 
\cite{Berkeletal2011,Fischer2011,Fischeretal2012,Fischeretal2013-1,Fischeretal2013-2,Guttingetal2017,Kontak2018,Kontaketal2018-2,Kontaketal2018-1,Michel2015-2,Micheletal2017-1,Micheletal2018-1,Micheletal2014,Micheletal2016-1,Leweke2018,Prakashetal2020,Schneider2020,Schneideretal2022,Telschow2014,Telschowetal2018}. 
We concentrate here on the LRFMP because it appears to be more suitable for travel time tomography due to its efficiency. For this method, we introduce here the characteristics relevant for understanding and those that are newly adjusted to the particular problem at hand. Note, however, that the transfer of the latter to the LROFMP is straightforward and an implementation of it is included in the corresponding source code. We also direct the reader to the notations given in \cref{sect:problem}.

\subsection{The (Learning) Regularized Functional Matching Pursuit}
\label{ssect:ipmps}
The RFMPs tackles inverse problems, such as the discretized, linearized travel time tomography problem $T_\daleth \delta S = \delta \psi_\daleth$, which are often ill-posed by nature. The inevitable regularization for such kind of problems is implemented as a Tikhonov--Phillips regularization in these methods. This is an established and well-performing choice for many ill-posed inverse problems, see e.g. \cite{Engletal1996,Hofmann1999,Kirsch1996,Louis1989,Rieder2003}. An approximation $f^\ast$ of the solution $f$ is then found as the minimizer of the Tikhonov--Phillips functional. In particular, with the $\Hs{1}$ Sobolev space introduced above, in the RFMP, we aim to determine $f^\ast$ such that 
\begin{align}
\left\| \delta \psi_\daleth - \T_\daleth f^\ast \right\|^2_{\real^\ell} + \lambda\left\|f^\ast\right\|^2_{\Hs{1}},\qquad \lambda >0,
\label{eq:TF}
\end{align}
is minimized. The first term is usually called the data misfit while the latter is the penalty term. Note that the corresponding Tikhonov--Phillips functional consists formally of only one penalty term instead of two (for smoothing and damping) as commonly used in seismology, see e.g. \cite{Charletyetal2013,Hosseinietal2020}. We choose the $\Hs{1}$-Sobolev space for regularization here because we 
considered FEHFs for the approximation and they are tightly connected to (classical) Sobolev spaces. As we saw in \cref{def:H1IP}, the respective inner product is a sum of the $\Lp{2}$- and the $\lp{2}$-inner product due to which we re-enact the trade-off between smoothing and damping. Thus, this difference is only minor. \\
However, there is a more important difference between  \cref{eq:TF} and the common seismological approach and it lies within the data misfit: due to uncertainties within the data, seismologists usually consider the (reduced) $\chi^2$ in the computation process as well as for model selection (via the L-curve) instead of the pure residual $\delta \psi_\daleth - \T_\daleth f^\ast$, see e.g. \cite{Hosseinietal2020}. Mathematically, they are connected as follows: 
\begin{align}	
\chi_{\mathrm{red}}^2 
&= \frac{1}{\ell} \sum_{i=1}^\ell \left( \frac{\left(\delta \psi_\daleth\right)_i - \T^i_\daleth f^\ast}{\sigma_i}\right)^2
= \frac{1}{\ell} \left\| \frac{\delta \psi_\daleth - \T_\daleth f^\ast}{\sigma}\right\|_{\real^\ell}^2
\intertext{with}
\frac{\delta \psi_\daleth - \T_\daleth f^\ast}{\sigma} 
&\coloneqq \left( \frac{\left(\delta \psi_\daleth\right)_i - \T^i_\daleth f^\ast}{\sigma_i} \right)_{i=1,...,\ell},
\end{align} 
where $\sigma_i \in \real_+$ denotes the (known) uncertainty with respect to the $i$-th measurement. As the uncertainty within the data cannot be circumvented, we, therefore, adjust the Tikhonov--Phillips functional considered in the RFMP for the case of travel time tomography as well. Therefore, in the sequel, we consider the noise-cognizant Tikhonov-Phillips functional 
\begin{align}
\left\| \frac{\delta \psi_\daleth - \T_\daleth f^\ast}{\sigma} \right\|^2_{\real^\ell} + \lambda\left\|f^\ast\right\|^2_{\Hs{1}},\qquad \lambda >0,
\label{eq:ncTF}
\end{align}
which shall be minimized in the RFMP. The minimizer $f^*$ is then obtained iteratively as a linear combination of weighted dictionary elements $d_n \in \dic$. Let $N$ denote the current (or final) iteration. Then we have
\begin{align}
f_N &= f_0 + \sum_{n=1}^N \alpha_n d_n
\label{eq:fn}
\end{align}
in the case of the RFMP.
Here, $f_0$ denotes a first approximation which is often the zero approximation in practice. As the dictionary is made of global polynomials and local FEHFs, also the approximation consists -- in all probability -- of both types of trial functions. \\
The noise-cognizant Tikhonov--Phillips functional of the $N$-th step transfers then to 
\begin{align}
(\alpha, d) &\mapsto \left\| \frac{R^N - \alpha\T_\daleth d}{\sigma} \right\|^2_{\real^\ell} + \lambda\left\|f_N+\alpha d\right\|^2_{\Hs{1}},\qquad \lambda >0,
\label{eq:TFNO}\\
&R^{N+1} \coloneqq R^N - \alpha_{N+1}\T_\daleth d_{N+1} = \delta \psi - \T_\daleth f_{N+1}
\label{def:RN}
\end{align}
for the RFMP and with
$R^0 = \delta \psi_\daleth-\T_\daleth f_0$ which yields $R^0 = y$ if $f_0 \equiv 0$. 
The main question is how to choose the dictionary element $d_{N+1} \in \dic$ and the corresponding weight $\alpha_{N+1} \in \real$ 
such that the corresponding Tikhonov--Phillips functional is minimized. 
Similarly as in the literature on the (L)IPMPs, we can exchange the minimization of the noise-cognizant Tikhonov--Phillips functional by an equivalent maximization of the objective functions, see \cref{sect:app:OF_IPMPs} for details:
\begin{align}
\RFMP(d;N) &\coloneqq \frac{\left( \left\langle \frac{R^N}{\sigma}, \frac{\T_\daleth d}{\sigma}\right\rangle_{\real^\ell} - \lambda\left\langle f_N,d \right\rangle_{\Hs{1}} \right)^2}{\left\|\frac{\T_\daleth d}{\sigma}\right\|_{\real^\ell}^2 + \lambda\left\|d\right\|^2_{\Hs{1}}}
\eqqcolon \frac{\left(A_N(d)\right)^2}{B_N(d)}.
\label{def:RFMP(d;N)}
\end{align}
The weights are then easily obtained via 
\begin{align}
\alpha_{N+1} &\coloneqq \frac{A_N(d_{N+1})}{B_N(d_{N+1})}.
\label{def:alphan}
\end{align}

In practice, the IPMPs need termination and model selection criteria. As implemented for the Tikhonov--Phillips functional, we adjust them here as well. Usually, in seismological experiments, we would strive to let $\chi^2_{\mathrm{red}}$ reach 1. This should yield the best trade-off between data (mis)fit, accuracy and smoothing. For the contrived data we use here, however, we have the uncertainty is assumed to be $\sigma\equiv 1\,$s. This enables us to consider the relative data error $\|R^N\|_{\real^\ell}/\|R^0\|_{\real^\ell}$ (as usually done in an IPMP) instead. In our experiments, we additionally perturb the delay vector with simulated noise. This again allows us to terminate the algorithm if the relative data error falls below this noise level. To avoid endless iterations for inappropriate parameters, we also set a maximum number of iterations. Among those models $f_N$ 
obtained for diverse regularization parameter $\lambda$, we select the one (i.e. the regularization parameter) which yields the lowest relative root mean square error
\begin{align}
\left(\frac{\sum_{i=1}^\kappa \left(f\left(x^i\right)-f_N\left(x^i\right)\right)^2}{\sum_{i=1}^\kappa f\left(x^i\right)^2}\right)^{1/2}
\label{def:rrmse}
\end{align}
for $\kappa\in\nat$, where $f$ is the (exact) solution which we use for our test and is given on the points $x^i$.

The IPMPs use by definition a finite dictionary $\dic^{\mathrm{fin}} \subset \dic = \dic^{\mathrm{Inf}}$. In this case, the maximization of \cref{def:RFMP(d;N)} 
can be done by pairwise comparisons. However, this means that $\dic^{\mathrm{fin}}$ must be chosen a-priori either automatically or manually. The latter cannot be recommended in the case of the travel time tomography because a) we are inexperienced regarding which trial functions are needed as this is a novel application for the methods; b) the size of $\dic^{\mathrm{fin}}$ grows tremendously due to the six characteristic parameters of the FEHFs; and c) a possible bias by the choice of $\dic^{\mathrm{fin}}$ cannot be quantified. 
In search of an automation of the a-priori choice, the LIPMPs were developed, see \cite{Micheletal2018-1,Schneider2020,Schneideretal2022}. They follow the same routine as the IPMPs but include an additional learning add-on. Then, the a-priori dictionary choice is negligible. Though they do produce a learnt dictionary which can be used as an automatically chosen one in the IPMPs, they also proved to be useful as standalone approximation algorithms. Thus, in our experiment, we include the learning add-on, that is, we use the LRFMP.

\subsection{The learning add-on}
\label{ssect:lipmps}
The idea is to allow all possible trial functions and, thus, use the infinite dictionary $\dic$. As we saw before, $\dic$ includes infinitely many trial functions of  those types with continuous characteristic parameters, i.e. here the FEHFs. If the characteristic parameters are discrete as here with the polynomials, we still allow only a finite set. In each iteration, we determine optimized dictionary elements (or candidates) for each type of trial functions separately. Together they form again a very small, finite dictionary of candidates from which we obtain the overall most suitable function. Thus, the learning add-on is the determination of the finite dictionary of candidates in each iteration.\\
Recall that we allow polynomials up to a maximally possible radial and angular degree. The global trial functions in the dictionary are usually chosen to reconstruct global trends and, thus, high degrees and orders are often negligible. Hence, it suffices to consider some maximum radial and angular degree. Theoretically, these maximally possible degrees can be chosen to be very high such that the methods can learn a maximally needed radial and angular degree, see \cite{Schneider2020,Schneideretal2022}. It remains to be seen whether this can be done in practice for travel time tomography due to efficiency reasons (the curse of dimensionality occurs here: the set of all (Cartesian) polynomials $G_{m,n,j}^{\mathrm{I}}$ with degree $\leq N$ has a size of $\mathcal{O}(N^3)$, see \cite{Michel1999}). Note that, due to the finiteness, we can still obtain the maximizer of \cref{def:RFMP(d;N)} 
by pairwise comparisons. Further, we can still use the common preprocessing routine of the IPMPs for the polynomials and improve efficiency in this way. Practically, this means we define a finite starting dictionary $\dic^{s}$ which contains at least the chosen polynomials from $\dic$. \\
In the case of the FEHFs, we use the truly infinite set of possible trial functions $[\mathcal{A}]_{\mathrm{FEHF}}$. The main challenge here is to maximize \cref{def:RFMP(d;N)} 
among all possible FEHFs. As a matter of fact, this is a non-linear constraint optimization problem with the objective function 
	$	\RFMP(N_{A,\Delta A};N) \rightarrow \max !$
For maximizing $\RFMP(N_{A,\Delta A};N)$, recall the corresponding constraints with respect to $A$ and $\Delta A$ given in \cref{def:N-A-DA-constraints}. We can solve this maximization with any kind of established optimization routine. For instance, methods from the NLOpt library, see \cite{NLOpt2019}, can be utilized. Note, however, that a gradient-based approach cannot be used here because it would necessitate the computation of the gradient of $\RFMP(N_{A,\Delta A};N)$ with respect to $(A,\Delta A)$. Unfortunately, this is not possible for practical purposes due to the absolute value in the definition of the function, see \cref{def:N-A-DA}. Further, from experience, we propose a 2-step-optimization procedure. First, we use a global method. Then we refine this solution by using a local counterpart starting from the former solution. This also enables us to soften the accuracy of the global optimization technique (i.e. its termination criteria) which decreases the runtime. Note that softening the termination criteria of the optimization is similar to including a weakness parameter as in the RWFMP. Also note that both solutions are inserted into the dictionary of candidates. Moreover, if the global method needs a starting point, we should insert a few FEHFs in the starting dictionary for this reason as well. Further note that this starting solution can also be included in the dictionary of candidates but should generally not be chosen, i.e. this starting solution should not have a major influence on the learnt approximation in general.

\subsection{Additional divide-and-conquer strategy for travel time tomography}

We experienced that the optimization within the learning add-on is slowed down -- amongst others -- by the use of many rays at once. Thus, we considered how to sensibly use fewer rays at once while taking a numerically significant number of rays into account. This yields an additional divide-and-conquer strategy for our challenge at hand which proved to be helpful in practice. The main idea is to start with a low number of rays. In our experiments, we started with 1000 rays. If the relative data error falls below a certain threshold such as $50\%$ then we add the next package of 1000 rays to our consideration. We consider then the first 2000 rays in our algorithm with the exception of optimization of the FEHFs where we consider only the latest package of 1000 rays. 

We are aware that this is a quite manual approach with certain seemingly arbitrary parameters. In the meantime, we also considered other aspects that concerned the efficiency such that it seems now to be possible to increase, for instance, the size of each ray package. This may form the basis of future research.

\section{Numerical implementation and tests}
\label{sect:numerics}
In this section, we show a numerical proof of concept for our TTIPMP code. We first introduce our experiment setting for reproducibility and afterwards our numerical results. Note that the corresponding code is available at \textcolor{blue}{\url{https://doi.org/10.5281/zenodo.8227888}} under licence the CC-BY-NC-SA 3.0 DE.

\subsection{Chosen Earthquake data and specific experiment setting}
\label{ssect:numerics:genset}

\begin{figure}[htbp]
\centering
\includegraphics[width=.32\textwidth]{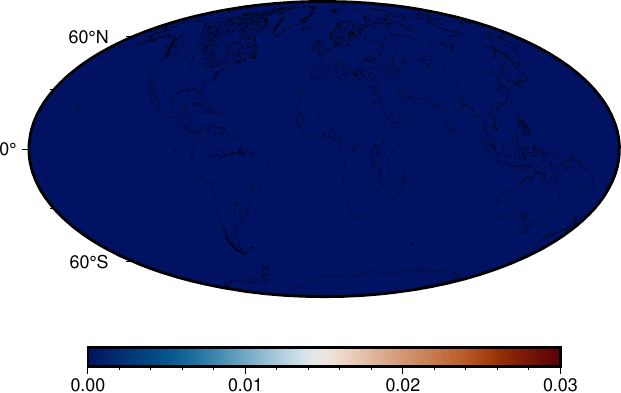}
\includegraphics[width=.32\textwidth]{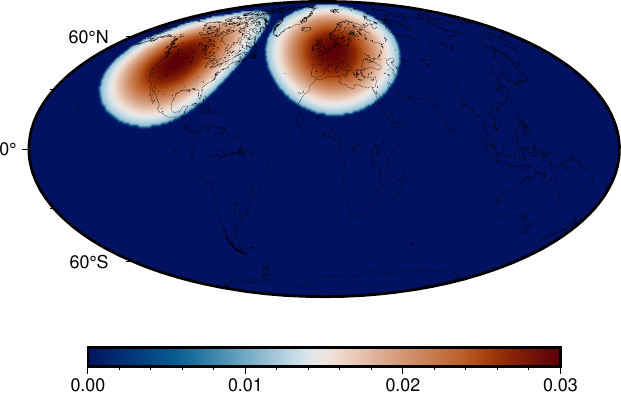}
\includegraphics[width=.32\textwidth]{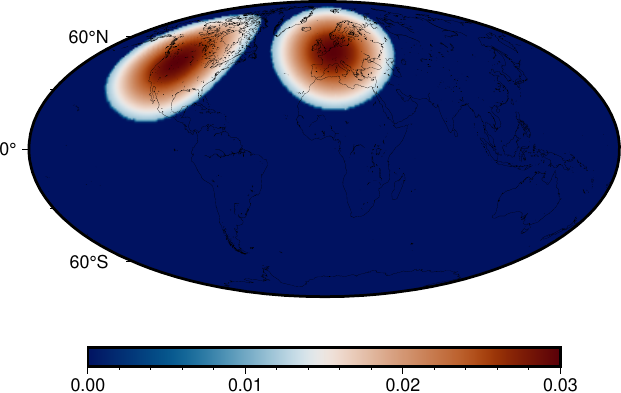}
\includegraphics[width=.32\textwidth]{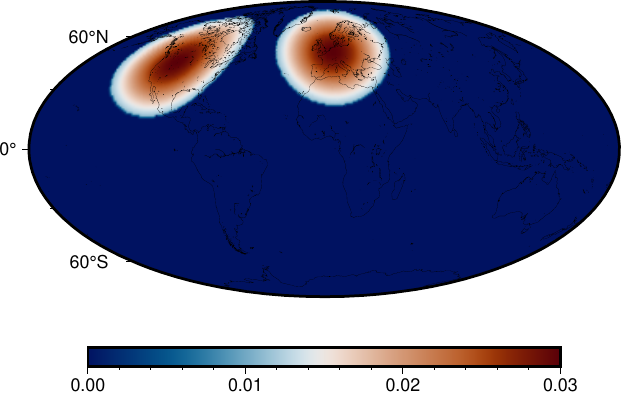}
\includegraphics[width=.32\textwidth]{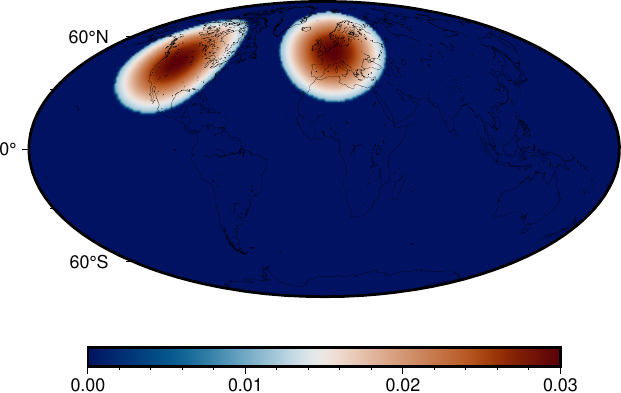}
\includegraphics[width=.32\textwidth]{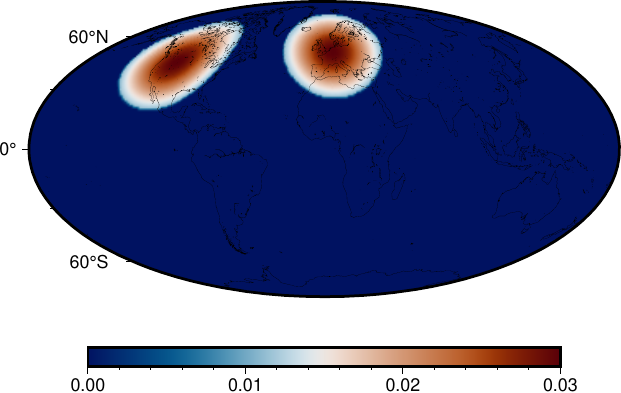}
\includegraphics[width=.32\textwidth]{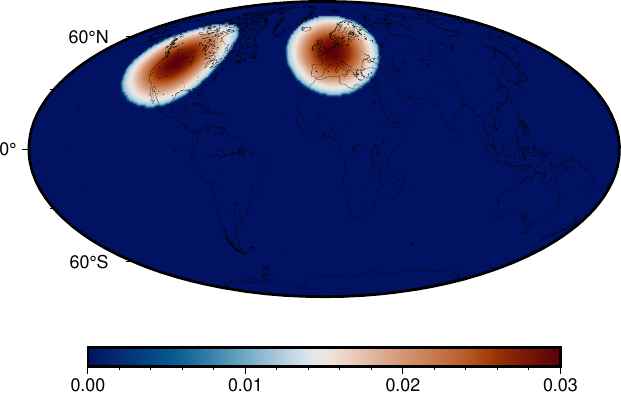}
\includegraphics[width=.32\textwidth]{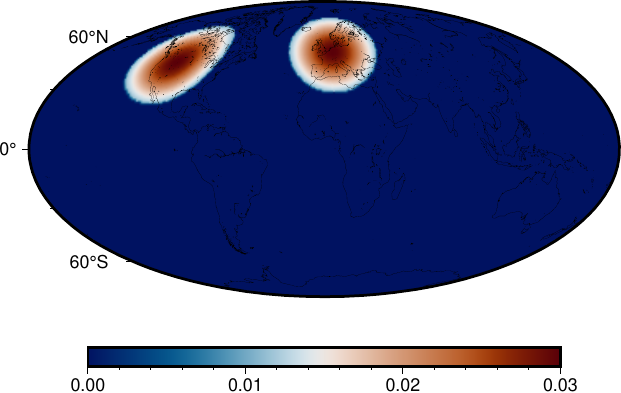}
\includegraphics[width=.32\textwidth]{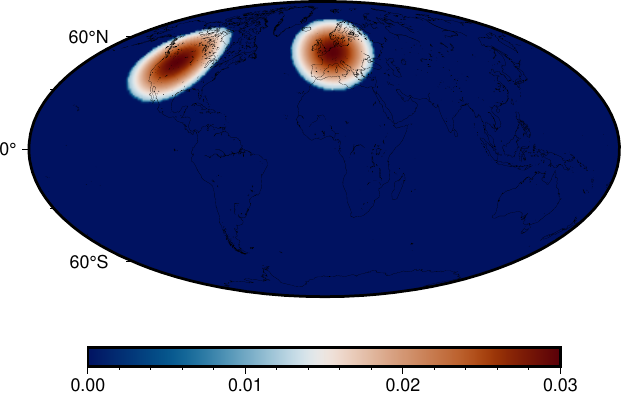}
\includegraphics[width=.32\textwidth]{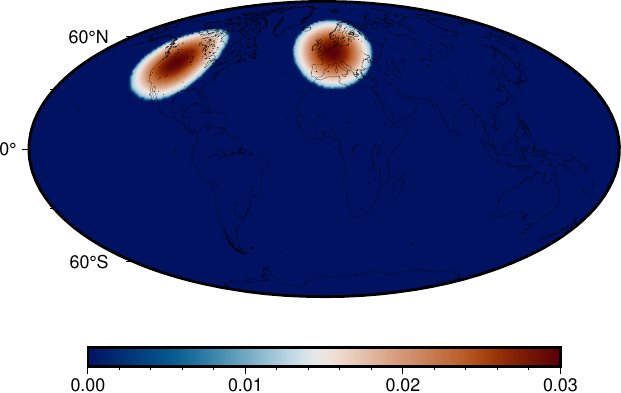}
\includegraphics[width=.32\textwidth]{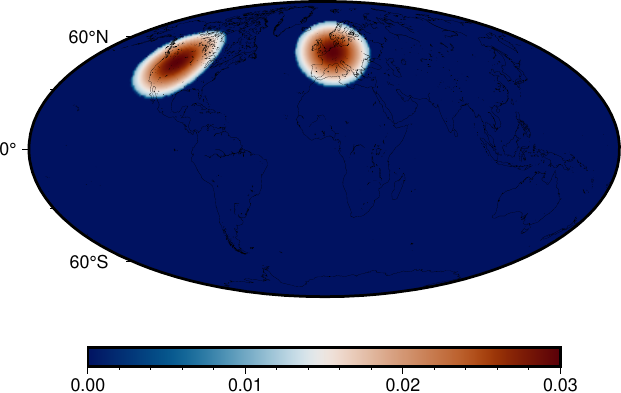}
\includegraphics[width=.32\textwidth]{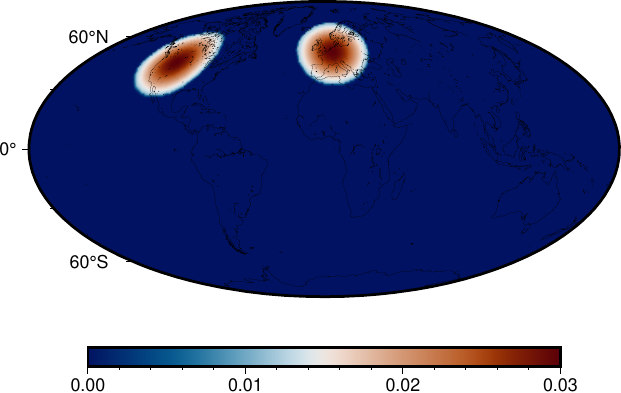}
\caption{Resolution test input model for P-velocity consisting of two plumes below the Volcanic Eifel and Yellowstone volcano. We show the depth slices (left, middle and right) at the following radial distances to the centre of the Earth in km: 3193.1, 3482.0 and 3770.9 (first row), 4059.8, 4348.7 and 4637.6 (second row), 4926.5, 5215.4 and 5504.3 (third row), 5793.2, 6082.1 and 6371.0 (last row).}
\label{fig:resolutiontests}
\end{figure}

\begin{figure}
\centering
\begin{subfigure}{.24\textwidth}
\includegraphics[width=\textwidth]{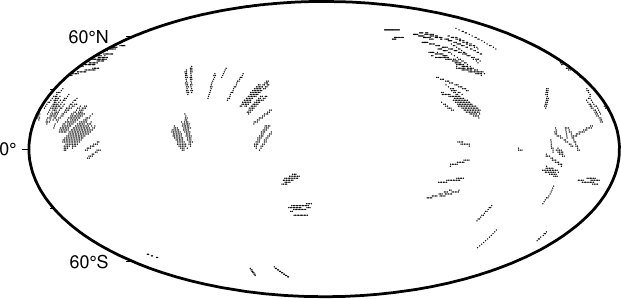}
\caption{3571-3771}
\end{subfigure}
\begin{subfigure}{.24\textwidth}
\includegraphics[width=\textwidth]{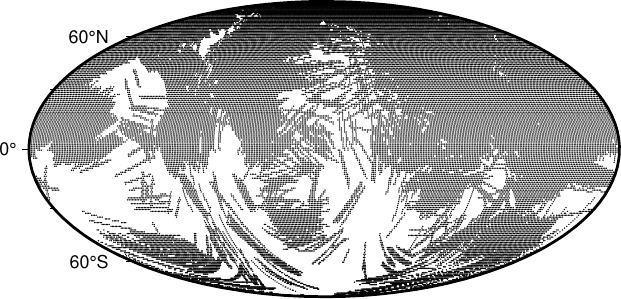}
\caption{3771-3971}
\end{subfigure}
\begin{subfigure}{.24\textwidth}
\includegraphics[width=\textwidth]{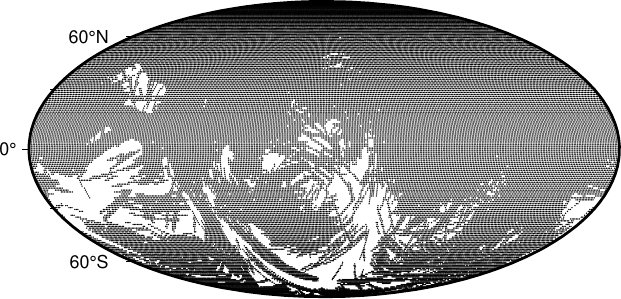}
\caption{3971-4171}
\end{subfigure}
\begin{subfigure}{.24\textwidth}
\includegraphics[width=\textwidth]{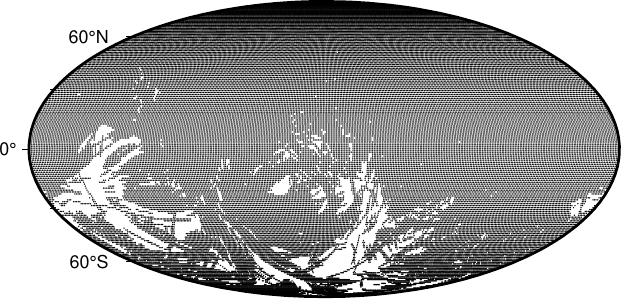}
\caption{4171-4371}
\end{subfigure}\\[\baselineskip]
\begin{subfigure}{.24\textwidth}
\includegraphics[width=\textwidth]{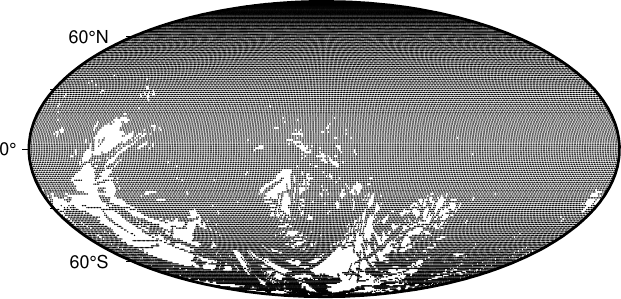}
\caption{4371-4571}
\end{subfigure}
\begin{subfigure}{.24\textwidth}
\includegraphics[width=\textwidth]{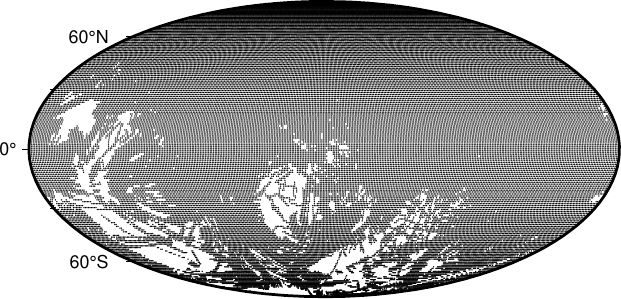}
\caption{4571-4771}
\end{subfigure}
\begin{subfigure}{.24\textwidth}
\includegraphics[width=\textwidth]{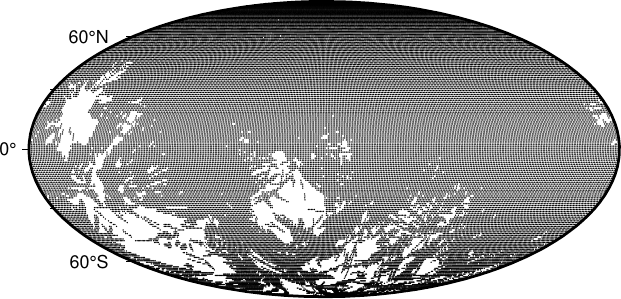}
\caption{4771-4971}
\end{subfigure}
\begin{subfigure}{.24\textwidth}
\includegraphics[width=\textwidth]{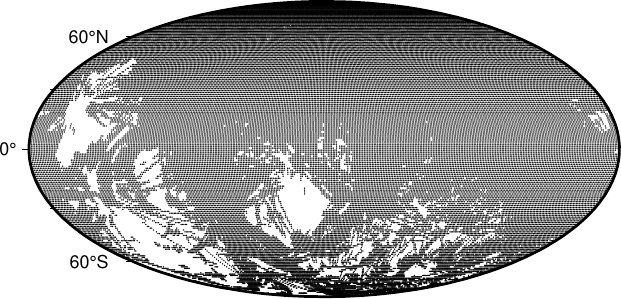}
\caption{4971-5171}
\end{subfigure}\\[\baselineskip]
\begin{subfigure}{.24\textwidth}
\includegraphics[width=\textwidth]{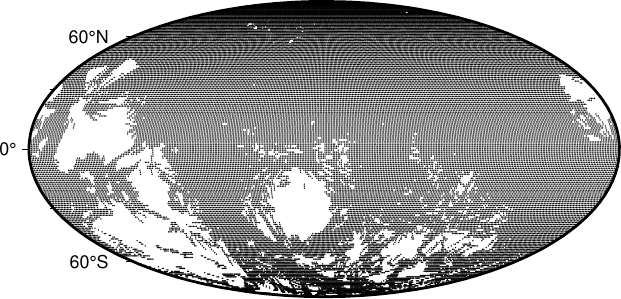}
\caption{5171-5371}
\end{subfigure}
\begin{subfigure}{.24\textwidth}
\includegraphics[width=\textwidth]{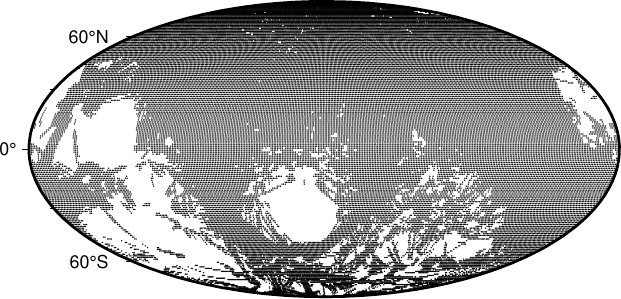}
\caption{5371-5571}
\end{subfigure}
\begin{subfigure}{.24\textwidth}
\includegraphics[width=\textwidth]{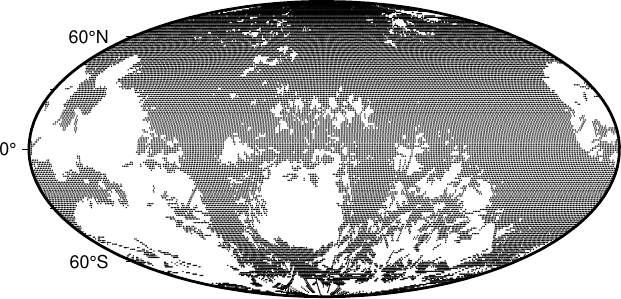}
\caption{5571-5771}
\end{subfigure}
\begin{subfigure}{.24\textwidth}
\includegraphics[width=\textwidth]{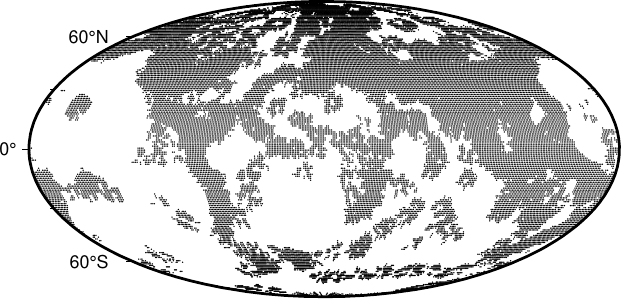}
\caption{5771-5971}
\end{subfigure}\\[\baselineskip]
\begin{subfigure}{.24\textwidth}
\includegraphics[width=\textwidth]{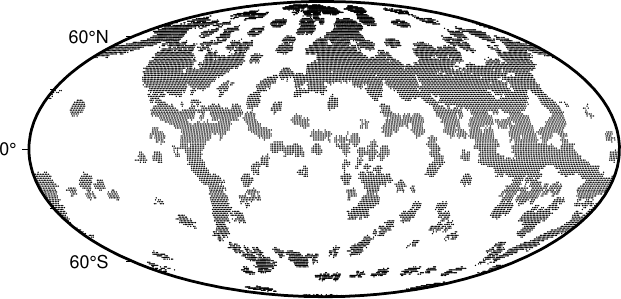}
\caption{5971-6171}
\end{subfigure}
\begin{subfigure}{.24\textwidth}
\includegraphics[width=\textwidth]{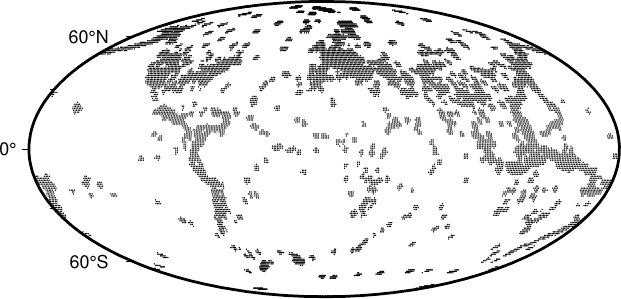}
\caption{6171-6371}
\end{subfigure}
\caption{Distribution of ISC-EHB meta-data rays in the given depth ranges. Depth in km as the distance to the Earth's centre.}
\label{fig:distrrays}
\end{figure}

We use the unit ball for our practical computations and show here the approximation of synthetic input structures  modelled on mantle plumes, i.e. seismically slow, vertical, columnar features that extend from the core-mantle boundary to the surface, as a model for the deviation of the slowness in the interior of the Earth. That is, we show here a resolution test. The plumes are constructed to extend from the core-mantle boundary up to the Earth's surface. They lie below the Volcanic Eifel and the Yellowstone volcano. Note that thin, vertically continuous plume conduits are among the structures suspected to occur in the real Earth. For a visualization of this model, see \cref{fig:resolutiontests}. \\
This enables us to compute a relative root mean square error as given in \cref{def:rrmse} where we use a grid of $\kappa = 12 \cdot 65341 = 784\,092$ data points. This grid is given as 12 layers of an equi-angular grid, commonly also called a Driscoll-Healy grid, see e.\,g.\, \cite{DriscollHealy1994,Michel2013}, of 361 $\cdot$ 181 grid points.
The data is perturbed with $5\%$ Gaussian noise, such that we have perturbed data $y^\delta$ given by 
\begin{align}
y^\delta_i = y_i \cdot \left( 1 + 0.05 \cdot \varepsilon_i\right)
\label{eq:noise}
\end{align}
for the unperturbed data $y_i = \left( \T f \right) \left(\sigma^i\eta^i\right)$ and a unit normally distributed random number $\varepsilon_i$. \\
As we aim for a proof of concept, it suffices to consider a ray theoretical setting and the ISC-EHB seismic meta-data, see \cite{ISC-EHB,Westonetal2018}. Note that these data were among others also used in \cite{Hosseinietal2020}. In particular, we are using the ISC-EHB meta-data from 1998 to 2016, respectively, but only the P waves. Note, however, due to the divide-and-conquer strategy, we might not use all rays from these years in practice. At most we have used 318\,542 rays. According to our strategy, we start with rays in 2016 and add further rays by moving back in time. \\
We correct these data with all necessary seismological corrections. The corrections are done in line with the tomography workflows exhibited in \cite{Hosseinietal2020,Mohammadzaherietal2021,Tsekhmistrenkoetal2021}. Moreover, with these meta-data, we have a constant uncertainty $\sigma_i = 1\, \mathrm{s}, i=1,...,\ell$. In general, we cannot expect the rays to illuminate the Earth evenly. That is, we have to take into account that all of our meta-data -- the receiving seismological stations as well as the rays themselves -- will be poorly-distributed to a certain extent. In particular, the path coverage is sparse in shallow depths under the largely uninstrumented oceans as well as in the deepest parts of the mantle since we exclude core-diffracted body-wave paths. For a visualization of the ISC-EHB meta-data rays, see \cref{fig:distrrays}.  \\
We solve the latitudinal integrals in the inner products of a polynomial and an FEHF with a Gauß-Legendre quadrature rule of $10^{6}$ points and use an adaptive Gauß-Kronrod quadrature rule with an integration error of $10^{-4}$ anywhere else (i.e. for the DSPO as well as for other integrals in the inner products). In particular, the latter is necessary for efficiency reasons. We set $\eps_R = \eps_\Phi = \eps_T = 10^{-2}$ and $\rho=3482/6371$. Recall that this sets the lower bound for the value of $R$ of an FEHF to the core-mantle boundary in our setting. We use the GN\_DIRECT\_L and the LN\_SBPLX algorithm for the global and local optimization in the learning add-on. We terminate the global algorithm if succeeding iterates vary less than $10^{-4}$ (i.e. xtol\_rel = $10^{-4}$) or succeeding function values vary less than $10^{0}$ (i.e. ftol\_rel = $10^{0}$). We terminate the local algorithm if successive iterates vary less than $10^{-8}$ or successive function values vary less than $10^{-4}$. In analogy to \cite{Schneider2020,Schneideretal2022}, we terminate the optimization after 10\,000 evaluations of the objective function or 600 s computation time. From experience, these termination criteria ensure that the optimization happens within a suitable time frame. Further, the loss in accuracy is negligible for our proof of concept. For more information on the termination criteria of the optimization methods, see \cite{NLOpt2019}. We terminate the LRFMP either after 300 iterations, or if $\|R^N\|_{\real^\ell}/\|y\|_{\real^\ell}$ is less than the noise level or greater than 2 or if $|\chi^2_{\mathrm{red}}-1|<10^{-8}$. \\
As a finite (starting) dictionary, we utilize 
\begin{align}
[\mathcal{G}_{5,5}]_{\mathrm{GI}} &= \left\{(m,n,j)\ \colon\ m,\ n \in \nat_0,\ m\leq 5,\ n\leq 5,\ j\in\intg_0,\ |j| \leq n\right\},\\
[\mathcal{A}]_{\mathrm{FEHF}} &=\left\{(A,\Delta A)\ \colon\ A \in D_p,\ \Delta A = \Delta D_p\right\},\\
D_p &= \left\{\frac{3482}{6371} + \frac{2889 i}{25484}\ \colon i=0,...,4\right\} \times  \left\{\frac{\pi i}{2}\ \colon i=0,...,4\right\} \\ &\qquad \times  \left\{-1+\eps_T + \frac{(1-\eps_T) i}{2}\ \colon i=0,...,4\right\},\\
\Delta D_p &= \left(  \frac{2889}{25484}, \frac{\pi}{2}, \frac{1-\eps_T}{2}\right)^\trans,\\
\dic &\coloneqq [\mathcal{G}]_{\mathrm{GI}} \cup [\mathcal{A}]_{\mathrm{FEHF}}.
\end{align}
Note that $\rho\RAD = 3482/6371$, $\RAD = 1.0$ and, thus, $\RAD-\rho\RAD = 2889/6371$ in our setting. 

\subsection{Synthetic inversion tests}
\label{ssect:numerics:exps}

\begin{figure}[htbp]
\centering
\includegraphics[width=.32\textwidth]{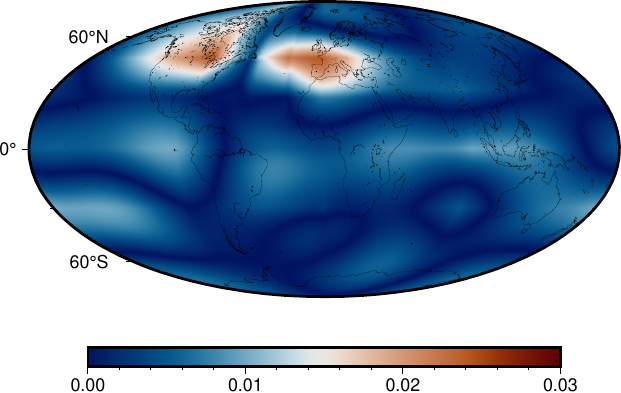}
\includegraphics[width=.32\textwidth]{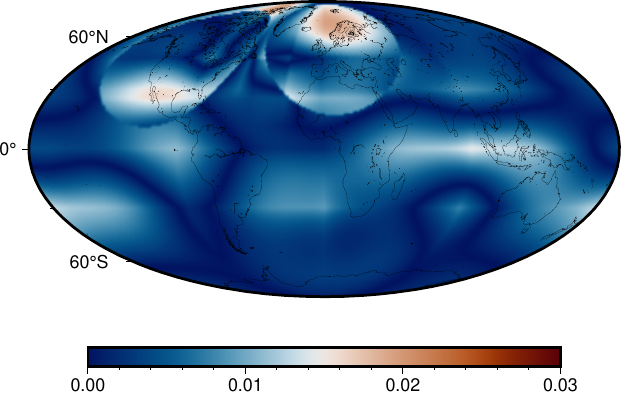}
\includegraphics[width=.32\textwidth]{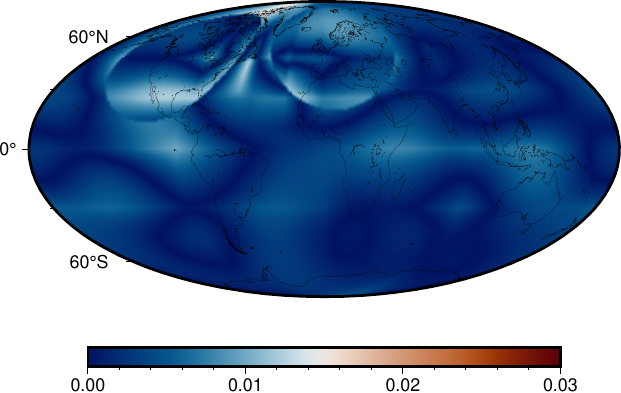}
\includegraphics[width=.32\textwidth]{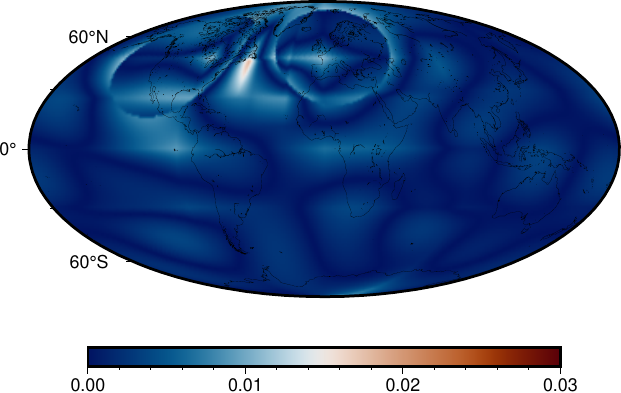}
\includegraphics[width=.32\textwidth]{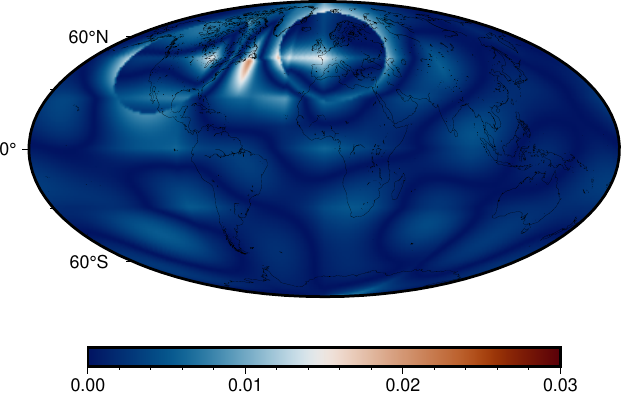}
\includegraphics[width=.32\textwidth]{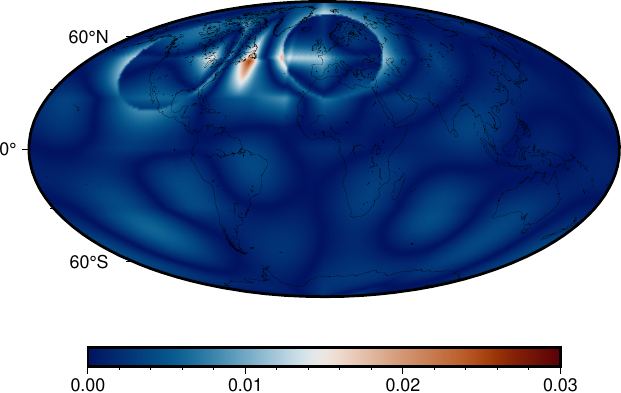}
\includegraphics[width=.32\textwidth]{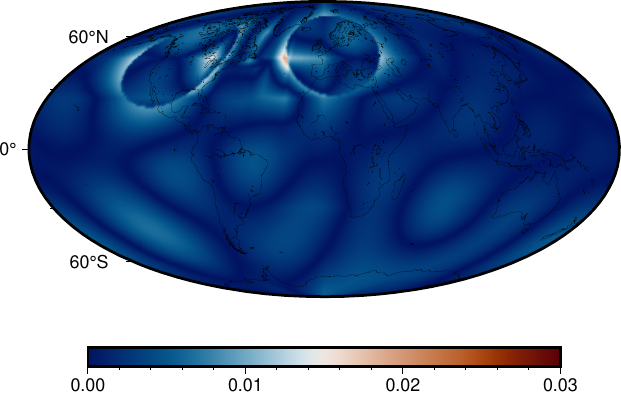}
\includegraphics[width=.32\textwidth]{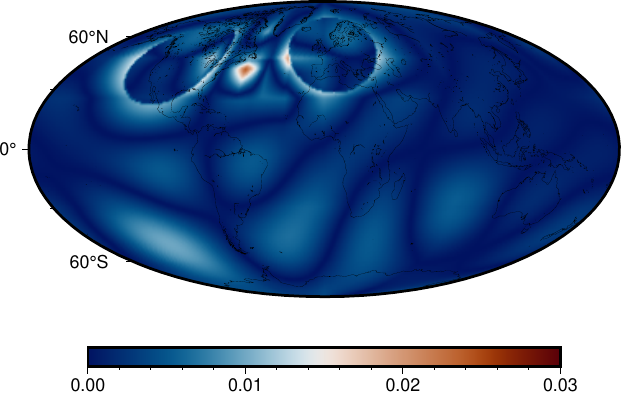}
\includegraphics[width=.32\textwidth]{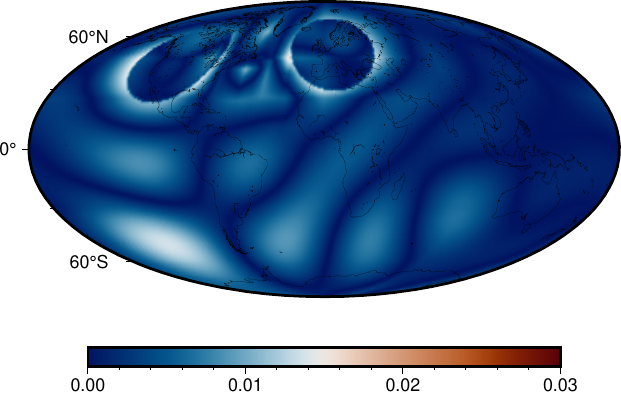}
\includegraphics[width=.32\textwidth]{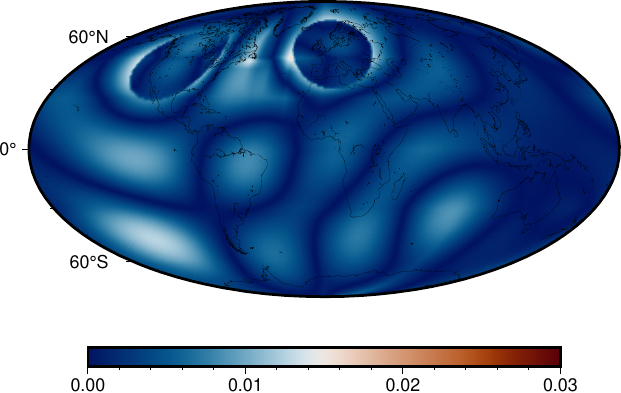}
\includegraphics[width=.32\textwidth]{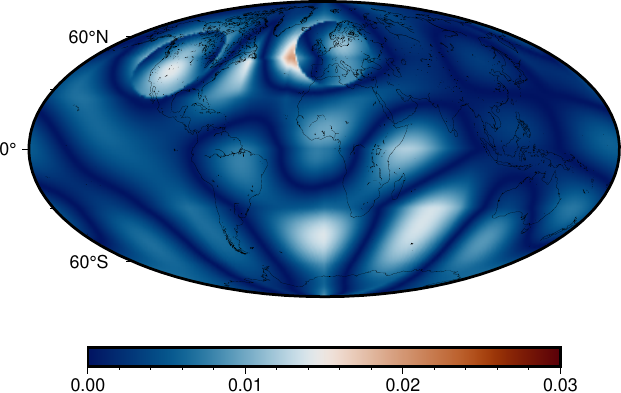}
\includegraphics[width=.32\textwidth]{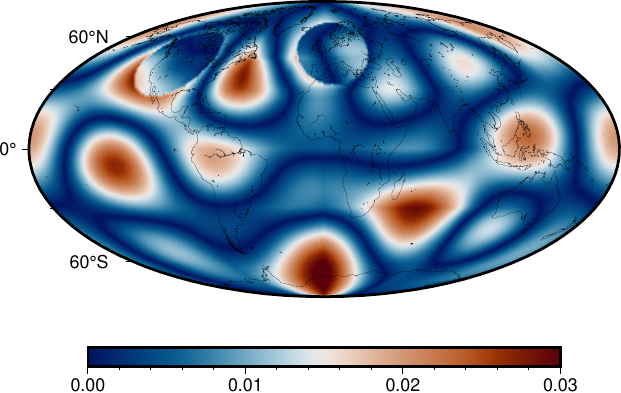}
\caption{Absolute approximation error for the approximation of the plumes model given in \cref{fig:resolutiontests}. We show the depth slices (left, middle and right) at the following radial distances to the centre of the Earth in km: 3193.1, 3482.0 and 3770.9 (first row), 4059.8, 4348.7 and 4637.6 (second row), 4926.5, 5215.4 and 5504.3 (third row), 5793.2, 6082.1 and 6371.0 (last row). The colour scales are adjusted for a better comparability.}
\label{fig:apprerr}
\end{figure}

\begin{figure}
\begin{subfigure}{.49\textwidth}
\includegraphics[width=\textwidth]{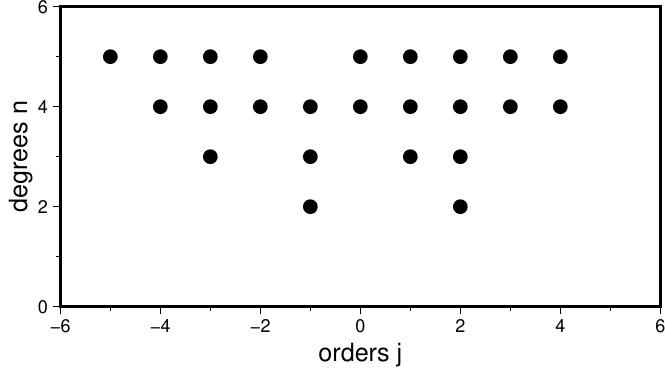}
\caption{$m=0$}
\end{subfigure}
\begin{subfigure}{.49\textwidth}
\includegraphics[width=\textwidth]{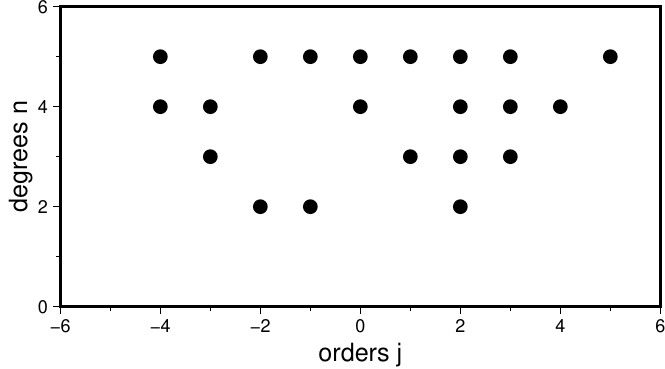}
\caption{$m=1$}
\end{subfigure}
\begin{subfigure}{.49\textwidth}
\includegraphics[width=\textwidth]{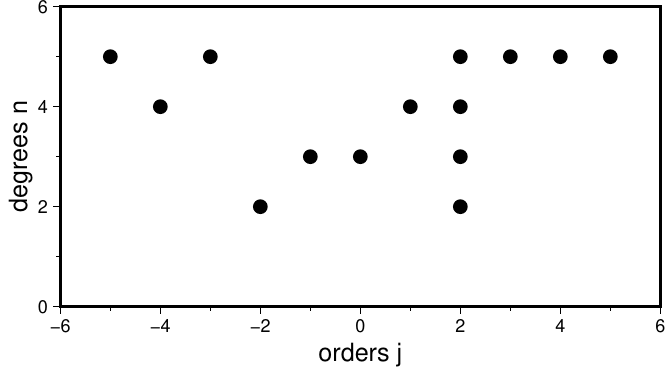}
\caption{$m=2$}
\end{subfigure}
\begin{subfigure}{.49\textwidth}
\includegraphics[width=\textwidth]{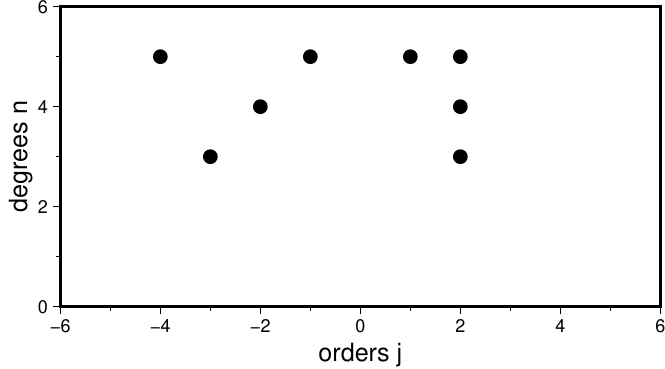}
\caption{$m=3$}
\end{subfigure}
\begin{subfigure}{.49\textwidth}
\includegraphics[width=\textwidth]{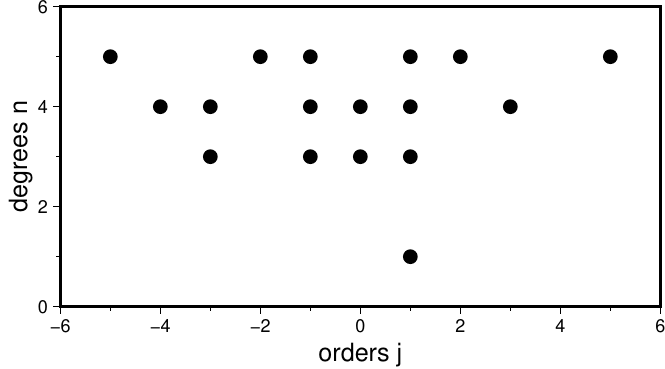}
\caption{$m=4$}
\end{subfigure}
\begin{subfigure}{.49\textwidth}
\includegraphics[width=\textwidth]{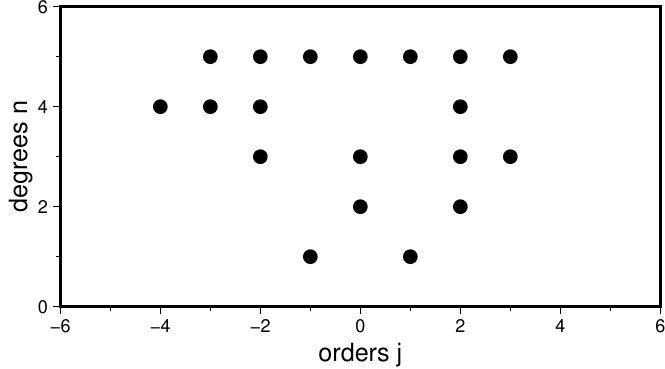}
\caption{$m=5$}
\end{subfigure}
\caption{Chosen polynomials for the given radial degrees $m$.}
\label{fig:Polys}
\end{figure}

\begin{figure}
\centering

\begin{subfigure}{\textwidth}
\includegraphics[width=.32\textwidth]{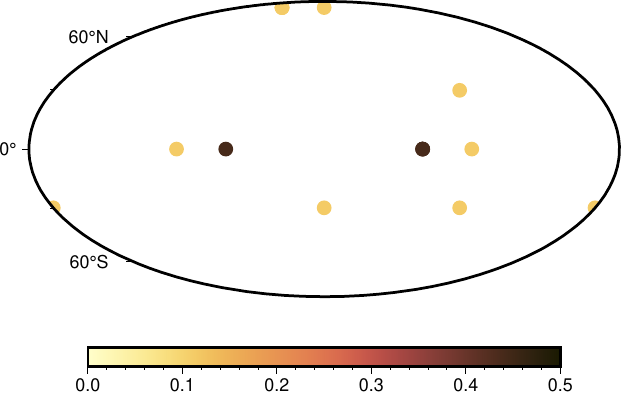}
\includegraphics[width=.32\textwidth]{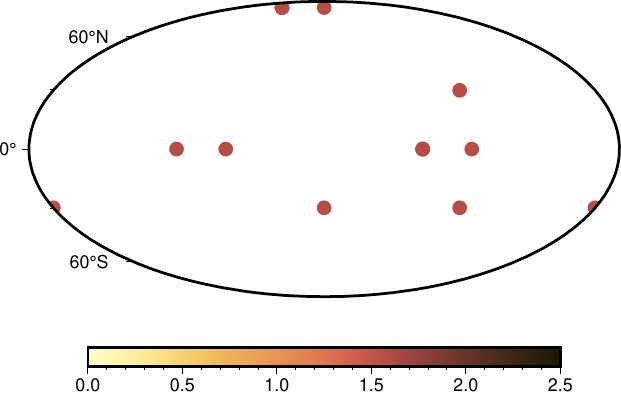}
\includegraphics[width=.32\textwidth]{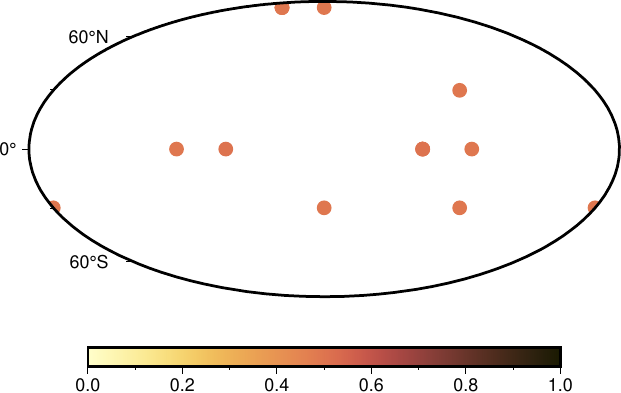}
\caption{3371 km - 3571 km}
\end{subfigure}

\begin{subfigure}{\textwidth}
\includegraphics[width=.32\textwidth]{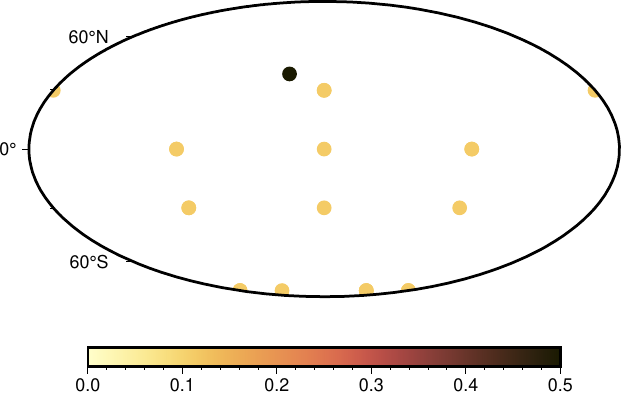}
\includegraphics[width=.32\textwidth]{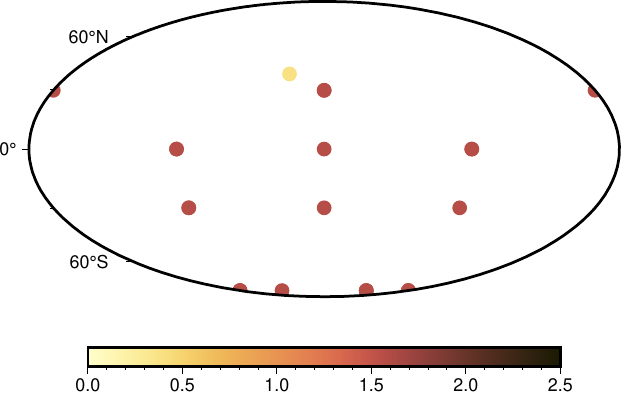}
\includegraphics[width=.32\textwidth]{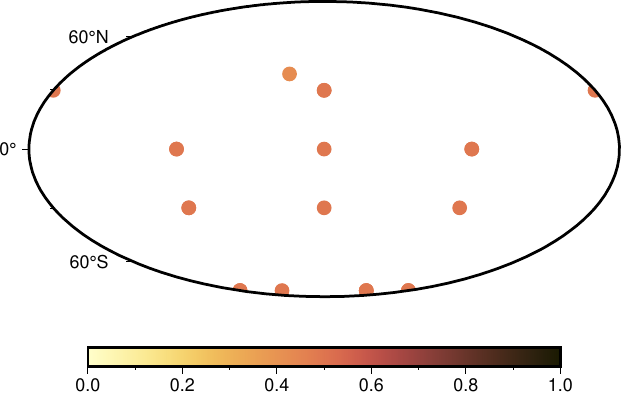}
\caption{4171 km - 4371 km}
\end{subfigure}

\begin{subfigure}{\textwidth}
\includegraphics[width=.32\textwidth]{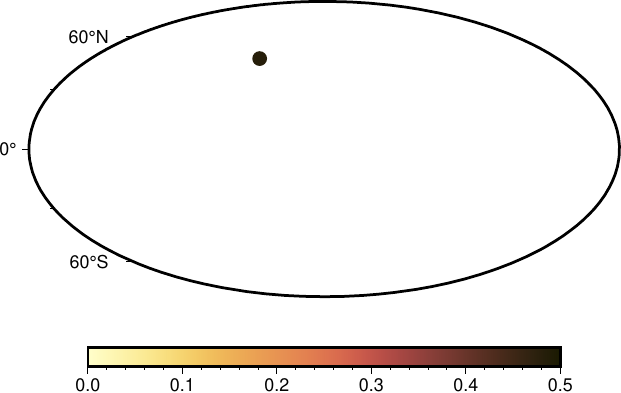}
\includegraphics[width=.32\textwidth]{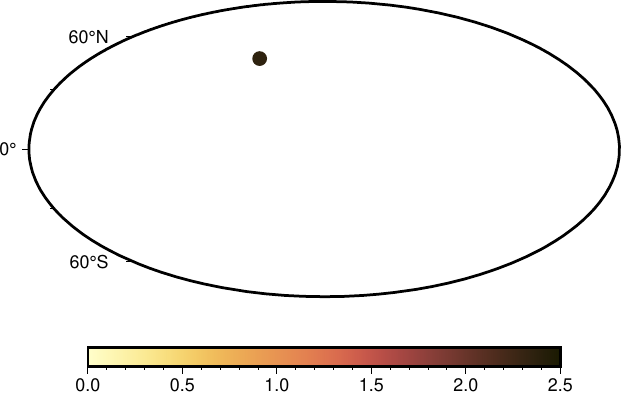}
\includegraphics[width=.32\textwidth]{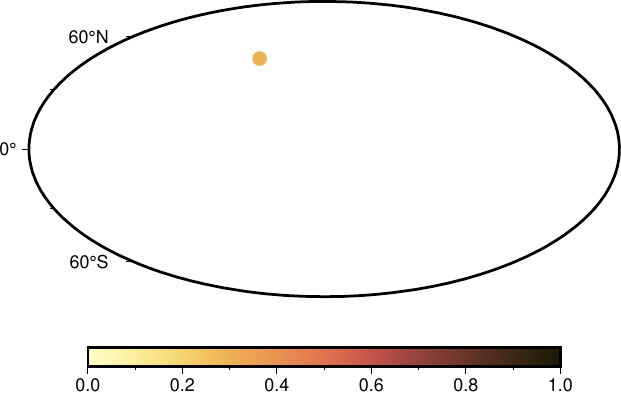}
\caption{4371 km - 4571 km}
\end{subfigure}

\begin{subfigure}{\textwidth}
\includegraphics[width=.32\textwidth]{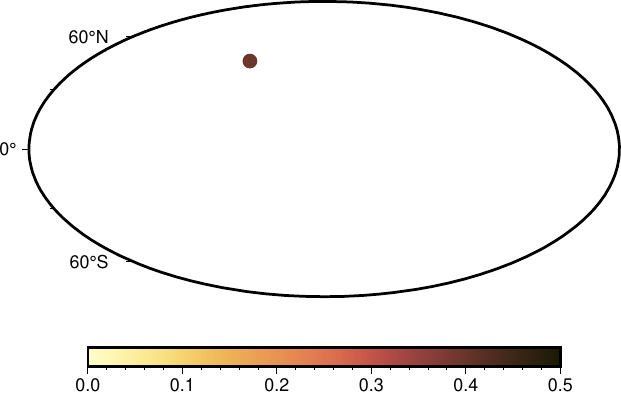}
\includegraphics[width=.32\textwidth]{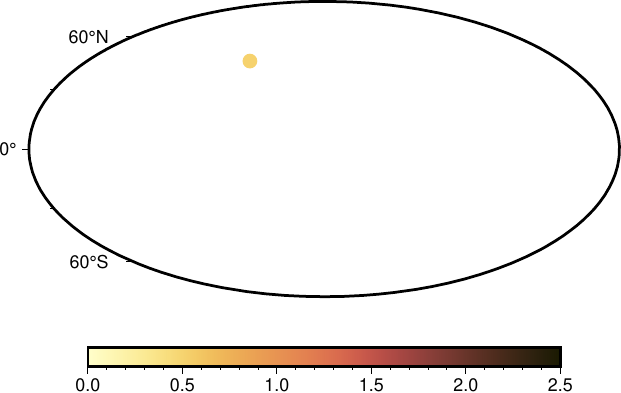}
\includegraphics[width=.32\textwidth]{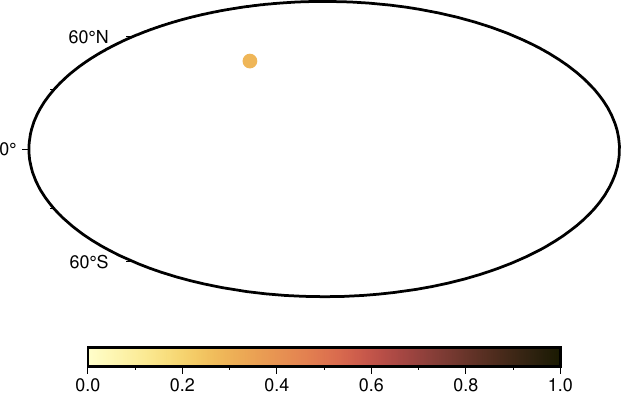}
\caption{4571 km - 4771 km}
\end{subfigure}

\begin{subfigure}{\textwidth}
\includegraphics[width=.32\textwidth]{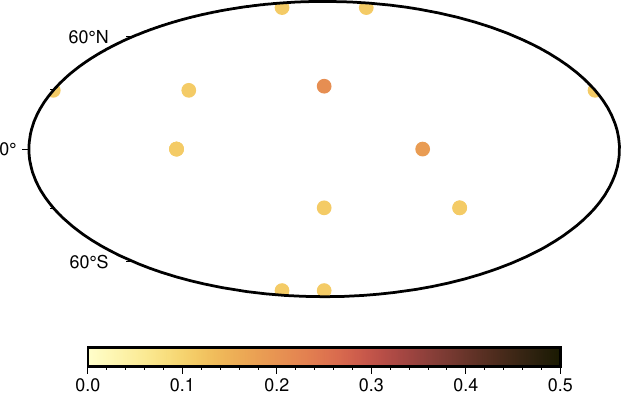}
\includegraphics[width=.32\textwidth]{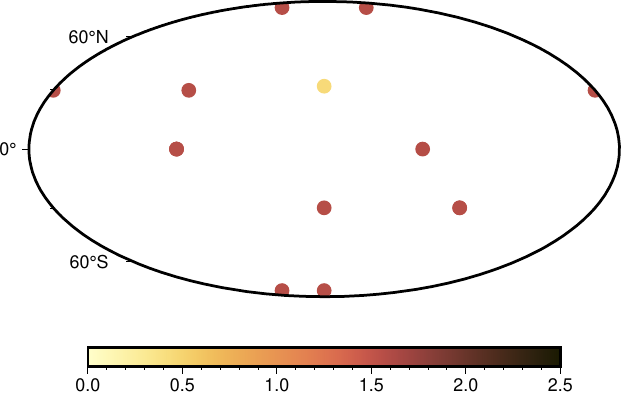}
\includegraphics[width=.32\textwidth]{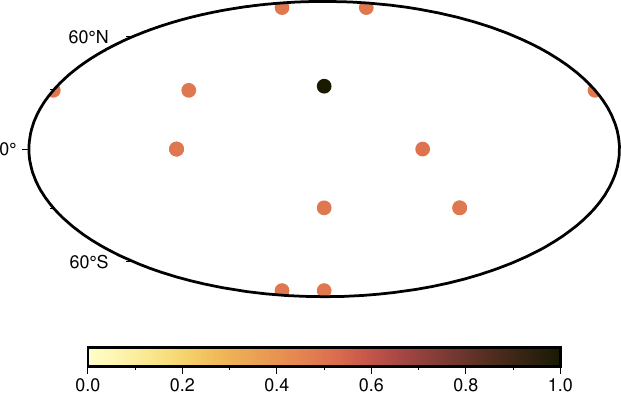}
\caption{4771 km - 4971 km}
\end{subfigure}
\caption{Chosen FEHFs within the given depth range (in terms of distances to the centre of the Earth): the positions of the dots refer to the positions of the tesseroids' centre and the colours refer to their size in the three polar-coordinate directions, namely $\Delta R$ (left column), $\Delta \Phi$ (middle) and $\Delta T$ (right).}
\label{fig:FEHFs1}
\end{figure}
\begin{figure}
\centering

\begin{subfigure}{\textwidth}
\includegraphics[width=.32\textwidth]{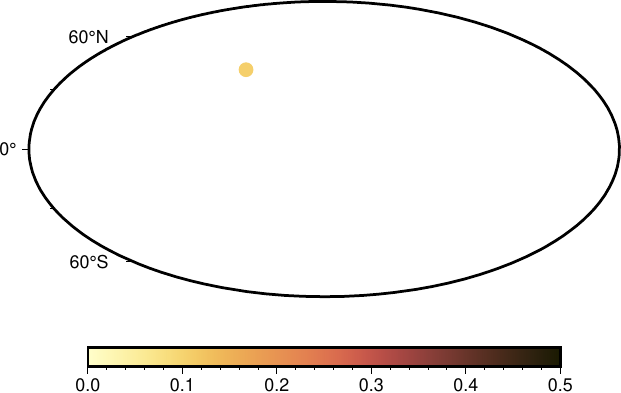}
\includegraphics[width=.32\textwidth]{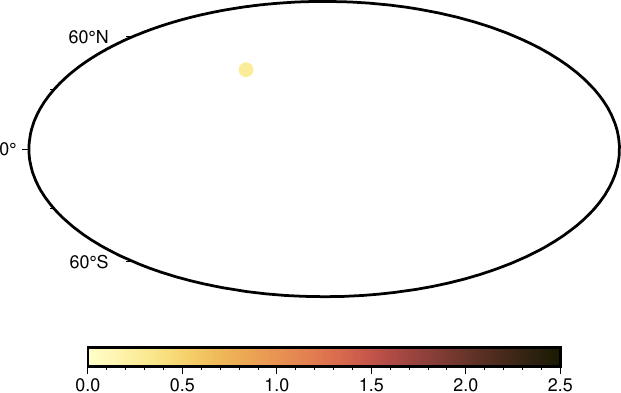}
\includegraphics[width=.32\textwidth]{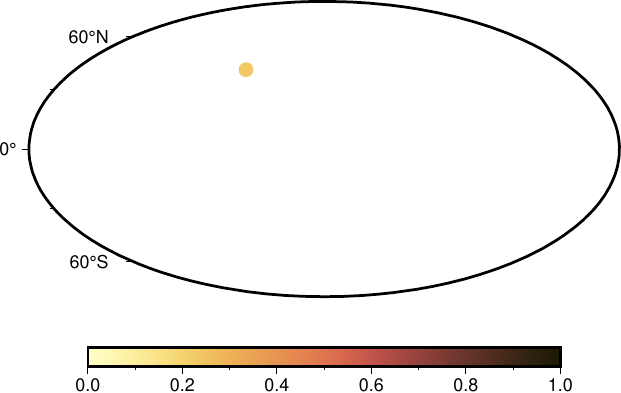}
\caption{5171 km - 5371 km}
\end{subfigure}

\begin{subfigure}{\textwidth}
\includegraphics[width=.32\textwidth]{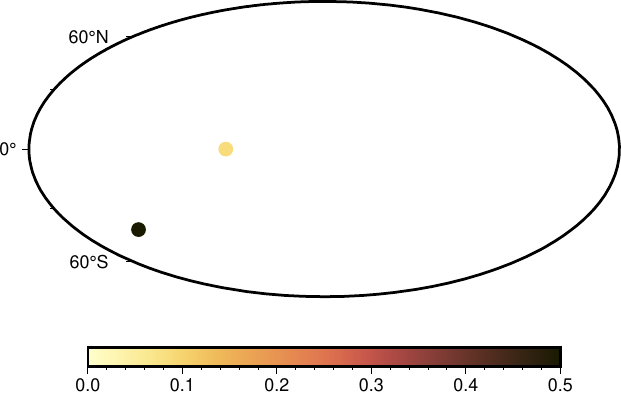}
\includegraphics[width=.32\textwidth]{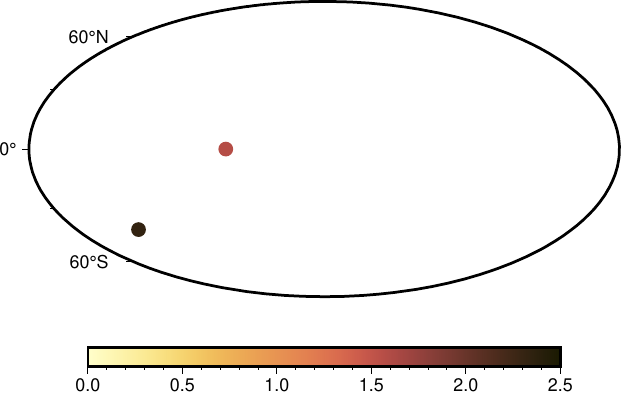}
\includegraphics[width=.32\textwidth]{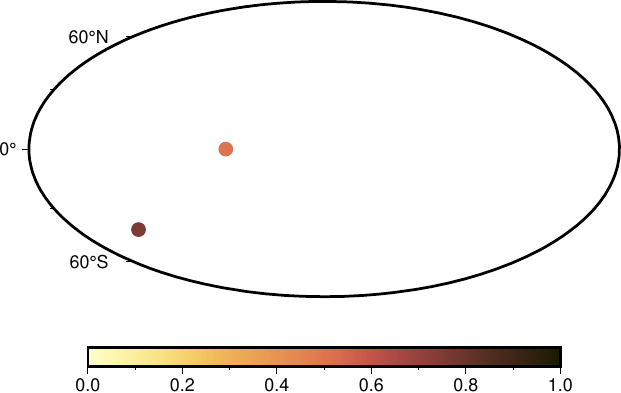}
\caption{5371 km - 5571 km}
\end{subfigure}

\begin{subfigure}{\textwidth}
\includegraphics[width=.32\textwidth]{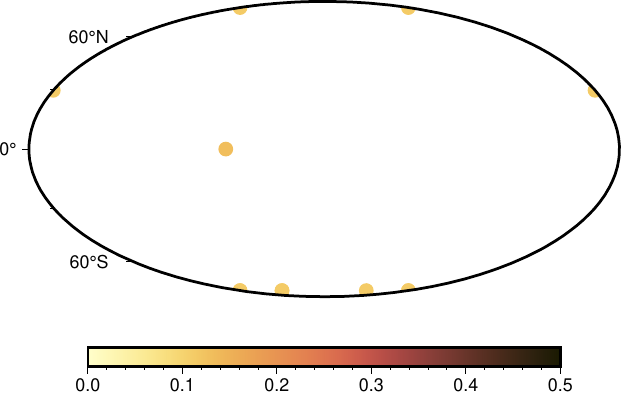}
\includegraphics[width=.32\textwidth]{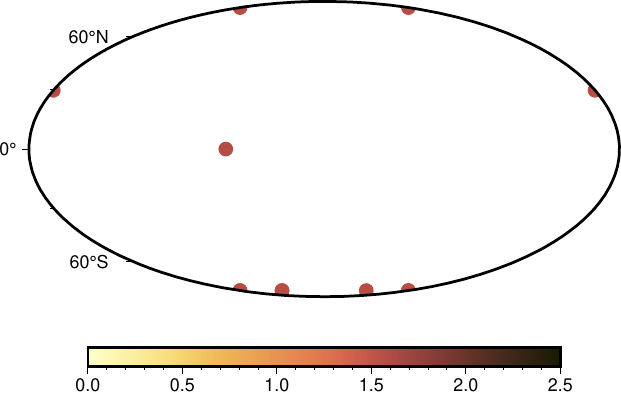}
\includegraphics[width=.32\textwidth]{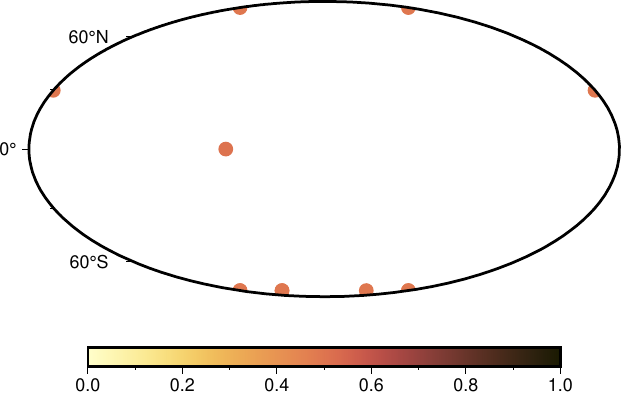}
\caption{5571 km - 5771 km}
\end{subfigure}

\begin{subfigure}{\textwidth}
\includegraphics[width=.32\textwidth]{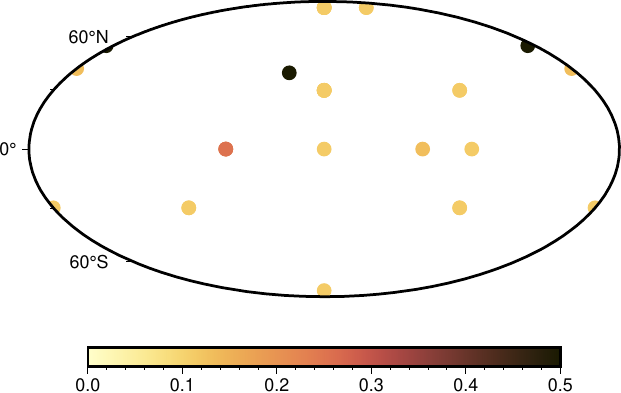}
\includegraphics[width=.32\textwidth]{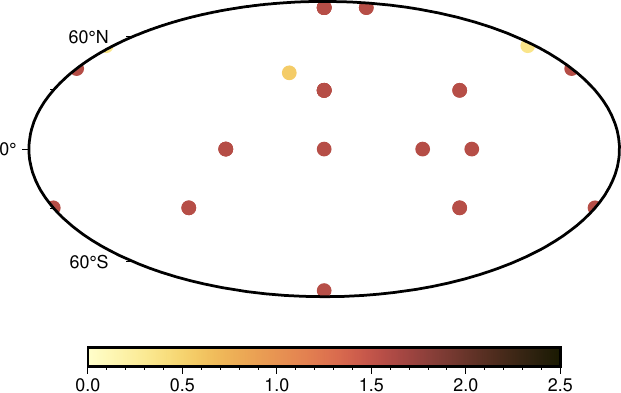}
\includegraphics[width=.32\textwidth]{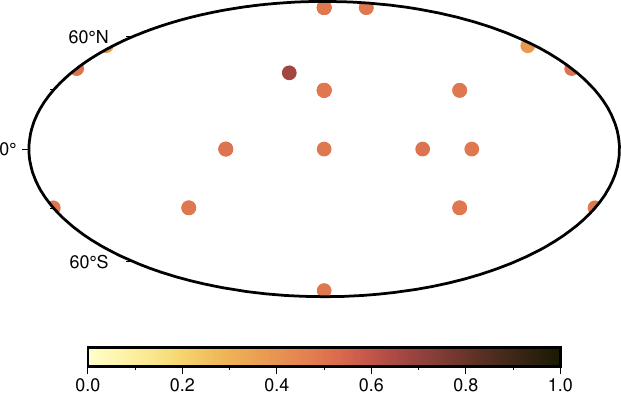}
\caption{6171 km - 6371 km}
\end{subfigure}
\caption{Illustration of the chosen FEHFs (see also \cref{fig:FEHFs1}).}
\label{fig:FEHFs2}
\end{figure}

\begin{figure}
\centering

\begin{subfigure}{\textwidth}
\includegraphics[width=.32\textwidth]{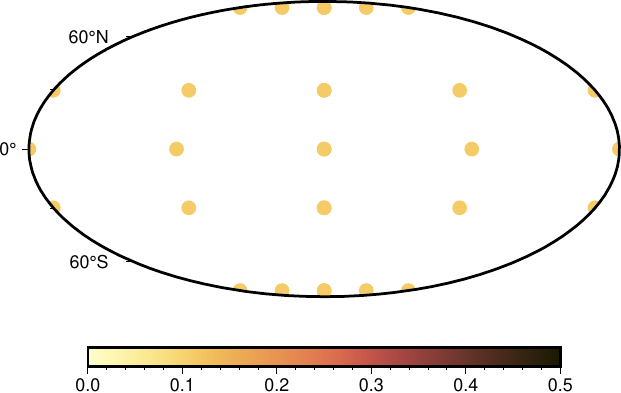}
\includegraphics[width=.32\textwidth]{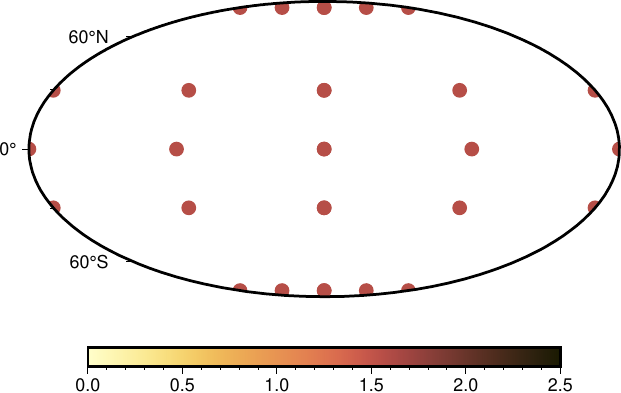}
\includegraphics[width=.32\textwidth]{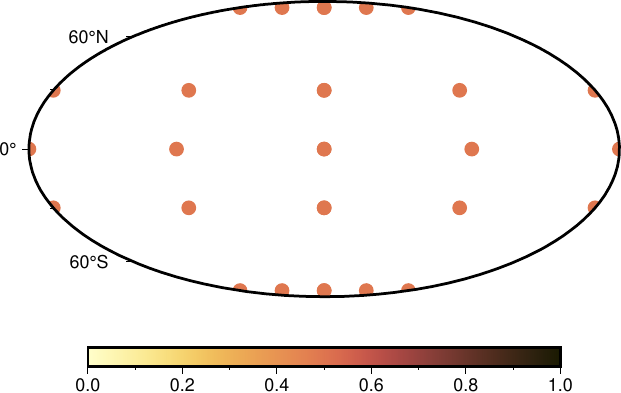}
\caption{3371 km - 3571 km}
\end{subfigure}

\begin{subfigure}{\textwidth}
\includegraphics[width=.32\textwidth]{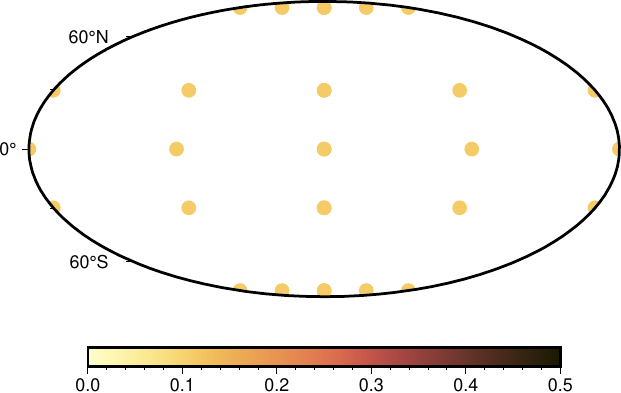}
\includegraphics[width=.32\textwidth]{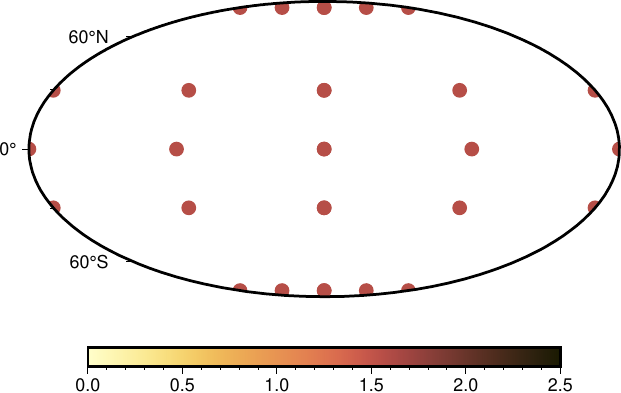}
\includegraphics[width=.32\textwidth]{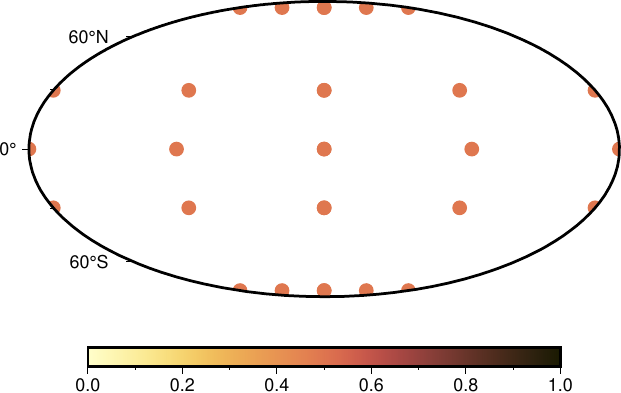}
\caption{4171 km - 4371 km}
\end{subfigure}

\begin{subfigure}{\textwidth}
\includegraphics[width=.32\textwidth]{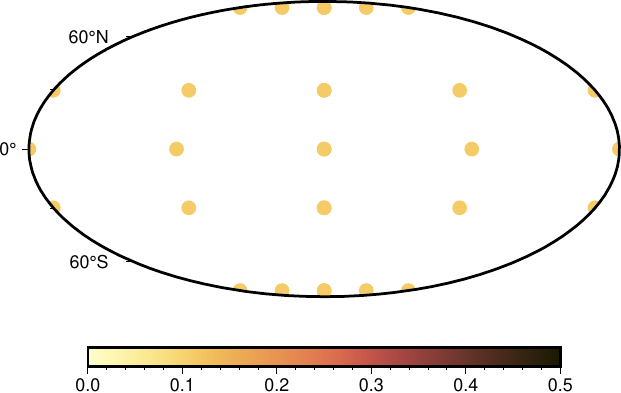}
\includegraphics[width=.32\textwidth]{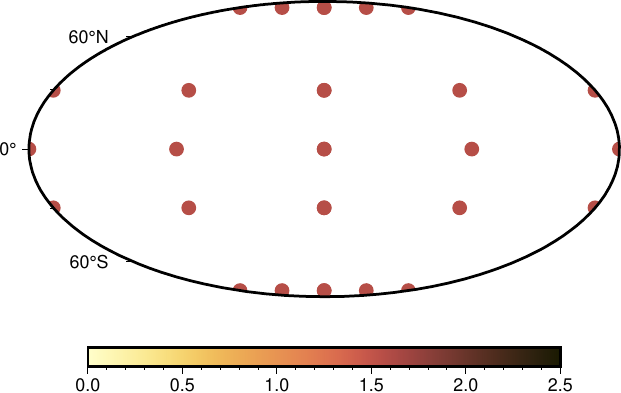}
\includegraphics[width=.32\textwidth]{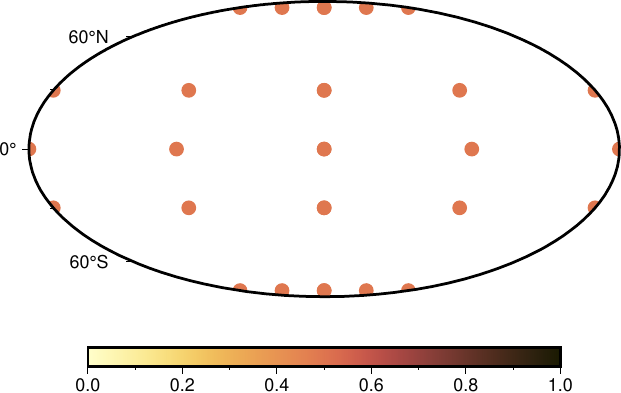}
\caption{4771 km - 4971 km}
\end{subfigure}

\begin{subfigure}{\textwidth}
\includegraphics[width=.32\textwidth]{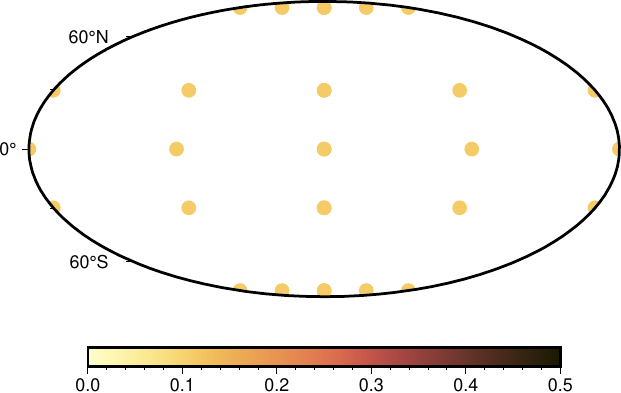}
\includegraphics[width=.32\textwidth]{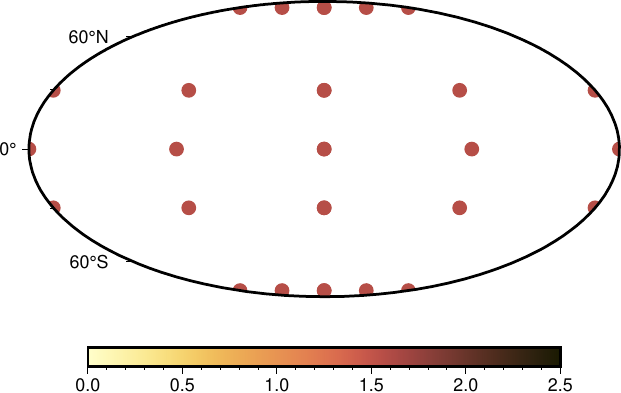}
\includegraphics[width=.32\textwidth]{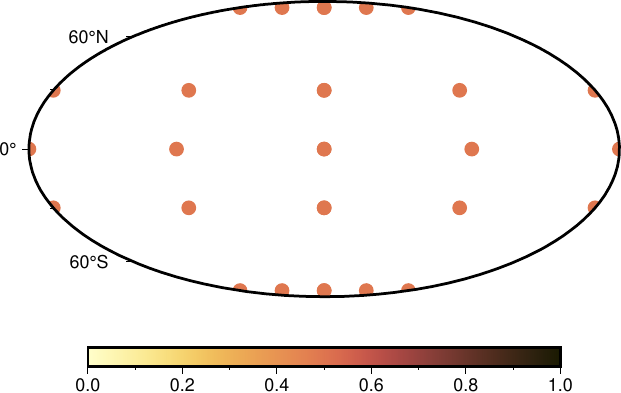}
\caption{5571 km - 5771 km}
\end{subfigure}

\begin{subfigure}{\textwidth}
\includegraphics[width=.32\textwidth]{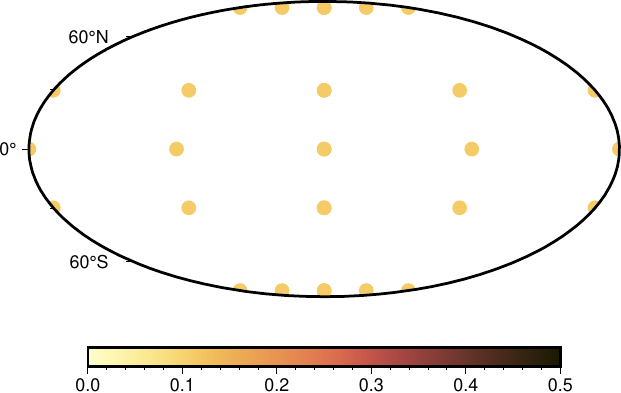}
\includegraphics[width=.32\textwidth]{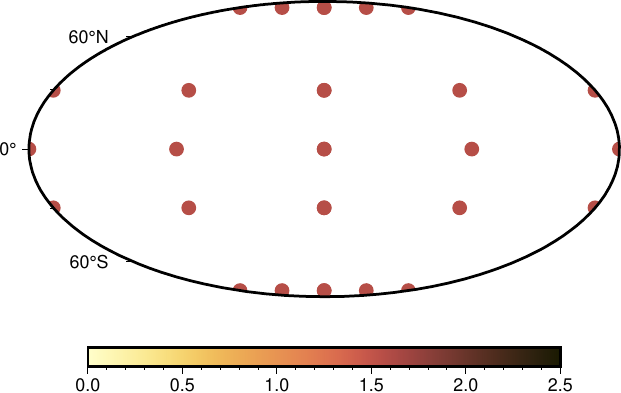}
\includegraphics[width=.32\textwidth]{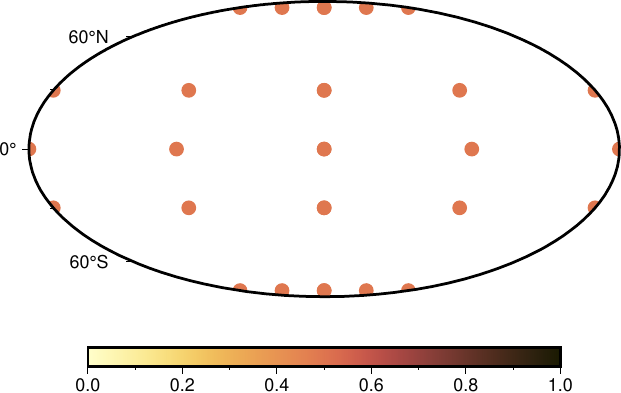}
\caption{6171 km - 6371 km}
\end{subfigure}
\caption{FEHFs of the starting dictionary within the given depth ranges (the illustration is analogous to \cref{fig:FEHFs1} and \cref{fig:FEHFs2}).}
\label{fig:DS_TFs}
\end{figure}

We chose the regularization parameter to be $10^{-3}\|y\|_{\real^\ell}$ as this produced the lowest RRMSE among the tested values of parameters. The experiment terminated after 300 iterations with an RRMSE of 0.881173 and a relative data error of 0.154958. The absolute approximation error is shown in \cref{fig:apprerr}. We scaled the figures to the size of the solution for better comparability (compare colour scales in \cref{fig:resolutiontests} and \cref{fig:apprerr}). We note that the errors are low for intermediate depths, where our teleseismic body wave paths sample the mantle most extensively, and higher for large and shallow depths. Recall that we have very unregularly distributed data (see \cref{fig:distrrays}).

In particular, at the distances from 3770.9 km to 4926.5 km to the centre, we see only very little remaining errors, and some of the errors are simply due to boundary effects since the method is unable to recover the precise geometry of the plumes. Moreover, those are close to the plumes, i.\,e.\, the region where our structure is given. Also in larger depths (distances 3193.1 km and 3482.0 km to the centre) the main errors are situated in the Northern hemisphere. Since there is practically no data, see \cref{fig:distrrays}, in these depths, the method cannot register that the plumes are cut of at the core-mantle-boundary and, further, it does not introduce random artefacts there but continues the structures. We assume that using additionally core-diffracted waves as well would erase the errors in these depths. In shallower depths (distances 5215.4 km to 6082.1 km to the centre), we also see a similar continuing behaviour of the approximation. Comparing the approximation in these depths with the data distribution, we observe that the errors do not increase as rapidly as the data becomes sparse. Hence, also here we have a continuation of structure within sparser data regions without many artefacts. The errors increase in particular below the Southern Pacific at radii 5504.3 km and 5793.2 km, which is typically not a very active seismic region and the data is not well-distributed there. Unfortunately, at the Earth's surface we have many artefacts. As the sparsity there is extremely high, we assume that a certain sparsity level can also be a limit to the method. However, note that the Earth's surface represents also the boundary of the considered regions and boundaries are known to possibly introduce additional challenges in approximation tasks.

In \cref{fig:Polys}, \cref{fig:FEHFs1} and \cref{fig:FEHFs2}, we show the chosen dictionary elements. Further, in \cref{fig:DS_TFs}, we provide a comparison to the FEHFs that were given in the starting dictionary. First of all, we note that the method chooses both local and global functions, much more polynomials than FEHFs to be precise. In particular, higher degrees and orders are preferred in general. This supports the idea that, due to the irregularly distributed data, the method uses global functions which fill empty blanks with similar structures as in regions with many data. This is at least the case, as long as there are enough data points nearby such as in medium depths. Thus, it would be interesting to investigate how the method would work with even higher degrees and orders. Maybe it would choose even less FEFHs then. Regarding the chosen FEHFs, we see that the LRFMP finds optimal ones in depths where more data is given though also there are not many selected. In regions where the data is sparse, the method falls back to the FEHFs given in the starting dictionary. In general this is not really desired. Note, however, that this also happens majorly near the Earth's surface which represents a specific challenge as discussed above. However, since the FEHFs are also generally less often chosen, this outcome suggests that FEHFs may not be the best choice within the LRFMP for travel time tomography but with only sparse data the polynomials are also not suitable. It remains a question for future research to determine which alternative types of trial functions in the dictionary can yield better numerical results. For other applications such as gravitational field modelling, the LRFMP proved to perform well for the combination of orthogonal polynomials (spherical harmonics in this case) and localized scaling functions and wavelets constructed from reproducing kernels, see e.g. \cite{Schneideretal2023-G}. However, the analogues of such kernels on balls are connected to essentially higher numerical costs, which is why we had not chosen them for our first experiments demonstrated here. Therefore, they might be promising but it remains open to improve the efficiency of their numerical calculation and integration.

Hence, we observe that the LRFMP can provide the possibility to reconstruct plumes within the Earth and constrain errors to their spatial locale, though the method certainly leaves the potential for further improvements.

\section{Discussion}
The underlying inverse problem, travel time tomography, is ill-posed. Unfortunately, we lack other, helpful theoretical insights into this problem such as a singular value decomposition of the corresponding operator which makes the numerical computations rather time-consuming. Due to the lack of numerically exploitable properties of the respective operator, we have to take special care on how to use a sensible amount of rays in practice. The amount we use here is, as far as we know, sufficient for global models. However, the number of rays has already provided a numerical challenge to the used inversion method.\\
Thus, for now and as we just started our development, we aimed for a proof of concept and not novel seismological insights here. Hence, the use of a ray theoretical setting is sufficient for now, in particular, as a finite frequency ansatz would even include more computational challenges. \\
All in all, the applicability of the LRFMP to travel time tomography is demonstrated quite fairly given the irregularly distributed data. Nonetheless, there are a number of open questions for future research.\\ 
First of all, we note that we should harmonize the used perturbation noise and the uncertainty $\sigma$ in future experiments. For use with real data, corrections for earthquake hypocentres ('source relocations') will need to be added in an efficient manner. Another crucial aspect that emerges from the result is the question which trial functions are best for this inverse problem. Though we tried FEHFs in order to obtain a model that is comparable with other approaches, the results point towards a different dictionary for future research. In general, we conclude that the LIPMPs can promise to provide an alternative numerical regularization method in comparison to other methods for travel time tomography. However, for achieving a competing status, the method still needs to be further enhanced.

\section{Conclusion and Outlook}
\label{sect:cons}
Here, we proposed to use the LRFMP algorithm which yields an approximation in a chosen best basis of dictionary elements. For the latter, we allowed polynomials as well as tesseroid-based linear finite element hat functions in order to include global and local trial functions in the dictionary. As this is the first application of the LRFMP to travel time tomography, we aimed for a proof of concept and follow the ray-theoretical approach. 

We presented here the method in general with an emphasis on aspects that have to be remodelled for travel time tomography. These are -- besides the choice of dictionary elements -- the choice of a regularization space and the practical use of a necessary amount of data. For the former, we chose a Sobolev space due to the deep mathematical connections between finite elements and these spaces. For the latter, we introduced an additional divide-and-conquer strategy which allowed us to consider packages of rays iteratively and, in this way, overall consider a suitable amount of rays to obtain a proof of concept. 

In our experiments, we considered a contrived Earth model consisting of two cone-shaped plumes between the core-mantle boundary and the Earth's surface. Our results showed that the LRFMP is able to reconstruct such structures. The approximation contains more errors where the data is rather sparse. 

Thus, future research could deal with a parameter study for the divide-and-conquer strategy in combination with other accuracy-lowering and efficiency-improving approaches as well as tackling more complicated contrived Earth models.

\section*{Declarations}

\paragraph{Funding} V.\, Michel gratefully acknowledges the support by the German Research Foundation (DFG; Deutsche Forschungsgemeinschaft), project MI 655/14-1. K. Sigloch was supported by the French government through the UCAJEDI Investments in the Future project, reference number ANR-15-IDEX-01. 
\paragraph{Conflicts of interest/Competing interests} Not applicable.
\paragraph{Availability of data, material and code} The code is available at \textcolor{blue}{\url{https://doi.org/10.5281/zenodo.8227888}} under licence the CC-BY-NC-SA 3.0 DE. The solution shown here can be computed from it. The used and corrected rays accompany the code in a compressed format. Note that the ISC-EHB meta-data is also freely available online.
\paragraph{Author's contributions} The research was carried out during N.\, Schneider's postdoc project. In this DFG-funded project, V.\, Michel was the principal investigator. K.\, Sigloch was the project partner from geosciences. Both supervised the project and assisted N.\, Schneider. E.\,J.\, Totten provided available data and software and assisted in the preparation of the tests.

\footnotesize
\bibliography{biblioneu}
\normalfont

\newpage
\appendix
\section{Mathematical derivation of specific terms}
\label{sect:app}

\subsection{Derivation of $\nabla N_{A,\Delta A}$}
\label{ssect:app:nablaFEHF}
As we are using the gradient in the $\lp{2}$-integrals, we need the version of $\nabla$ in spherical coordinates: 
\begin{align}
\nabla_{r\xi(\lon,t)} &= \era \pdervr + \frac{1}{r} \nabla^* 
= \era \pdervr + \frac{1}{r} \left( \ephi \frac{1}{\sqrt{1-t^2}}\pdervlon + \ete \sqrt{1-t^2}\pdervt \right)
\end{align}
see, e.g. \cite{Freedenetal1998,Freedenetal2013-1,Michel2013}.
Then, for an FEHF, we obtain 
\begin{align}
\nabla_{r\xi(\lon,t)} &N_{(R,\Phi,T),(\Delta R,\Delta \Phi,\Delta T)}(r,\lon,t) \\
&= 
\chi_{\mathrm{supp}_{[(R,\Phi,T)-(\Delta R,\Delta \Phi,\Delta T),(R,\Phi,T)+(\Delta R,\Delta \Phi,\Delta T)]}}(r,\lon,t)\\ 
&\qquad  \times \left( \era \pdervr N_{(R,\Phi,T),(\Delta R,\Delta \Phi,\Delta T)}(r,\lon,t)\right.\\ 
&\qquad \qquad  + \frac{1}{r}\ephi \frac{1}{\sqrt{1-t^2}}\pdervlon N_{(R,\Phi,T),(\Delta R,\Delta \Phi,\Delta T)}(r,\lon,t)\\
&\qquad \qquad \left. + \frac{1}{r}\ete \sqrt{1-t^2}\pdervt N_{(R,\Phi,T),(\Delta R,\Delta \Phi,\Delta T)}(r,\lon,t)
\right)\\
&= 
\chi_{\mathrm{supp}_{[(R,\Phi,T)-(\Delta R,\Delta \Phi,\Delta T),(R,\Phi,T)+(\Delta R,\Delta \Phi,\Delta T)]}}(r,\lon,t) \\
&\qquad\times \left(
\era \frac{[-\sgn(r-R)]}{\Delta R} \frac{\Delta \Phi-|\lon-\Phi|}{\Delta \Phi}\frac{\Delta T-|t-T|}{\Delta T}\right.\\ 
&\qquad\qquad  + \frac{1}{r}\ephi \frac{1}{\sqrt{1-t^2}} \frac{\Delta R-|r-R|}{\Delta R} \frac{[-\sgn(\lon-\Phi)]}{\Delta \Phi} \frac{\Delta T-|t-T|}{\Delta T}\\
&\qquad\qquad  \left. + \frac{1}{r}\ete \sqrt{1-t^2}\frac{\Delta R-|r-R|}{\Delta R}\frac{\Delta \Phi-|\lon-\Phi|}{\Delta \Phi}\frac{[-\sgn(t-T)]}{\Delta T}
\right)
\end{align}
almost everywhere because 
\begin{align}
\frac{\partial}{\partial a_k} \frac{\Delta A_k-|a_k-A_k|}{\Delta A_k}
= \frac{- \frac{\partial}{\partial a_k} |a_k-A_k|}{\Delta A_k}
= \frac{-\sgn(a_k-A_k)}{\Delta A_k}.
\end{align}
Note that this is only piecewise continuous with respect to the differentiated component but it is still continuous for the other components.

\subsection{Derivation of $\nabla G_{m,n,j}^\mathrm{I}$}
\label{ssect:app:nablaGI}

Similarly, for the polynomials with $n\geq 1$, we obtain 
\begin{align}
&\nabla_{r\xi(\lon,t)} G_{m,n,j}^{\mathrm{I}} (r\xi(\lon,t))\\
&= p_{m,n}\left[ \era \pdervr + \frac{1}{r} \nabla^*\right] \left[P_m^{(0,n+1/2)}\left(I(r)\right) \left(\frac{r}{\mathbf{R}}\right)^{n} Y_{n,j}(\xi(\lon,t))\right]\\
&= p_{m,n}\era \pdervr \left[ P_m^{(0,n+1/2)}\left(I(r)\right) \left(\frac{r}{\mathbf{R}}\right)^{n}\right] Y_{n,j}(\xi(\lon,t))\\
&\qquad + \frac{p_{m,n}}{r} \nabla^* \left[P_m^{(0,n+1/2)}\left(I(r)\right) \left(\frac{r}{\mathbf{R}}\right)^{n} Y_{n,j}(\xi(\lon,t))\right]\\
&=  p_{m,n}\left[ \left(P_m^{(0,n+1/2)}\left(I(r)\right)\right)'I'(r) \left(\frac{r}{\mathbf{R}}\right)^{n} + \frac{n}{\mathbf{R}} P_m^{(0,n+1/2)}\left(I(r)\right) \left(\frac{r}{\mathbf{R}}\right)^{n-1} \right]  \mu_n^{(1)}y^{(1)}_{n,j}(\xi(\lon,t))\\
&\qquad + \frac{p_{m,n}}{\mathbf{R}} P_m^{(0,n+1/2)}\left(I(r)\right) \left(\frac{r}{\mathbf{R}}\right)^{n-1} \mu_n^{(2)}y^{(2)}_{n,j}(\xi(\lon,t))
\label{eq:nablaGImnj}
\end{align}
with 
\begin{align}
\mu_n^{(i)} &\coloneqq \left\{\begin{matrix} 1,&i=1\\ \sqrt{n(n+1)},&i=2,3, \end{matrix} \right.
\label{def:mu}
\end{align}
and the vector spherical harmonics $y_{n,j}^{(i)}$, see e.g. \cite{DahlenTromp1998,Freedenetal2013-1,Freedenetal2009,Michel2020,MorseFeshbachI1953,MorseFeshbachII1953}.
In the case $n=0$, the angular derivative as well as the derivative of $(r/R)^{n}$ does not exist:
\begin{align}
\nabla_{r\xi(\lon,t)} G_{m,0,0}^{\mathrm{I}} (r\xi(\lon,t)) 
&= p_{m,0}\left[ \era \pdervr + \frac{1}{r} \nabla^*\right] P_m^{(0,1/2)}\left(I(r)\right) Y_{0,0}(\xi(\lon,t))\\
&= p_{m,0}\left[ \era \pdervr + \frac{1}{r} \nabla^*\right] P_m^{(0,1/2)}\left(I(r)\right) \frac{1}{\sqrt{4\pi}}\\
&= p_{m,0}\left(P_m^{(0,1/2)}\left(I(r)\right)\right)'I'(r) \left(\mu_0^{(1)}\right)y^{(1)}_{0,0}(\xi(\lon,t)).
\label{eq:nablaGIm00}
\end{align}
Note that \cref{eq:nablaGIm00} can be written in the form \cref{eq:nablaGImnj} due to the multiplication with $0$ (2nd term) and by defining $y_{0,0}^{(2)} \coloneqq 0$ (3rd term). Hence, the gradient \cref{eq:nablaGImnj} is well-defined for all $r$ and all possible $m,\ n$ and $j$. For practical purposes, we need to specify the gradient in more detail:
\begin{align}
&\nabla_{r\xi(\lon,t)} G_{m,n,j}^{\mathrm{I}}(r\xi(\lon,t)) \\
&=p_{m,n}\left[ \left(P_m^{(0,n+1/2)}\left(I(r)\right)\right)'I'(r) \left(\frac{r}{\mathbf{R}}\right)^{n} + \frac{n}{\mathbf{R}} P_m^{(0,n+1/2)}\left(I(r)\right) \left(\frac{r}{\mathbf{R}}\right)^{n-1} \right]  y^{(1)}_{n,j}(\xi(\lon,t))\\
&\qquad + \frac{p_{m,n}}{\mathbf{R}} P_m^{(0,n+1/2)}\left(I(r)\right) \left(\frac{r}{\mathbf{R}}\right)^{n-1}  \mu_n^{(2)}y^{(2)}_{n,j}(\xi(\lon,t))\label{eq:l2IPGI_1}\\
&= p_{m,n} \left[ \left(P_m^{(0,n+1/2)}\left(I(r)\right)\right)'I'(r) \left(\frac{r}{\mathbf{R}}\right)^{n} q_{n,j} P_{n,|j|}(t)\mathrm{Trig}(j\lon) \xi(\lon,t) \right.\\
&\qquad + \frac{n}{\mathbf{R}} P_m^{(0,n+1/2)}\left(I(r)\right) \left(\frac{r}{\mathbf{R}}\right)^{n-1} q_{n,j} P_{n,|j|}(t)\mathrm{Trig}(j\lon) \xi(\lon,t)\\
&\qquad + \left. \frac{1}{\mathbf{R}} P_m^{(0,n+1/2)}\left(I(r)\right)\left(\frac{r}{\mathbf{R}}\right)^{n-1} q_{n,j} \nabla^* \left(P_{n,|j|}(t)\mathrm{Trig}(j\lon)\right) \right]\\
&= p_{m,n}q_{n,j}\left[\left(P_m^{(0,n+1/2)}\left(I(r)\right)\right)'I'(r) \left(\frac{r}{\mathbf{R}}\right)^{n} P_{n,|j|}(t)\mathrm{Trig}(j\lon) \xi(\lon,t)\right.\\
&\qquad + \frac{n}{\mathbf{R}} P_m^{(0,n+1/2)}\left(I(r)\right) \left(\frac{r}{\mathbf{R}}\right)^{n-1} P_{n,|j|}(t)\mathrm{Trig}(j\lon) \xi(\lon,t)\\
&\qquad + \frac{1}{\mathbf{R}} P_m^{(0,n+1/2)}\left(I(r)\right) \left(\frac{r}{\mathbf{R}}\right)^{n-1} \frac{1}{\sqrt{1-t^2}} P_{n,|j|}(t)\left(\mathrm{Trig}(j\lon)\right)'\ephi(\lon,t)\\
&\qquad + \left. \frac{1}{\mathbf{R}} P_m^{(0,n+1/2)}\left(I(r)\right) \left(\frac{r}{\mathbf{R}}\right)^{n-1} \sqrt{1-t^2} P'_{n,|j|}(t)\mathrm{Trig}(j\lon)\ete(\lon,t) \right]
\end{align}
\begin{align}
&= p_{m,n}q_{n,j} \left[ \left(P_m^{(0,n+1/2)}\left(I(r)\right)\right)'I'(r) \left(\frac{r}{\mathbf{R}}\right)^{n} P_{n,|j|}(t)\mathrm{Trig}(j\lon) \xi(\lon,t)\right.\\
&\qquad + \frac{n}{\mathbf{R}} P_m^{(0,n+1/2)}\left(I(r)\right) \left(\frac{r}{\mathbf{R}}\right)^{n-1}  P_{n,|j|}(t)\mathrm{Trig}(j\lon) \xi(\lon,t)\\
&\qquad + \frac{j}{\mathbf{R}} P_m^{(0,n+1/2)}\left(I(r)\right) \left(\frac{r}{\mathbf{R}}\right)^{n-1} \frac{1}{\sqrt{1-t^2}} P_{n,|j|}(t)\mathrm{Trig}(-j\lon)\ephi(\lon,t)\\
&\qquad + \left.\frac{1}{\mathbf{R}} P_m^{(0,n+1/2)}\left(I(r)\right) \left(\frac{r}{\mathbf{R}}\right)^{n-1}  \sqrt{1-t^2} P'_{n,|j|}(t)\mathrm{Trig}(j\lon)\ete(\lon,t)\right]\\
&\eqqcolon p_{m,n}q_{n,j} \sum_{p=1}^4 G_{m,n,j;p}^{\mathrm{I}} (r\xi(\lon,t)).
\end{align}
Note that the well-definedness of the terms
\begin{align}
\frac{P_{n,|j|}(t)}{\sqrt{1-t^2}} \qquad \textrm{and} \qquad \sqrt{1-t^2} P'_{n,|j|}(t)
\end{align}
used only for $n\geq 1$ was already discussed in \cite{Schneider2020}. The latter can also be computed with the respective algorithm given there. We discuss the former here in a bit more detail as this increases the efficiency in our implementation. Away from the poles, the term can be calculated straightforwardly. In a neighbourhood of the North and South Pole (i.e. for $t\to \pm 1$), we obtain 
\begin{align}
\lim_{t\to\pm 1} \frac{P_{n,|j|}(t)}{\sqrt{1-t^2}} 
= \lim_{t\to\pm 1} \left(1-t^2\right)^{(|j|-1)/2} P^{(|j|)}_{n}(t)
= \left\{ \begin{matrix}
\lim_{t\to\pm 1} \left(1-t^2\right)^{0} P'_{n}(t) 
= P'_{n}(\pm 1), & |j|=1\\
0, & |j|>1,
\end{matrix}\right\}
\end{align}
where the lower row holds because every derivative of a Legendre polynomial is bounded in $[-1,1]$ for trivial reasons.

\subsection{Derivation of regularization terms}
\label{ssect:app:H1IPs}

At first, we list certain trigonometric identities that we will need hereafter: 
\begin{align}
&\sin(\alpha) + \sin(\beta) = 2\sin\left(\frac{\alpha+\beta}{2}\right)\cos\left(\frac{\alpha-\beta}{2}\right),& 
&\sin(\alpha) - \sin(\beta) = 2\cos\left(\frac{\alpha+\beta}{2}\right)\sin\left(\frac{\alpha-\beta}{2}\right),\\
&\cos(\alpha) + \cos(\beta) = 2\cos\left(\frac{\alpha+\beta}{2}\right)\cos\left(\frac{\alpha-\beta}{2}\right),&
&\cos(\alpha) - \cos(\beta) = -2\sin\left(\frac{\alpha+\beta}{2}\right)\cos\left(\frac{\alpha-\beta}{2}\right),\\
&\sin(\alpha)\sin(\beta) = \frac{1}{2}\left(\cos\left(\alpha-\beta\right) - \cos\left(\alpha + \beta\right) \right),&
&\cos(\alpha)\cos(\beta) = \frac{1}{2}\left(\cos\left(\alpha-\beta\right) + \cos\left(\alpha + \beta\right) \right),\\
&\sin(\alpha)\cos(\beta) = \frac{1}{2}\left(\sin\left(\alpha-\beta\right) - \sin\left(\alpha + \beta\right) \right),&
&\int \sin(x)\cos(x) \intd x = \frac{\sin^2(x)}{2},\\
&\int x\sin(ax) \intd x = \frac{\sin(ax)}{a^2} - \frac{x\cos(ax)}{a},&
&\int x\cos(ax) \intd x = \frac{\cos(ax)}{a^2} + \frac{x\sin(ax)}{a},\\
&\int \sin^2(x) \intd x = \frac{x}{2} - \frac{\sin(2x)}{4},&
&\int \cos^2(x) \intd x = \frac{x}{2} + \frac{\sin(2x)}{4}.
\end{align}
We have to discuss the following inner products 
\begin{align}
&\left\langle G_{m,n,j}^{\mathrm{I}}, G_{m',n',j'}^{\mathrm{I}} \right\rangle_{\mathcal{H}^1},&
&\left\langle N_{A,\Delta A}, N_{A',(\Delta A)'}  \right\rangle_{\mathcal{H}^1}&
&\left\langle N_{A,\Delta A}, G_{m,n,j}^{\mathrm{I}} \right\rangle_{\mathcal{H}^1}
\intertext{which practically means computing}
&\left\langle G_{m,n,j}^{\mathrm{I}}, G_{m',n',j'}^{\mathrm{I}} \right\rangle_{\Lp{2}},&
&\left\langle N_{A,\Delta A}, N_{A',(\Delta A)'}  \right\rangle_{\Lp{2}},&
&\left\langle N_{A,\Delta A}, G_{m,n,j}^{\mathrm{I}} \right\rangle_{\Lp{2}},\\
&\left\langle \nabla G_{m,n,j}^{\mathrm{I}}, \nabla G_{m',n',j'}^{\mathrm{I}} \right\rangle_{\lp{2}},&
&\left\langle \nabla N_{A,\Delta A}, \nabla N_{A',(\Delta A)'}  \right\rangle_{\lp{2}},&
&\left\langle \nabla N_{A,\Delta A}, \nabla G_{m,n,j}^{\mathrm{I}} \right\rangle_{\lp{2}}.
\end{align}
We start with two polynomials. We have
\begin{align}
\left\langle G_{m,n,j}^{\mathrm{I}}, G_{m',n',j'}^{\mathrm{I}} \right\rangle_{\Lp{2}}
= \delta_{m,m'}\delta_{n,n'}\delta_{j,j'} 
\end{align}
due to the fact that they constitute a basis system in $\Lp{2}$. Next, we consider the vectorial inner product of their gradients. Note that the vector spherical harmonics are orthonormal with respect to their degrees, order and type. With \cref{eq:gradIr}, \cref{def:mu} and the substitution $r=\mathbf{R}\sqrt{(1+u)/2}$, we obtain from \cref{eq:l2IPGI_1}
\begin{align}
&\left\langle \nabla G_{m,n,j}^{\mathrm{I}}, \nabla G_{m',n',j'}^{\mathrm{I}}\right\rangle_{\lp{2}}\\
&= \delta_{n,n'}\delta_{j,j'} p_{m,n}p_{m',n} \\
&\qquad \times \left[\int_0^{\mathbf{R}} \left[ \left(P_m^{(0,n+1/2)}\left(I(r)\right)\right)'I'(r) \left(\frac{r}{\mathbf{R}}\right)^{n} + \frac{n}{\mathbf{R}} P_m^{(0,n+1/2)}\left(I(r)\right) \left(\frac{r}{\mathbf{R}}\right)^{n-1} \right]\right.\\
&\qquad\qquad\qquad \times \left[ \left(P_{m'}^{(0,n+1/2)}\left(I(r)\right)\right)'I'(r) \left(\frac{r}{\mathbf{R}}\right)^{n} + \frac{n}{\mathbf{R}} P_{m'}^{(0,n+1/2)}\left(I(r)\right) \left(\frac{r}{\mathbf{R}}\right)^{n-1} \right] r^2 \intd r \\
&\qquad\qquad + \left. \left(\mu_n^{(2)}\right)^2\int_0^{\mathbf{R}} \left[\frac{1}{\mathbf{R}} P_m^{(0,n+1/2)}\left(I(r)\right) \left(\frac{r}{\mathbf{R}}\right)^{n-1} \right]\left[\frac{1}{\mathbf{R}} P_{m'}^{(0,n+1/2)}\left(I(r)\right) \left(\frac{r}{\mathbf{R}}\right)^{n-1} \right] r^2 \intd r \right]\\
&= \delta_{n,n'}\delta_{j,j'} p_{m,n}p_{m',n} \\
&\qquad \times \left[\left(\mu_n^{(1)}\right)^2\int_0^{\mathbf{R}} \left(P_m^{(0,n+1/2)}\left(I(r)\right)\right)'I'(r) \left(\frac{r}{\mathbf{R}}\right)^{n} \left(P_{m'}^{(0,n+1/2)}\left(I(r)\right)\right)'I'(r) \left(\frac{r}{\mathbf{R}}\right)^{n} r^2 \intd r \right. \\
&\qquad\qquad + \left(\mu_n^{(1)}\right)^2\int_0^{\mathbf{R}} \left(P_m^{(0,n+1/2)}\left(I(r)\right)\right)'I'(r) \left(\frac{r}{\mathbf{R}}\right)^{n}\frac{n}{\mathbf{R}} P_{m'}^{(0,n+1/2)}\left(I(r)\right) \left(\frac{r}{\mathbf{R}}\right)^{n-1} r^2 \intd r \\
&\qquad\qquad + \left(\mu_n^{(1)}\right)^2\int_0^{\mathbf{R}}\frac{n}{\mathbf{R}} P_m^{(0,n+1/2)}\left(I(r)\right)\left(\frac{r}{\mathbf{R}}\right)^{n-1} \left(P_{m'}^{(0,n+1/2)}\left(I(r)\right)\right)' I'(r) \left(\frac{r}{\mathbf{R}}\right)^{n} r^2 \intd r \\
&\qquad\qquad  + \left(\mu_n^{(1)}\right)^2\int_0^{\mathbf{R}} \frac{n}{\mathbf{R}} P_m^{(0,n+1/2)}\left(I(r)\right) \left(\frac{r}{\mathbf{R}}\right)^{n-1} \frac{n}{\mathbf{R}} P_{m'}^{(0,n+1/2)}\left(I(r)\right) \left(\frac{r}{\mathbf{R}}\right)^{n-1} r^2 \intd r\\
&\qquad\qquad + \left. \left(\mu_n^{(2)}\right)^2\int_0^{\mathbf{R}} \frac{1}{\mathbf{R}} P_m^{(0,n+1/2)}\left(I(r)\right) \left(\frac{r}{\mathbf{R}}\right)^{n-1}\frac{1}{\mathbf{R}} P_{m'}^{(0,n+1/2)}\left(I(r)\right) \left(\frac{r}{\mathbf{R}}\right)^{n-1}  r^2 \intd r \right]\\
&= \delta_{n,n'}\delta_{j,j'} p_{m,n}p_{m',n} \\
&\qquad \times \left[\int_0^{\mathbf{R}} \left(P_m^{(0,n+1/2)}\left(I(r)\right)\right)'\left(P_{m'}^{(0,n+1/2)}\left(I(r)\right)\right)'(I'(r))^2 \left(\frac{r}{\mathbf{R}}\right)^{2n} r^2 \intd r \right. \\
&\qquad\qquad\qquad + \frac{n}{\mathbf{R}}\int_0^{\mathbf{R}} \left[ \left(P_m^{(0,n+1/2)}\left(I(r)\right)\right)'P_{m'}^{(0,n+1/2)}\left(I(r)\right)\right.\\
&\qquad\qquad\qquad\qquad\qquad \left. + P_m^{(0,n+1/2)}\left(I(r)\right)\left(P_{m'}^{(0,n+1/2)}\left(I(r)\right)\right)'\right] I'(r) \left(\frac{r}{\mathbf{R}}\right)^{2n-1} r^2 \intd r \\
&\qquad\qquad\qquad  +\left. \left( n^2 + n(n+1)\right)\int_0^{\mathbf{R}} P_m^{(0,n+1/2)}\left(I(r)\right)P_{m'}^{(0,n+1/2)}\left(I(r)\right) \left(\frac{r}{\mathbf{R}}\right)^{2n}  \intd r \right]
\end{align}
\begin{align}
&= \delta_{n,n'}\delta_{j,j'} p_{m,n}p_{m',n} \\
&\qquad \times \left[\int_0^{\mathbf{R}} \left(P_m^{(0,n+1/2)}\left(I(r)\right)\right)'\left(P_{m'}^{(0,n+1/2)}\left(I(r)\right)\right)'\frac{16r^2}{\mathbf{R}^4} \left(\frac{r}{\mathbf{R}}\right)^{2n} r^2 \intd r \right. \\
&\qquad\qquad\qquad + \frac{n}{\mathbf{R}}\int_0^{\mathbf{R}} \left[ \left(P_m^{(0,n+1/2)}\left(I(r)\right)\right)'P_{m'}^{(0,n+1/2)}\left(I(r)\right)\right.\\
&\qquad\qquad\qquad\qquad\qquad \left. + P_m^{(0,n+1/2)}\left(I(r)\right)\left(P_{m'}^{(0,n+1/2)}\left(I(r)\right)\right)'\right] \frac{4r}{\mathbf{R}^2} \left(\frac{r}{\mathbf{R}}\right)^{2n-1} r^2 \intd r \\
&\qquad\qquad\qquad  +\left. \left( n^2 + n^2+n\right)\int_0^{\mathbf{R}} P_m^{(0,n+1/2)}\left(I(r)\right)P_{m'}^{(0,n+1/2)}\left(I(r)\right) \left(\frac{r}{\mathbf{R}}\right)^{2n}  \intd r \right]\\
&= \delta_{n,n'}\delta_{j,j'} p_{m,n}p_{m',n} \\
&\qquad \times \left[16\int_0^{\mathbf{R}} \left(P_m^{(0,n+1/2)}\left(I(r)\right)\right)'\left(P_{m'}^{(0,n+1/2)}\left(I(r)\right)\right)' \left(\frac{r}{\mathbf{R}}\right)^{2n+4} \intd r \right. \\
&\qquad\qquad\qquad + 4n\int_0^{\mathbf{R}} \left[ \left(P_m^{(0,n+1/2)}\left(I(r)\right)\right)'P_{m'}^{(0,n+1/2)}\left(I(r)\right)\right.\\
&\qquad\qquad\qquad\qquad\qquad \left. + P_m^{(0,n+1/2)}\left(I(r)\right)\left(P_{m'}^{(0,n+1/2)}\left(I(r)\right)\right)'\right] \left(\frac{r}{\mathbf{R}}\right)^{2n+2} \intd r \\
&\qquad\qquad\qquad  +\left. \left( 2n^2 + n\right)\int_0^{\mathbf{R}} P_m^{(0,n+1/2)}\left(I(r)\right)P_{m'}^{(0,n+1/2)}\left(I(r)\right) \left(\frac{r}{\mathbf{R}}\right)^{2n}  \intd r \right]\\
&= \delta_{n,n'}\delta_{j,j'} p_{m,n}p_{m',n} \\
&\qquad \times \left[16\int_{-1}^{1} \left(P_m^{(0,n+1/2)}(u)\right)'\left(P_{m'}^{(0,n+1/2)}(u)\right)' \left(\frac{1+u}{2}\right)^{n+2} \frac{\mathbf{R}}{4}\sqrt{\frac{2}{1+u}}\intd u \right. \\
&\qquad\qquad\qquad + 4n\int_{-1}^{1} \left[ \left(P_m^{(0,n+1/2)}(u)\right)'P_{m'}^{(0,n+1/2)}(u)\right.\\
&\qquad\qquad\qquad\qquad\qquad \left. + P_m^{(0,n+1/2)}(u)\left(P_{m'}^{(0,n+1/2)}(u)\right)'\right] \left(\frac{1+u}{2}\right)^{n+1} \frac{\mathbf{R}}{4}\sqrt{\frac{2}{1+u}}\intd u \\
&\qquad\qquad\qquad  +\left. \left( n(2n+1)\right)\int_{-1}^{1} P_m^{(0,n+1/2)}(u)P_{m'}^{(0,n+1/2)}(u) \left(\frac{1+u}{2}\right)^{n}  \frac{\mathbf{R}}{4}\sqrt{\frac{2}{1+u}}\intd u \right]\\
&= \delta_{n,n'}\delta_{j,j'} p_{m,n}p_{m',n} \\
&\quad \times \left[4\mathbf{R}\int_{-1}^{1} \left(P_m^{(0,n+1/2)}(u)\right)'\left(P_{m'}^{(0,n+1/2)}(u)\right)' \left(\frac{1+u}{2}\right)^{n+3/2} \intd u \right. \\
&\qquad\qquad\qquad + \mathbf{R}n\int_{-1}^{1} \left[ \left(P_m^{(0,n+1/2)}(u)\right)'P_{m'}^{(0,n+1/2)}(u)\right.\\
&\qquad\qquad\qquad\qquad\qquad \left. + P_m^{(0,n+1/2)}(u)\left(P_{m}^{(0,n+1/2)}(u)\right)'\right] \left(\frac{1+u}{2}\right)^{n+1/2} \intd u \\
&\qquad\qquad\qquad  +\left. \frac{\mathbf{R} n(2n+1)}{4}\int_{-1}^{1} P_m^{(0,n+1/2)}(u)P_{m'}^{(0,n+1/2)}(u) \left(\frac{1+u}{2}\right)^{n-1/2} \intd u \right]
\end{align}
\begin{align}
&= \delta_{n,n'}\delta_{j,j'} p_{m,n}p_{m',n} \\
&\quad \times \left[\frac{\mathbf{R}\sqrt{2}}{2^{n}}\int_{-1}^{1} \left(P_m^{(0,n+1/2)}(u)\right)'\left(P_{m'}^{(0,n+1/2)}(u)\right)' \left(1+u\right)^{n+3/2} \intd u \right. \\
&\qquad\qquad\qquad + \frac{\mathbf{R}n}{2^{n}\sqrt{2}}\int_{-1}^{1} \left[ \left(P_m^{(0,n+1/2)}(u)\right)'P_{m'}^{(0,n+1/2)}(u)\right.\\
&\qquad\qquad\qquad\qquad\qquad \left. + P_m^{(0,n+1/2)}(u)\left(P_{m'}^{(0,n+1/2)}(u)\right)'\right] \left(1+u\right)^{n+1/2} \intd u \\
&\qquad\qquad\qquad  +\left. \frac{\mathbf{R} n(2n+1)}{2^{n+1}\sqrt{2}}\int_{-1}^{1} P_m^{(0,n+1/2)}(u)P_{m'}^{(0,n+1/2)}(u) \left(1+u\right)^{n-1/2} \intd u \right].
\end{align}

Next, we discuss the inner products of two finite element hat functions. Generally, we have 
\begin{align}
&\left\langle N_{A,\Delta A}, N_{A',(\Delta A)'}\right\rangle_{\Lp{2}}\\
&\qquad = \int_\ball \chi_{\mathrm{supp}_{A,\Delta A}}(a) \chi_{\mathrm{supp}_{A',(\Delta A)'}}(a)\prod_{j=1}^{3} \frac{[\Delta A_j-|a_j-A_j|][(\Delta A_j)'-|a_j-A'_j|]}{\Delta A_j (\Delta A_j)'} \intd a\\
&\qquad =\prod_{j=1}^{3} \int_{lb_{a_j}}^{ub_{a_j}} \frac{[\Delta A_j-|a_j-A_j|][(\Delta A_j)'-|a_j-A'_j|]}{\Delta A_j (\Delta A_j)'} \left\{\begin{matrix} a_1^2,j=1,\\ 1,\ j=2,3\end{matrix} \right\} \intd a_j
\end{align}
where $lb_a$ denotes the lower and $ub_a$ the upper bound of the respective $a_j$-integral. That is, principally, we have (informally speaking)
\begin{align}
[lb,ub] = \mathrm{supp}_{A,\Delta A} \cap \mathrm{supp}_{A',(\Delta A)'}
\end{align}
though we have to take a detailed look on the case $a_2$ if $P\not=P'$ (see below).
As a generalization, we consider the integrals 
\begin{align}
\int_{lb_{x}}^{ub_{x}} \frac{[\Delta X-|x-X|][(\Delta X)'-|x-X'|]}{\Delta X (\Delta X)'} x^q \intd x
\label{def:intNNa}
\end{align}
for $q\in\nat_0$ (here, in particular, $q=0$ and $q=2$)
in the sequel. Then, the $\Lp{2}$-integral is obtained as the product of this integral for all variables $x=a_j,\ j=1,2,3$. At first, we need to consider how the lower and upper bounds are defined. Furthermore, because the FEHFs are piecewise defined, we have to determine these critical points in between $lb_x$ and $ub_x$ as well. At last, we need to discuss the remaining integrals. \\
For $x=a_1=r$ and $x=a_3=t$, there are three cases depending on $X,\ \Delta X,\ X'$ and $(\Delta X)'$ that we have to consider: 
\begin{compactitem}
\item[(a)] $\mathrm{supp}_{X,\Delta X}= [X-\Delta X, X+\Delta X]$
	\begin{compactitem}
	\item[(a.1)] $\mathrm{supp}_{X',(\Delta X)'}= [X'-(\Delta X)', X+(\Delta X)']$ 
	\item[(a.2)] $\mathrm{supp}_{X',(\Delta X)'}= [X_{\mathrm{min}}, X'+(\Delta X)']$ 
	\item[(a.3)] $\mathrm{supp}_{X',(\Delta X)'}= [X'-(\Delta X)', X_{\mathrm{max}}]$ 
	\end{compactitem}
\item[(b)] $\mathrm{supp}_{X,\Delta X}= [X_{\mathrm{min}}, X+\Delta X]$ 
	\begin{compactitem}
	\item[(b.1)] $\mathrm{supp}_{X',(\Delta X)'}= [X'-(\Delta X)', X'+(\Delta X)']$
	\item[(b.2)] $\mathrm{supp}_{X',(\Delta X)'}= [X_{\mathrm{min}}, X'+(\Delta X)']$ 
	\item[(c.3)] $\mathrm{supp}_{X',(\Delta X)'}= [X'-(\Delta X)', X_{\mathrm{max}}]$ 
	\end{compactitem}
\item[(c)] $\mathrm{supp}_{X,\Delta X}= [X-\Delta X, X_{\mathrm{max}}]$ 
	\begin{compactitem}
	\item[(c.1)] $\mathrm{supp}_{X',(\Delta X)'}= [X'-(\Delta X)', X'+(\Delta X)']$
	\item[(c.2)] $\mathrm{supp}_{X',(\Delta X)'}= [X_{\mathrm{min}}, X'+(\Delta X)']$ 
	\item[(c.3)] $\mathrm{supp}_{X',(\Delta X)'}= [X'-(\Delta X)', X_{\mathrm{max}}]$ 
	\end{compactitem}
\end{compactitem}
In any of these cases, we obtain the lower bound as 
\begin{align}
lb_x &= \max[ \max(X_{\mathrm{min}},X-\Delta X), \max(X_{\mathrm{min}},X'-(\Delta X)')]
\label{def:lbxgen}
\end{align}
and the upper bound as 
\begin{align}
ub_x &= \min [\min(X_{\mathrm{max}},X+\Delta X), \min(X_{\mathrm{max}},X'+(\Delta X)')]
\label{def:ubxgen}
\end{align}
Note that if $lb_x \geq ub_x$ (componentwise), the intersection is the empty set and the integral vanishes. For $x=a_2=\lon$, we also have three relevant cases: 
\begin{compactitem}
\item[(a)] $\sgn(P)=-1 \Rightarrow \mathrm{supp}_{X,\Delta X}= [\max(-\pi,X-\Delta X), \min(\pi,X+\Delta X)]$
	\begin{compactitem}
	\item[(a.1)] $\sgn(P')=-1 \Rightarrow \mathrm{supp}_{X',(\Delta X)'}= [\max(-\pi,X'-(\Delta X)'), \min(\pi,X+(\Delta X)')]$
	\item[(a.2)] $\sgn(P')=0 \Rightarrow \mathrm{supp}_{X',(\Delta X)'}= [\max(0,X-(\Delta X)'), \min(2\pi,X+(\Delta X)')]$
	\item[(a.3)] $\sgn(P')=1 \Rightarrow \mathrm{supp}_{X',(\Delta X)'}= [\max(\pi,X'-(\Delta X)'), \min(3\pi,X+(\Delta X)')]$
	\end{compactitem}
\item[(b)] $\sgn(P)=0 \Rightarrow \mathrm{supp}_{X,\Delta X}= [\max(0,X-\Delta X), \min(2\pi,X+\Delta X)]$
	\begin{compactitem}
	\item[(b.1)] $\sgn(P')=-1 \Rightarrow \mathrm{supp}_{X',(\Delta X)'}= [\max(-\pi,X'-(\Delta X)'), \min(\pi,X+(\Delta X)')]$
	\item[(b.2)] $\sgn(P')=0 \Rightarrow \mathrm{supp}_{X',(\Delta X)'}= [\max(0,X-(\Delta X)'), \min(2\pi,X+(\Delta X)')]$
	\item[(b.3)] $\sgn(P')=1 \Rightarrow \mathrm{supp}_{X',(\Delta X)'}= [\max(\pi,X'-(\Delta X)'), \min(3\pi,X+(\Delta X)')]$
	\end{compactitem}
\item[(c)] $\sgn(P)=1 \Rightarrow \mathrm{supp}_{X,\Delta X}= [\max(\pi,X-\Delta X), \min(3\pi,X+\Delta XA)]$
	\begin{compactitem}
	\item[(c.1)] $\sgn(P')=-1 \Rightarrow \mathrm{supp}_{X',(\Delta X)'}= [\max(-\pi,X'-(\Delta X)'), \min(\pi,X+(\Delta X)')]$
	\item[(c.2)] $\sgn(P')=0 \Rightarrow \mathrm{supp}_{X',(\Delta X)'}= [\max(0,X-(\Delta X)'), \min(2\pi,X+(\Delta X)')]$
	\item[(c.3)] $\sgn(P')=1 \Rightarrow \mathrm{supp}_{X',(\Delta X)'}= [\max(\pi,X'-(\Delta X)'), \min(3\pi,X+(\Delta X)')]$
	\end{compactitem}
\end{compactitem}
If $\sgn(P)=\sgn(P')$ (i.e. the cases $(-1,-1),\ (0,0)$ and $(1,1)$), we obtain the same lower and upper bound as in \cref{def:lbxgen} and \cref{def:ubxgen}, respectively, but with the respective $X_{\mathrm{min}}  \in \{-\pi,0,\pi\}$. If $\sgn(P)=-\sgn(P')$ while $\sgn(P)\not=\sgn(P')$ (i.e. the cases $(1,-1)$ and $(-1,1)$), we shift one of them about $2\pi$ and obtain that $\sgn(P)=\sgn(P')$ (and $P=P'$). Last but not least, we have the cases $(0,1),\ (1,0),\ (0,-1)$ and $(-1,0)$, i.e. where one, let it be $\sgn(P)$, equals zero and the other one, here $\sgn(P')$, is $\pm 1$. If $\sgn(P')=-1$, then we can cut the support $\mathrm{supp}_{X',(\Delta X)'} \subset [-\pi,\pi]$ from the left-hand side at 0 and shift the part that is in $[-\pi,0]$ into $[\pi,2\pi]$. Then, we obtain $\mathrm{supp}_{X',(\Delta X)'}^{\mathrm{shifted}} = [X'-(\Delta X)'+2\pi,\min(X'+(\Delta X)' + 2\pi,2\pi)]\cup[\max(0,X'-(\Delta X)'),X'+(\Delta X)']$. As a consequence, we obtain possibly two lower and upper bounds via 
\begin{align}
[X'-(\Delta X)'+2\pi,\min(X'+(\Delta X)' + 2\pi,2\pi)] \cap \mathrm{supp}_{X,\Delta X}
\intertext{and}
[\max(0,X'-(\Delta X)'),X'+(\Delta X)'] \cap \mathrm{supp}_{X,\Delta X}.
\end{align}
In particular, we have 
\begin{align}
lb_x &= \left( \begin{matrix} 
\max[ \max(0,X-\Delta X), X'-(\Delta X)'+2\pi]\\
\max[ \max(0,X-\Delta X), \max(0,X'-(\Delta X)')]
\end{matrix} \right)
\end{align}
and the upper bound as 
\begin{align}
ub_x &= \left( \begin{matrix}
\min [\min(2\pi,X+\Delta X), \min(X'+(\Delta X)' + 2\pi,2\pi)]\\
\min [\min(2\pi,X+\Delta X), X'+(\Delta X)']
\end{matrix} \right).
\end{align}
If $\sgn(P')=1$, we cut its domain at $2\pi$ and shift the part that is in $[2\pi,3\pi]$ into $[0,\pi]$. Thus, analogously, we obtain $\mathrm{supp}_{X',(\Delta X)'}^{\mathrm{shifted}} = [X'-(\Delta X)',\min(2\pi,X'+(\Delta X)')]\cup[\max(0,X'-(\Delta X)'-2\pi),X'+(\Delta X)'-2\pi]$ and  
\begin{align}
lb_x &= \left( \begin{matrix} 
\max[ \max(0,X-\Delta X), X'-(\Delta X)']\\
\max[ \max(0,X-\Delta X), \max(0,X'-(\Delta X)'-2\pi)]
\end{matrix} \right)
\end{align}
and the upper bound as 
\begin{align}
ub_x = \left( \begin{matrix}
\min [\min(2\pi,X+\Delta X), \min(2\pi,X'+(\Delta X)')]\\
\min [\min(2\pi,X+\Delta X), X'+(\Delta X)'-2\pi]
\end{matrix} \right).
\end{align}
For the cases where $\sgn(P')=0$ and $\sgn(P)=\pm 1$, we obtain analogous solutions (exchange $X$ with $X'$, $\Delta X$ with $(\Delta X)'$ and vice versa). Note again that, in all cases, if $lb_x \geq ub_x$ (componentwise), the intersection is the empty set and the integral vanishes.

We explain the following steps of determining critical points and deriving the remaining integrals only for the case where we have one integration interval $[lb_x,ub_x]$. If we have two, we can execute these steps for both separately and add the integral values to obtain the value of \cref{def:intNNa} with respect to $x=a_2=\lon$.\\
Critical points here are those points between the lower and upper bound(s) where one of the FEHF turn from increase to decrease of the hat (always with respect to a fixed dimension). Here is how to determine them. The critical points can obviously only be from $\{X,\ X\pm\Delta X,\ X',\ X'\pm(\Delta X)'\}.$ Thus, we sort them for increasing order and then check each value whether it is in $[lb_x,ub_x]$. For practical purposes, it is sensible to count how many critical points -- including $lb_x$ and $ub_x$ -- we have. Let this count be $I\leq 6$. In the sequel, we set $p_i$ for a critical point: $lb_x \leq p_i \leq ub_x,\ i=1,...,I,\ I\leq 6$. Note that it cannot be more than 6 different critical points because of the definition of $lb$ and $ub$. The integrals in \cref{def:intNNa} are, thus, equal to 
\begin{align}
\sum_{i=1}^{I-1} \int_{p_i}^{p_{i+1}} \frac{[\Delta X-|x-X|][(\Delta X)'-|x-X'|]}{\Delta X (\Delta X)'}x^q \intd x.
\label{eq:intNNasum}
\end{align}
Note, that if two critical points $p_i$ and $p_{i+1}$ coincide, then the respective integral from $p_i$ to $p_{i+1}$ vanishes. Thus, we do not have to take care of this situation by ourselves. We first consider the following  
\begin{align}
\sum_{i=1}^{I-1} &\int_{p_i}^{p_{i+1}} \frac{[\Delta X-|x-X|][(\Delta X)'-|x-X'|]}{\Delta X (\Delta X)'}x^q \intd x\\
&= 
\sum_{i=1}^{I-1} \int_{p_i}^{p_{i+1}} \frac{[\Delta X-\sgn(x-X)(x-X)][(\Delta X)'-\sgn(x-X')(x-X')]}{\Delta X (\Delta X)'}x^q \intd x\\
&= 
\sum_{i=1}^{I-1}\int_{p_i}^{p_{i+1}} \frac{\left[\Delta X \left[\begin{matrix} +\\-\\-\\+ \end{matrix} \right](x-X)\right]\left[(\Delta X)'\left[\begin{matrix} +\\-\\+\\- \end{matrix}\right](x-X')\right]}{\Delta X (\Delta X)'} x^q \intd x,\\
\end{align}
where the last two cases must coincide. Note that we need to determine the values of $-\sgn(x-X)$ and $-\sgn(x-X')$ for this. For $x=a_1=r$ and $x=a_3=t$, the sign value can be obtained straightforwardly. Similarly, this holds for $x=a_2=\lon$ if $(P,P') \not \in \{(0,1),\ (1,0),\ (0,-1),\ (-1,0)\}$, that is, if we exclude the cases where we shifted the support of one FEHF in order to obtain the upper and lower bound. In the latter cases, we have to shift $P$ or $P'$ here accordingly to the shift done before to obtain the correct sign values. In practice, we can also shift $P$ or $P'$ once in the beginning of the computation of the inner product as well.

For readability, we only consider the integration of $x^qW\left(x;X,\Delta X,X',(\Delta X)'\right)$ defined in
\begin{align}
\int \frac{\left[\Delta X \left[\begin{matrix} +\\-\\-\end{matrix} \right](x-X)\right]\left[(\Delta X)'\left[\begin{matrix} +\\-\\+ \end{matrix}\right](x-X')\right]}{\Delta X (\Delta X)'}x^q\intd x
\eqqcolon 
&\int x^q\frac{W\left(x;X,\Delta X,X',(\Delta X)'\right)}{\Delta X (\Delta X)'}\intd x\\
= \frac{1}{\Delta X (\Delta X)'}
&\int x^qW\left(x;X,\Delta X,X',(\Delta X)'\right)\intd x
\label{eq:intNNafinal}
\end{align}
in the sequel. We obtain
\begin{align}
&\int x^qW\left(x;X,\Delta X,X',(\Delta X)'\right)\intd x\\ 
&= \int
\left[\begin{matrix} +\\+\\-\end{matrix} \right]
x^{q+2}
\left[\begin{matrix} +\\-\\+\end{matrix} \right]
x^{q+1}
\left(\left[\begin{matrix} -\\+\\+\end{matrix} \right] X\left[\begin{matrix} -\\+\\+\end{matrix} \right] X' \left[\begin{matrix} +\\+\\+\end{matrix} \right] \Delta X \left[\begin{matrix} +\\+\\-\end{matrix} \right](\Delta X)' \right)\\
&\qquad \qquad
\left[\begin{matrix} +\\+\\-\end{matrix} \right]
x^q
\left(X \left[\begin{matrix} -\\+\\+\end{matrix} \right]\Delta X \right) \left( X' \left[\begin{matrix} -\\+\\-\end{matrix} \right](\Delta X') \right)
\intd x\\ \\
&= 
\left[\begin{matrix} +\\+\\-\end{matrix} \right]
\frac{x^{q+3}}{q+3}
\left[\begin{matrix} +\\-\\+\end{matrix} \right]
\frac{x^{q+2}}{q+2} 
\left(\left[\begin{matrix} -\\+\\+\end{matrix} \right] X\left[\begin{matrix} -\\+\\+\end{matrix} \right] X' \left[\begin{matrix} +\\+\\+\end{matrix} \right] \Delta X \left[\begin{matrix} +\\+\\-\end{matrix} \right](\Delta X)' \right)\\
&\qquad \qquad
\left[\begin{matrix} +\\+\\-\end{matrix} \right]
\frac{x^{q+1}}{q+1}
\left(X \left[\begin{matrix} -\\+\\+\end{matrix} \right]\Delta X \right) \left( X' \left[\begin{matrix} -\\+\\-\end{matrix} \right](\Delta X') \right)\\
&= \left(
\left[\begin{matrix} +\\+\\-\end{matrix} \right]
(q+2)(q+1)x^{q+3}
\left[\begin{matrix} +\\-\\+\end{matrix} \right]
(q+3)(q+1)x^{q+2} 
\left(\left[\begin{matrix} -\\+\\+\end{matrix} \right] X\left[\begin{matrix} -\\+\\+\end{matrix} \right] X' \left[\begin{matrix} +\\+\\+\end{matrix} \right] \Delta X \left[\begin{matrix} +\\+\\-\end{matrix} \right](\Delta X)' \right)\right.\\
&\qquad \qquad \left.
\left[\begin{matrix} +\\+\\-\end{matrix} \right]
(q+3)(q+2)x^{q+1}
\left(X \left[\begin{matrix} -\\+\\+\end{matrix} \right]\Delta X \right)\left( X' \left[\begin{matrix} -\\+\\-\end{matrix} \right](\Delta X') \right)\right)\frac{1}{(q+3)(q+2)(q+1)}\\
&= \left(
\left[\begin{matrix} +\\+\\-\end{matrix} \right]
(q+2)(q+1)x^{2}
\left[\begin{matrix} +\\-\\+\end{matrix} \right]
(q+3)(q+1)x 
\left(\left[\begin{matrix} -\\+\\+\end{matrix} \right] X\left[\begin{matrix} -\\+\\+\end{matrix} \right] X' \left[\begin{matrix} +\\+\\+\end{matrix} \right] \Delta X \left[\begin{matrix} +\\+\\-\end{matrix} \right](\Delta X)' \right)\right.\\
&\qquad \qquad \left.
\left[\begin{matrix} +\\+\\-\end{matrix} \right]
(q+3)(q+2)
\left(X \left[\begin{matrix} -\\+\\+\end{matrix} \right]\Delta X \right) \left( X' \left[\begin{matrix} -\\+\\-\end{matrix} \right](\Delta X') \right)\right)\frac{x^{q+1}}{(q+3)(q+2)(q+1)}.\\
\end{align}
With these values the $\Lp{2}$ inner product of two FEHFs is fully discussed. For the respective $\lp{2}$ inner product, we obtain
\begin{landscape}
\begin{align}
&\left\langle \nabla_{r\xi(\lon,t)} N_{(R,\Phi,T),(\Delta R,\Delta \Phi,\Delta T)}, \nabla_{r\xi(\lon,t)} N_{(R',\Phi',T'),((\Delta R)',(\Delta \Phi)',(\Delta T)')}\right\rangle_{\lp{2}}\\
%
%
&=\int_{\ball} \chi_{\mathrm{supp}_{[(R,\Phi,T)-(\Delta R,\Delta \Phi,\Delta T),(R,\Phi,T)+(\Delta R,\Delta \Phi,\Delta T)]}}(r,\lon,t)\\
&\qquad\qquad\times 
\chi_{\mathrm{supp}_{[(R',\Phi',T')-((\Delta R)',(\Delta \Phi)',(\Delta T)'),(R',\Phi',T')+((\Delta R)',(\Delta \Phi)',(\Delta T)')]}}(r,\lon,t)\\
&\qquad\qquad\qquad\times \left(
 \frac{[-\sgn(r-R)]}{\Delta R} \frac{\Delta \Phi-|\lon-\Phi|}{\Delta \Phi}\frac{\Delta T-|t-T|}{\Delta T}\frac{[-\sgn(r-R')]}{(\Delta R)'} \frac{(\Delta \Phi)'-|\lon-\Phi'|}{(\Delta \Phi)'}\frac{(\Delta T)'-|t-T'|}{(\Delta T)'}\right.\\ 
&\qquad\qquad\qquad\qquad + 
\frac{1}{r^2} \frac{1}{1-t^2} \frac{\Delta R-|r-R|}{\Delta R} \frac{[-\sgn(\lon-\Phi)]}{\Delta \Phi} \frac{\Delta T-|t-T|}{\Delta T} \frac{(\Delta R)'-|r-R'|}{(\Delta R)'} \frac{[-\sgn(\lon-\Phi')]}{(\Delta \Phi)'} \frac{(\Delta T)'-|t-T'|}{(\Delta T)'}\\
&\qquad\qquad\qquad\qquad
\left. + \frac{1}{r^2} (1-t^2)\frac{\Delta R-|r-R|}{\Delta R}\frac{\Delta \Phi-|\lon-\Phi|}{\Delta \Phi}\frac{[-\sgn(t-T)]}{\Delta T}\frac{(\Delta R)'-|r-R'|}{(\Delta R)'}\frac{(\Delta \Phi)'-|\lon-\Phi'|}{(\Delta \Phi)'}\frac{[-\sgn(t-T')]}{(\Delta T)'}
\right)\intd x(r,\lon,t)\\
&=\int_{\ball} \chi_{\mathrm{supp}_{[(R,\Phi,T)-(\Delta R,\Delta \Phi,\Delta T),(R,\Phi,T)+(\Delta R,\Delta \Phi,\Delta T)]}}(r,\lon,t)\\
&\qquad\qquad\times 
\chi_{\mathrm{supp}_{[(R',\Phi',T')-((\Delta R)',(\Delta \Phi)',(\Delta T)'),(R',\Phi',T')+((\Delta R)',(\Delta \Phi)',(\Delta T)')]}}(r,\lon,t)\\
&\qquad\qquad\qquad\times 
\frac{[-\sgn(r-R)]}{\Delta R} \frac{\Delta \Phi-|\lon-\Phi|}{\Delta \Phi}\frac{\Delta T-|t-T|}{\Delta T}\frac{[-\sgn(r-R')]}{(\Delta R)'} \frac{(\Delta \Phi)'-|\lon-\Phi'|}{(\Delta \Phi)'}\frac{(\Delta T)'-|t-T'|}{(\Delta T)'}\intd x(r,\lon,t)\\ 
&\qquad + \int_{\ball} \chi_{\mathrm{supp}_{[(R,\Phi,T)-(\Delta R,\Delta \Phi,\Delta T),(R,\Phi,T)+(\Delta R,\Delta \Phi,\Delta T)]}}(r,\lon,t)\\
&\qquad\qquad\qquad\times 
\chi_{\mathrm{supp}_{[(R',\Phi',T')-((\Delta R)',(\Delta \Phi)',(\Delta T)'),(R',\Phi',T')+((\Delta R)',(\Delta \Phi)',(\Delta T)')]}}(r,\lon,t)\\
&\qquad\qquad\qquad\qquad\times 
\frac{1}{r^2} \frac{1}{1-t^2} \frac{\Delta R-|r-R|}{\Delta R} \frac{[-\sgn(\lon-\Phi)]}{\Delta \Phi} \frac{\Delta T-|t-T|}{\Delta T} \frac{(\Delta R)'-|r-R'|}{(\Delta R)'} \frac{[-\sgn(\lon-\Phi')]}{(\Delta \Phi)'} \frac{(\Delta T)'-|t-T'|}{(\Delta T)'}\intd x(r,\lon,t)\\ 
&\qquad + \int_{\ball} \chi_{\mathrm{supp}_{[(R,\Phi,T)-(\Delta R,\Delta \Phi,\Delta T),(R,\Phi,T)+(\Delta R,\Delta \Phi,\Delta T)]}}(r,\lon,t)\\
&\qquad\qquad\qquad\times 
\chi_{\mathrm{supp}_{[(R',\Phi',T')-((\Delta R)',(\Delta \Phi)',(\Delta T)'),(R',\Phi',T')+((\Delta R)',(\Delta \Phi)',(\Delta T)')]}}(r,\lon,t)\\
&\qquad\qquad\qquad\qquad\times \frac{1}{r^2} (1-t^2)\frac{\Delta R-|r-R|}{\Delta R}\frac{\Delta \Phi-|\lon-\Phi|}{\Delta \Phi}\frac{[-\sgn(t-T)]}{\Delta T}\frac{(\Delta R)'-|r-R'|}{(\Delta R)'}\frac{(\Delta \Phi)'-|\lon-\Phi'|}{(\Delta \Phi)'}\frac{[-\sgn(t-T')]}{(\Delta T)'}\intd x(r,\lon,t)
\end{align}
\begin{align}
&=\int_{lb_r}^{ub_r} \int_{lb_\lon}^{ub_\lon} \int_{lb_t}^{ub_t}
\frac{[-\sgn(r-R)]}{\Delta R} \frac{\Delta \Phi-|\lon-\Phi|}{\Delta \Phi}\frac{\Delta T-|t-T|}{\Delta T}\frac{[-\sgn(r-R')]}{(\Delta R)'} \frac{(\Delta \Phi)'-|\lon-\Phi'|}{(\Delta \Phi)'}\frac{(\Delta T)'-|t-T'|}{(\Delta T)'} r^2 \intd r \intd \lon \intd t\\ 
&\quad + \int_{lb_r}^{ub_r} \int_{lb_\lon}^{ub_\lon} \int_{lb_t}^{ub_t}
\frac{1}{r^2} \frac{1}{1-t^2} \frac{\Delta R-|r-R|}{\Delta R} \frac{[-\sgn(\lon-\Phi)]}{\Delta \Phi} \frac{\Delta T-|t-T|}{\Delta T} \frac{(\Delta R)'-|r-R'|}{(\Delta R)'} \frac{[-\sgn(\lon-\Phi')]}{(\Delta \Phi)'} \frac{(\Delta T)'-|t-T'|}{(\Delta T)'}r^2 \intd r \intd \lon \intd t\\
&\quad + \int_{lb_r}^{ub_r} \int_{lb_\lon}^{ub_\lon} \int_{lb_t}^{ub_t}
\frac{1}{r^2} (1-t^2)\frac{\Delta R-|r-R|}{\Delta R}\frac{\Delta \Phi-|\lon-\Phi|}{\Delta \Phi}\frac{[-\sgn(t-T)]}{\Delta T}\frac{(\Delta R)'-|r-R'|}{(\Delta R)'}\frac{(\Delta \Phi)'-|\lon-\Phi'|}{(\Delta \Phi)'}\frac{[-\sgn(t-T')]}{(\Delta T)'}r^2 \intd r \intd \lon \intd t\\
&=\int_{lb_r}^{ub_r} 
\frac{\sgn(r-R)}{\Delta R} 
\frac{\sgn(r-R')}{(\Delta R)'} 
r^2 
\intd r 
\int_{lb_\lon}^{ub_\lon} 
\frac{\Delta \Phi-|\lon-\Phi|}{\Delta \Phi}
\frac{(\Delta \Phi)'-|\lon-\Phi'|}{(\Delta \Phi)'}
\intd \lon 
\int_{lb_t}^{ub_t}
\frac{\Delta T-|t-T|}{\Delta T}
\frac{(\Delta T)'-|t-T'|}{(\Delta T)'} 
\intd t\\ 
&\quad + 
\int_{lb_r}^{ub_r} 
\frac{\Delta R-|r-R|}{\Delta R} 
\frac{(\Delta R)'-|r-R'|}{(\Delta R)'} 
\intd r 
\int_{lb_\lon}^{ub_\lon} 
\frac{\sgn(\lon-\Phi)}{\Delta \Phi} 
\frac{\sgn(\lon-\Phi')}{(\Delta \Phi)'} 
\intd \lon 
\int_{lb_t}^{ub_t}
\frac{1}{1-t^2} 
\frac{\Delta T-|t-T|}{\Delta T} 
\frac{(\Delta T)'-|t-T'|}{(\Delta T)'}
\intd t\\
&\quad + 
\int_{lb_r}^{ub_r} 
\frac{\Delta R-|r-R|}{\Delta R}
\frac{(\Delta R)'-|r-R'|}{(\Delta R)'}
\intd r 
\int_{lb_\lon}^{ub_\lon} 
\frac{\Delta \Phi-|\lon-\Phi|}{\Delta \Phi}
\frac{(\Delta \Phi)'-|\lon-\Phi'|}{(\Delta \Phi)'}
\intd \lon 
\int_{lb_t}^{ub_t}
(1-t^2)
\frac{\sgn(t-T)}{\Delta T}
\frac{\sgn(t-T')}{(\Delta T)'}
\intd t\\
&= \sum_{k=1}^3\int_{lb_{a_k}}^{ub_{a_k}} \frac{\sgn(a_k-A_k)}{\Delta A_k} \frac{\sgn(a_k-A_k')}{(\Delta A_k)'} \left\{\begin{matrix} a_k^2, &k=1,\\ 1,& k=2,\\ 1-a_k^2, & k=3 \end{matrix} \right\} \intd a_k\\ &\qquad \times \prod_{j=1,\ j\not=k}^3 \int_{lb_{a_j}}^{ub_{a_j}} \frac{\Delta A_j - |a_j-A_j|}{\Delta A_j} \frac{(\Delta A_j)' - |a_j-A_j'|}{(\Delta A_j)'} \left\{\begin{matrix} \frac{1}{1-a_j^2}, &j=3,k=2,\\ 1,& \text{else} \end{matrix} \right\} \intd a_j
\end{align}
\end{landscape}
Note that $lb_r = lb_{a_1}, ub_r = ub_{a_1}, lb_\lon = lb_{a_2}, ub_\lon = ub_{a_2}, lb_t = lb_{a_3}$ and $ub_t = ub_{a_3}$. We see that also this integral fragments with respect to the variables $a_j,\ j=1,2,3,$ and the only integrals left for a discussion are of the form 
\begin{align}
\int \frac{\sgn(x-X)\sgn(x-X')}{\Delta X (\Delta X)'} x^q \intd x,
\end{align}
again for $q\in\nat_0$, here in particular $q=0$ and $q=2$, and 
\begin{align}
\int \frac{1}{1-t^2} 
\frac{\Delta T-|t-T|}{\Delta T} 
\frac{(\Delta T)'-|t-T'|}{(\Delta T)'} \intd t.
\end{align}
For each integral between two critical points, we obtain
\begin{align}
\int \frac{\sgn(x-X)\sgn(x-X')}{\Delta X (\Delta X)'} x^q \intd x
= \pm \frac{1}{\Delta X (\Delta X)'}\int x^q \intd x
= \pm \frac{x^{q+1}}{(q+1)\Delta X (\Delta X)'} 
\end{align}
in the first case and 
\begin{align}
&\int \frac{W\left(t;T,\Delta T,T',(\Delta T)'\right)}{1-t^2}\intd t\\ 
&= \int
\left[\begin{matrix} +\\+\\-\end{matrix} \right]
\frac{t^{2}}{1-t^2}
\left[\begin{matrix} +\\-\\+\end{matrix} \right]
\frac{t}{1-t^2}
\left(\left[\begin{matrix} -\\+\\+\end{matrix} \right] T\left[\begin{matrix} -\\+\\+\end{matrix} \right] T' \left[\begin{matrix} +\\+\\+\end{matrix} \right] \Delta T \left[\begin{matrix} +\\+\\-\end{matrix} \right](\Delta T)' \right)\\
&\qquad \qquad
\left[\begin{matrix} +\\+\\-\end{matrix} \right]
\frac{1}{1-t^2}
\left(T \left[\begin{matrix} -\\+\\+\end{matrix} \right]\Delta T \right) \left( T' \left[\begin{matrix} -\\+\\-\end{matrix} \right](\Delta T') \right)
\intd t\\ \\
&= 
\left[\begin{matrix} +\\+\\-\end{matrix} \right]
(\atanh(t)-t)
\left[\begin{matrix} -\\+\\-\end{matrix} \right]
\frac{\log\left(1-t^2\right)}{2}
\left(\left[\begin{matrix} -\\+\\+\end{matrix} \right] T\left[\begin{matrix} -\\+\\+\end{matrix} \right] T' \left[\begin{matrix} +\\+\\+\end{matrix} \right] \Delta T \left[\begin{matrix} +\\+\\-\end{matrix} \right](\Delta T)' \right)\\
&\qquad \qquad
\left[\begin{matrix} +\\+\\-\end{matrix} \right]
\atanh(t)
\left(T \left[\begin{matrix} -\\+\\+\end{matrix} \right]\Delta T \right) \left( T' \left[\begin{matrix} -\\+\\-\end{matrix} \right](\Delta T') \right)\\
\end{align}
in the second case. Note that the latter is obviously not well-defined for $t=\pm 1$ due to the logarithm which is why we have to diminish the domain for this variable in practice.

At last, we consider the mixed cases. We obtain
\begin{align}
&\left\langle N_{A,\Delta A}, G_{m,n,j}^{\mathrm{I}} \right\rangle_{\Lp{2}}\\
&\qquad =\prod_{j=1}^3 \int_{\max(A_{\mathrm{min}},A_j-\Delta A_j)}^{\min(A_{\mathrm{max}},A_j+\Delta A_j)} \frac{\Delta A_j-|a_j-A_j|}{\Delta A_j} G_{m,n,j}^{\mathrm{I}}(x(a_1,a_2,a_3))\left\{\begin{matrix}a^2_1,&j=1\\1,&j=2,3\end{matrix}\right\} \intd a_j\\
&\qquad = p_{m,n}q_{n,j}
\int_{\max(R_{\mathrm{min}},R-\Delta R)}^{\min(R_{\mathrm{max}},R+\Delta R)} \frac{\Delta R-|r-R|}{\Delta R} P_m^{(0,n+1/2)}(I(r))\left(\frac{r}{\mathbf{R}}\right)^nr^2 \intd r\\
&\qquad\qquad  \times \int_{\Phi-\Delta \Phi}^{\Phi+\Delta \Phi} \frac{\Delta \Phi-|\lon-\Phi|}{\Delta \Phi} \mathrm{Trig}(j\lon) \intd \lon\\
&\qquad\qquad  \times \int_{\max(T_{\mathrm{min}},T-\Delta T)}^{\min(T_{\mathrm{max}},T+\Delta T)} \frac{\Delta T-|t-T|}{\Delta T} P_{n,|j|}(t) \intd t.
\end{align} 
We immediately see that the integrals with respect to $r$ and $t$ cannot be calculated analytically. However, as they are one-dimensional integrals, they can easily be integrated numerically, e.g. with suitable software libraries. With respect to the $\lon$-integral, we obtain 
\begin{align}
\frac{\Delta \Phi-|\lon-\Phi|}{\Delta \Phi} \mathrm{Trig}(j\lon)
&=\frac{\Delta \Phi-\sgn(\lon-\Phi)(\lon-\Phi)}{\Delta \Phi} \mathrm{Trig}(j\lon)\\
&= \left(1 +\frac{\Phi \sgn(\lon-\Phi)}{\Delta \Phi}\right)\mathrm{Trig}(j\lon) - \frac{\sgn(\lon-\Phi)}{\Delta \Phi}\ \lon \mathrm{Trig}(j\lon). 
\end{align}
Thus, the following cases remain:
\begin{align}
I_1(j,\lon) \coloneqq \int \lon\mathrm{Trig}(j\lon) \intd \lon
&= 
\int\left\{ \begin{matrix} \sqrt{2}\lon \cos(j\lon),& j<0\\ \lon,&j=0\\\sqrt{2}\lon\sin(j\lon),&j>0 \end{matrix} \right\} \intd \lon 
= 
\left\{ \begin{matrix} 
\sqrt{2} \left[ \frac{\cos(j\lon)}{j^2} + \frac{\lon\sin(j\lon)}{j} \right],& j<0\\ 
\frac{1}{2} \lon^2,&j=0\\
\sqrt{2}\left[ \frac{\sin(j\lon)}{j^2} - \frac{\lon\cos(j\lon)}{j}\right],&j>0 
\end{matrix} \right\}   
\intertext{and}
I_2(j,\lon) \coloneqq \int\mathrm{Trig}(j\lon) \intd \lon
&= 
\int\left\{ \begin{matrix} \sqrt{2} \cos(j\lon),& j<0\\ 1,&j=0\\\sqrt{2}\sin(j\lon),&j>0 \end{matrix} \right\} \intd \lon 
= 
\left\{ \begin{matrix} \frac{\sqrt{2}}{j} \sin(j\lon),& j<0\\ \lon,&j=0\\ -\frac{\sqrt{2}}{j}\cos(j\lon),&j>0 \end{matrix} \right\}.
\end{align}
This yields
\begin{align}
&\int_{\Phi-\Delta\Phi}^{\Phi+\Delta\Phi}\frac{\Delta \Phi-|\lon-\Phi|}{\Delta \Phi} \mathrm{Trig}(j\lon)\intd \lon\\
&= \int_{\Phi-\Delta\Phi}^{\Phi}\frac{\Delta \Phi-|\lon-\Phi|}{\Delta \Phi} \mathrm{Trig}(j\lon)\intd \lon + \int_{\Phi}^{\Phi+\Delta\Phi}\frac{\Delta \Phi-|\lon-\Phi|}{\Delta \Phi} \mathrm{Trig}(j\lon)\intd \lon\\
&= \int_{\Phi-\Delta\Phi}^{\Phi}\frac{\Delta \Phi + \lon-\Phi}{\Delta \Phi} \mathrm{Trig}(j\lon)\intd \lon 
+ \int_{\Phi}^{\Phi+\Delta\Phi}\frac{\Delta \Phi-\lon+\Phi}{\Delta \Phi} \mathrm{Trig}(j\lon)\intd\lon\\
&= \int_{\Phi-\Delta\Phi}^{\Phi}\left(1 - \frac{\Phi}{\Delta \Phi} + \frac{\lon}{\Delta \Phi}\right)\mathrm{Trig}(j\lon)\intd \lon 
+ \int_{\Phi}^{\Phi+\Delta\Phi}\left(1 + \frac{\Phi}{\Delta \Phi} - \frac{\lon}{\Delta \Phi}\right) \mathrm{Trig}(j\lon)\intd\lon\\
&= \frac{1}{\Delta\Phi}I_1(j,\lon) |_{\Phi-\Delta\Phi}^{\Phi} + \left(1-\frac{\Phi}{\Delta\Phi}\right)I_2(j,\lon) |_{\Phi-\Delta\Phi}^{\Phi} - \frac{1}{\Delta\Phi}I_1(j,\lon) |_{\Phi}^{\Phi+\Delta\Phi} + \left(1+\frac{\Phi}{\Delta\Phi}\right)I_2(j,\lon) |_{\Phi}^{\Phi+\Delta\Phi}\\
&= 
-\frac{1}{\Delta\Phi}I_1(j,\Phi-\Delta\Phi) 
- \left(1-\frac{\Phi}{\Delta\Phi}\right)I_2(j,\Phi-\Delta\Phi)\\
&\qquad 
+ \frac{1}{\Delta\Phi}I_1(j,\Phi) 
+ \left(1-\frac{\Phi}{\Delta\Phi}\right)I_2(j,\Phi)
+ \frac{1}{\Delta\Phi}I_1(j,\Phi) 
- \left(1+\frac{\Phi}{\Delta\Phi}\right)I_2(j,\Phi)\\
&\qquad 
- \frac{1}{\Delta\Phi}I_1(j,\Phi+\Delta\Phi) 
+ \left(1+\frac{\Phi}{\Delta\Phi}\right)I_2(j,\Phi+\Delta\Phi)\\
&= 
-\frac{1}{\Delta\Phi}I_1(j,\Phi-\Delta\Phi) 
- \left(1-\frac{\Phi}{\Delta\Phi}\right)I_2(j,\Phi-\Delta\Phi)
+ \frac{2}{\Delta\Phi}I_1(j,\Phi) 
- \frac{2\Phi}{\Delta\Phi}I_2(j,\Phi)\\
&\qquad 
- \frac{1}{\Delta\Phi}I_1(j,\Phi+\Delta\Phi) 
+ \left(1+\frac{\Phi}{\Delta\Phi}\right)I_2(j,\Phi+\Delta\Phi).
\end{align}
For the gradients, we have similarly
\begin{landscape}
\begin{align}
&\left\langle \nabla_{r\xi(\lon,t)} N_{(R,\Phi,T),(\Delta R,\Delta \Phi,\Delta T)}, \nabla_{r\xi(\lon,t)} G_{m,n,j}^{\mathrm{I}} \right\rangle_{\lp{2}}\\
&= \int_{\ball} \chi_{\mathrm{supp}_{[(R,\Phi,T)-(\Delta R,\Delta \Phi,\Delta T),(R,\Phi,T)+(\Delta R,\Delta \Phi,\Delta T)]}}(r,\lon,t) \left(
\era \frac{[-\sgn(r-R)]}{\Delta R} \frac{\Delta \Phi-|\lon-\Phi|}{\Delta \Phi}\frac{\Delta T-|t-T|}{\Delta T} \right.\\
&\qquad\qquad + \left. \frac{1}{r}\ephi \frac{1}{\sqrt{1-t^2}} \frac{\Delta R-|r-R|}{\Delta R} \frac{[-\sgn(\lon-\Phi)]}{\Delta \Phi} \frac{\Delta T-|t-T|}{\Delta T}
+ \frac{1}{r}\ete \sqrt{1-t^2}\frac{\Delta R-|r-R|}{\Delta R}\frac{\Delta \Phi-|\lon-\Phi|}{\Delta \Phi}\frac{[-\sgn(t-T)]}{\Delta T}
\right)\\
&\qquad\qquad\cdot\left( p_{m,n}q_{n,j} \sum_{p=1}^4 G_{m,n,j;p}^{\mathrm{I}} (r\xi(\lon,t))\right)\intd x(r,\lon,t) \\
&= p_{m,n}q_{n,j} \int_{\ball} \chi_{\mathrm{supp}_{[(R,\Phi,T)-(\Delta R,\Delta \Phi,\Delta T),(R,\Phi,T)+(\Delta R,\Delta \Phi,\Delta T)]}}(r,\lon,t) \\
&\qquad\qquad\qquad \left[ \left(
\era \frac{[-\sgn(r-R)]}{\Delta R} \frac{\Delta \Phi-|\lon-\Phi|}{\Delta \Phi}\frac{\Delta T-|t-T|}{\Delta T}\right) \cdot\left( \sum_{p=1}^2G_{m,n,j;p}^{\mathrm{I}} (r\xi(\lon,t)) \right) \right.\\
&\qquad\qquad\qquad\qquad  + \left( \frac{1}{r}\ephi \frac{1}{\sqrt{1-t^2}} \frac{\Delta R-|r-R|}{\Delta R} \frac{[-\sgn(\lon-\Phi)]}{\Delta \Phi} \frac{\Delta T-|t-T|}{\Delta T}\right) \cdot\left(G_{m,n,j;3}^{\mathrm{I}} (r\xi(\lon,t))\right)\\
&\qquad\qquad\qquad\qquad  + \left. \left( \frac{1}{r}\ete \sqrt{1-t^2}\frac{\Delta R-|r-R|}{\Delta R}\frac{\Delta \Phi-|\lon-\Phi|}{\Delta \Phi}\frac{[-\sgn(t-T)]}{\Delta T}
\right) \cdot\left(G_{m,n,j;4}^{\mathrm{I}} (r\xi(\lon,t)) \right)\right]\intd x(r,\lon,t)\\
&= p_{m,n}q_{n,j} \int_{R-\Delta R}^{R+\Delta R}
\int_{\Phi-\Delta \Phi}^{\Phi+\Delta \Phi}
\int_{T-\Delta T}^{T+\Delta T}
\left(
\era \frac{[-\sgn(r-R)]}{\Delta R} \frac{\Delta \Phi-|\lon-\Phi|}{\Delta \Phi}\frac{\Delta T-|t-T|}{\Delta T}\right) \cdot\left( \sum_{p=1}^2 G_{m,n,j;p}^{\mathrm{I}} (r\xi(\lon,t))\right)r^2 \intd r \intd \lon \intd t\\ 
&\qquad + p_{m,n}q_{n,j} \int_{R-\Delta R}^{R+\Delta R}
\int_{\Phi-\Delta \Phi}^{\Phi+\Delta \Phi}
\int_{T-\Delta T}^{T+\Delta T} 
\left( \frac{1}{r}\ephi \frac{1}{\sqrt{1-t^2}} \frac{\Delta R-|r-R|}{\Delta R} \frac{[-\sgn(\lon-\Phi)]}{\Delta \Phi} \frac{\Delta T-|t-T|}{\Delta T}\right) \cdot\left(G_{m,n,j;3}^{\mathrm{I}} (r\xi(\lon,t)) \right)r^2 \intd r \intd \lon \intd t\\
&\qquad + p_{m,n}q_{n,j} \int_{R-\Delta R}^{R+\Delta R}
\int_{\Phi-\Delta \Phi}^{\Phi+\Delta \Phi}
\int_{T-\Delta T}^{T+\Delta T} 
\left( \frac{1}{r}\ete \sqrt{1-t^2}\frac{\Delta R-|r-R|}{\Delta R}\frac{\Delta \Phi-|\lon-\Phi|}{\Delta \Phi}\frac{[-\sgn(t-T)]}{\Delta T}
\right)\cdot\left(G_{m,n,j;4}^{\mathrm{I}} (r\xi(\lon,t)) \right) r^2 \intd r \intd \lon \intd t\\
\end{align}
\begin{align}
&= p_{m,n}q_{n,j} \int_{R-\Delta R}^{R+\Delta R}
\int_{\Phi-\Delta \Phi}^{\Phi+\Delta \Phi}
\int_{T-\Delta T}^{T+\Delta T}
\left(
\era \frac{[-\sgn(r-R)]}{\Delta R} \frac{\Delta \Phi-|\lon-\Phi|}{\Delta \Phi}\frac{\Delta T-|t-T|}{\Delta T}\right) \\
&\qquad\qquad\qquad\qquad \cdot\left( \left(P_m^{(0,n+1/2)}\left(I(r)\right)\right)'I'(r) \left(\frac{r}{\mathbf{R}}\right)^{n} P_{n,|j|}(t)\mathrm{Trig}(j\lon) \xi(\lon,t) + \frac{n}{\mathbf{R}} P_m^{(0,n+1/2)}\left(I(r)\right) \left(\frac{r}{\mathbf{R}}\right)^{n-1} P_{n,|j|}(t)\mathrm{Trig}(j\lon) \xi(\lon,t)\right)r^2 \intd r \intd \lon \intd t\\ 
&\qquad + p_{m,n}q_{n,j} \int_{R-\Delta R}^{R+\Delta R}
\int_{\Phi-\Delta \Phi}^{\Phi+\Delta \Phi}
\int_{T-\Delta T}^{T+\Delta T} 
\left( \frac{1}{r}\ephi \frac{1}{\sqrt{1-t^2}} \frac{\Delta R-|r-R|}{\Delta R} \frac{[-\sgn(\lon-\Phi)]}{\Delta \Phi} \frac{\Delta T-|t-T|}{\Delta T}\right)\\
&\qquad\qquad\qquad\qquad \cdot\left( \frac{j}{\mathbf{R}} P_m^{(0,n+1/2)}\left(I(r)\right) \left(\frac{r}{\mathbf{R}}\right)^{n-1} \frac{1}{\sqrt{1-t^2}} P_{n,|j|}(t) \mathrm{Trig}(-j\lon)\ephi(\lon,t)\right)r^2 \intd r \intd \lon \intd t\\
&\qquad + p_{m,n}q_{n,j} \int_{R-\Delta R}^{R+\Delta R}
\int_{\Phi-\Delta \Phi}^{\Phi+\Delta \Phi}
\int_{T-\Delta T}^{T+\Delta T} 
\left( \frac{1}{r}\ete \sqrt{1-t^2}\frac{\Delta R-|r-R|}{\Delta R}\frac{\Delta \Phi-|\lon-\Phi|}{\Delta \Phi}\frac{[-\sgn(t-T)]}{\Delta T}
\right)\\
&\qquad\qquad\qquad\qquad \cdot\left( \frac{1}{\mathbf{R}} P_m^{(0,n+1/2)}\left(I(r)\right) \left(\frac{r}{\mathbf{R}}\right)^{n-1} \sqrt{1-t^2} P'_{n,|j|}(t)\mathrm{Trig}(j\lon)\ete(\lon,t)\right) r^2 \intd r \intd \lon \intd t\\
&= p_{m,n}q_{n,j} 
\int_{R-\Delta R}^{R+\Delta R}
\frac{[-\sgn(r-R)]}{\Delta R} \left(P_m^{(0,n+1/2)}\left(I(r)\right)\right)'I'(r) \left(\frac{r}{\mathbf{R}}\right)^{n} r^2 \intd r 
\int_{\Phi-\Delta \Phi}^{\Phi+\Delta \Phi}
\frac{\Delta \Phi-|\lon-\Phi|}{\Delta \Phi} \mathrm{Trig}(j\lon)
\intd \lon  
\int_{T-\Delta T}^{T+\Delta T}
\frac{\Delta T-|t-T|}{\Delta T}  
P_{n,|j|}(t)\intd t\\ 
&\qquad + p_{m,n}q_{n,j}n 
\int_{R-\Delta R}^{R+\Delta R}
\frac{[-\sgn(r-R)]}{\Delta R} r P_m^{(0,n+1/2)}\left(I(r)\right) \left(\frac{r}{\mathbf{R}}\right)^{n}\intd r
\int_{\Phi-\Delta \Phi}^{\Phi+\Delta \Phi}
\frac{\Delta \Phi-|\lon-\Phi|}{\Delta \Phi}\mathrm{Trig}(j\lon) \intd \lon
\int_{T-\Delta T}^{T+\Delta T}
\frac{\Delta T-|t-T|}{\Delta T} P_{n,|j|}(t) \intd t\\
&\qquad + p_{m,n}q_{n,j} \int_{R-\Delta R}^{R+\Delta R}\frac{\Delta R-|r-R|}{\Delta R}P_m^{(0,n+1/2)}\left(I(r)\right) \left(\frac{r}{\mathbf{R}}\right)^{n} \intd r 
\int_{\Phi-\Delta \Phi}^{\Phi+\Delta \Phi}\frac{[-\sgn(\lon-\Phi)]}{\Delta \Phi} j \mathrm{Trig}(-j\lon) \intd \lon 
\int_{T-\Delta T}^{T+\Delta T} 
\frac{\Delta T-|t-T|}{\Delta T} \frac{1}{1-t^2} P_{n,|j|}(t)
\intd t\\
&\qquad + p_{m,n}q_{n,j}
\int_{R-\Delta R}^{R+\Delta R}
\frac{\Delta R-|r-R|}{\Delta R}
P_m^{(0,n+1/2)}\left(I(r)\right) \left(\frac{r}{\mathbf{R}}\right)^{n} \intd r 
\int_{\Phi-\Delta \Phi}^{\Phi+\Delta \Phi}
\frac{\Delta \Phi-|\lon-\Phi|}{\Delta \Phi}\mathrm{Trig}(j\lon) \intd \lon
\int_{T-\Delta T}^{T+\Delta T} 
\frac{[-\sgn(t-T)]}{\Delta T}(1-t^2) P'_{n,|j|}(t)\intd t\\
&= p_{n,m}q_{n,j}
\sum_{k=1}^3 
\int_{A_k-\Delta A_k}^{A_k+\Delta A_k}
\frac{-\sgn(a_k-A_k)}{\Delta A_k}
\left\{\begin{matrix} 
\left(P_m^{(0,n+1/2)}\left(I(a_k)\right)\right)'I'(a_k) \left(\frac{a_k}{\mathbf{R}}\right)^{n} a_k^2 + 
P_m^{(0,n+1/2)}\left(I(a_k)\right) \left(\frac{a_k}{\mathbf{R}}\right)^{n} n a_k,  &k=1\\ 
j\mathrm{Trig}(-j a_k),&k=2\\ 
(1-a_k^2)P_{n,|j|}'(a_k),&k=3 
\end{matrix}\right\}
\intd a_k\\
&\qquad \times
\prod_{i=1,i\not=k}^3
\int_{A_i-\Delta A_i}^{A_i+\Delta A_i}
\frac{\Delta A_i-|a_i-A_i|}{\Delta A_i}
\left\{\begin{matrix} 
P_m^{(0,n+1/2)}\left(I(a_i)\right) \left(\frac{a_i}{\mathbf{R}}\right)^{n}, &i=1\\ 
\mathrm{Trig}(j a_i), &i=2\\ 
P_{n,|j|}(a_i),&i=3,k=1\\
\frac{1}{1-a_i^2} P_{n,|j|}(a_i),&i=3,k=2
\end{matrix}\right\}
\intd a_i.
\end{align}
\end{landscape}
Note that, for practical purposes, we still have to compute 8 different integrals for this term. Moreover, the integrals with respect to $a_1=r$ and $a_3=t$ can still only be computed via numerical integration. For the case $a_2=\lon$, we obtain 
\begin{align}
\int_{\Phi-\Delta\Phi}^{\Phi+\Delta\Phi} -\frac{\sgn(\lon-\Phi)}{\Delta \Phi} j\mathrm{Trig}(-j\lon) \intd \lon
&=
\int_{\Phi-\Delta\Phi}^{\Phi} \frac{1}{\Delta \Phi} j\mathrm{Trig}(-j\lon) \intd \lon
-\int_{\Phi}^{\Phi+\Delta\Phi} \frac{1}{\Delta \Phi} j\mathrm{Trig}(-j\lon) \intd \lon\\ 
&=
 \frac{j}{\Delta \Phi} \left(\int_{\Phi-\Delta\Phi}^{\Phi} \mathrm{Trig}(-j\lon) \intd \lon
-\int_{\Phi}^{\Phi+\Delta\Phi} \mathrm{Trig}(-j\lon) \intd \lon\right)\\
&= 
 \frac{j}{\Delta \Phi} \left(2I_2(-j;\Phi) - I_2(-j;\Phi-\Delta \Phi) - I_2(-j;\Phi+\Delta \Phi)\right)
\intertext{and}
&\int_{\Phi-\Delta\Phi}^{\Phi+\Delta\Phi} \frac{\Delta \Phi-|\lon-\Phi|}{\Delta \Phi} \mathrm{Trig}(j\lon) \intd \lon,
\end{align}
which we have already discussed above.

\subsection{Derivation of objective functions $\RFMP(\cdot;\cdot)$ and $\ROFMP(\cdot;\cdot)$ and related coefficients}
\label{sect:app:OF_IPMPs}
We start with RFMP. The respective noise-cognizant Tikhonov-Phillips functional is given in \cref{eq:TFNO} by
\begin{align}
\mathcal{J}^{\mathrm{SM}} \left( f_N + \alpha d; \T_\daleth, \lambda, \delta \psi, \sigma\right) &\coloneqq 
\left\| \frac{R^N - \alpha\T_\daleth d}{\sigma} \right\|^2_{\real^\ell} + \lambda\left\|f_N+\alpha d\right\|^2_{\Hs{1}},\qquad \lambda >0,\\
R^{N+1} &\coloneqq R^N - \alpha_{N+1} \T_\daleth d_{N+1} = \delta \psi - \T_\daleth f_{N+1}.
\end{align}
We now aim to determine $(\alpha,d),\ \alpha \in \real,\ d\in\dic$ which minimizes $\mathcal{J}^{\mathrm{SM}}$. We start as follows: 
\begin{align}
0 &= \frac{\partial}{\partial \alpha} \mathcal{J}^{\mathrm{SM}} \left( f_N + \alpha d; \T_\daleth, \lambda, \delta \psi, \sigma\right)
= \frac{\partial}{\partial \alpha} \left[ \left\| \frac{R^N - \alpha\T_\daleth d}{\sigma} \right\|^2_{\real^\ell} + \lambda\left\|f_N+\alpha d\right\|^2_{\Hs{1}}\right]
\end{align}
\begin{align}
&= \frac{\partial}{\partial \alpha} \left[ 
\left\| \frac{R^N}{\sigma} \right\|^2_{\real^\ell} 
-2\alpha\left\langle \frac{R^N}{\sigma},\frac{\T_\daleth d}{\sigma} \right\rangle^2_{\real^\ell} 
+ \alpha^2\left\| \frac{\T_\daleth d}{\sigma} \right\|^2_{\real^\ell} 
\right. \\ &\qquad\qquad \left.
+ \lambda\left\|f_N\right\|^2_{\Hs{1}}
+ 2\alpha\lambda\left\langle f_N,d\right\rangle^2_{\Hs{1}}
+\alpha^2\lambda\left\|d\right\|^2_{\Hs{1}}
\right]\\
&=  
-2\left\langle \frac{R^N}{\sigma},\frac{\T_\daleth d}{\sigma} \right\rangle^2_{\real^\ell} 
+ 2\alpha\left\| \frac{\T_\daleth d}{\sigma} \right\|^2_{\real^\ell} 
+ 2\lambda\left\langle f_N,d\right\rangle^2_{\Hs{1}}
+2\alpha\lambda\left\|d\right\|^2_{\Hs{1}}.
\end{align}
This yields 
\begin{align}
\alpha_{N+1} = \frac{\left\langle \frac{R^N}{\sigma},\frac{\T_\daleth d}{\sigma} \right\rangle^2_{\real^\ell}-\lambda\left\langle f_N,d\right\rangle^2_{\Hs{1}}}{\left\| \frac{\T_\daleth d}{\sigma} \right\|^2_{\real^\ell}+\lambda\left\|d\right\|^2_{\Hs{1}}}.
\end{align}
Inserting this value into $\mathcal{J}^{\mathrm{SM}} \left( f_N + \alpha d; \T_\daleth, \lambda, \delta \psi, \sigma\right)$, we obtain 
\begin{align}
&\mathcal{J}^{\mathrm{SM}} \left( f_N + \alpha_{N+1} d; \T_\daleth, \lambda, \delta \psi, \sigma\right)\\
&= \left\| \frac{R^N}{\sigma} \right\|^2_{\real^\ell} 
-2\alpha_{N+1}\left\langle \frac{R^N}{\sigma},\frac{\T_\daleth d}{\sigma} \right\rangle^2_{\real^\ell} 
+ \alpha_{N+1}^2\left\| \frac{\T_\daleth d}{\sigma} \right\|^2_{\real^\ell} \\ 
&\qquad\qquad + \lambda\left\|f_N\right\|^2_{\Hs{1}}
+ 2\alpha_{N+1}\lambda\left\langle f_N,d\right\rangle^2_{\Hs{1}}
+\alpha_{N+1}^2\lambda\left\|d\right\|^2_{\Hs{1}}\\
&= \left\| \frac{R^N}{\sigma} \right\|^2_{\real^\ell} 
+ \lambda\left\|f_N\right\|^2_{\Hs{1}}\\
&\qquad\qquad -2\frac{\left\langle \frac{R^N}{\sigma},\frac{\T_\daleth d}{\sigma} \right\rangle^2_{\real^\ell}-\lambda\left\langle f_N,d\right\rangle^2_{\Hs{1}}}{\left\| \frac{\T_\daleth d}{\sigma} \right\|^2_{\real^\ell}+\lambda\left\|d\right\|^2_{\Hs{1}}} \left[\left\langle \frac{R^N}{\sigma},\frac{\T_\daleth d}{\sigma} \right\rangle^2_{\real^\ell} - \lambda\left\langle f_N,d\right\rangle^2_{\Hs{1}} \right]\\
&\qquad\qquad +\left[ \frac{\left\langle \frac{R^N}{\sigma},\frac{\T_\daleth d}{\sigma} \right\rangle^2_{\real^\ell}-\lambda\left\langle f_N,d\right\rangle^2_{\Hs{1}}}{\left\| \frac{\T_\daleth d}{\sigma} \right\|^2_{\real^\ell}+\lambda\left\|d\right\|^2_{\Hs{1}}}\right]^2 \left[\left\| \frac{\T_\daleth d}{\sigma} \right\|^2_{\real^\ell} + \lambda\left\|d\right\|^2_{\Hs{1}} \right]\\
&= \left\| \frac{R^N}{\sigma} \right\|^2_{\real^\ell} 
+ \lambda\left\|f_N\right\|^2_{\Hs{1}} 
-\frac{\left(\left\langle \frac{R^N}{\sigma},\frac{\T_\daleth d}{\sigma} \right\rangle^2_{\real^\ell}-\lambda\left\langle f_N,d\right\rangle^2_{\Hs{1}}\right)^2}{\left\| \frac{\T_\daleth d}{\sigma} \right\|^2_{\real^\ell}+\lambda\left\|d\right\|^2_{\Hs{1}}} \\
&= \mathcal{J}^{\mathrm{SM}} \left( f_{N-1} + \alpha_N d_N; \T_\daleth, \lambda, \delta \psi, \sigma\right)
- \RFMP(d;N). \\
\end{align}
Note that $\mathcal{J}^{\mathrm{SM}} \left( f_{N-1} + \alpha_N d_N; \T_\daleth, \lambda, \delta \psi, \sigma\right)$ is fixed in the $(N+1)$-th iteration. Thus, we see that, if we maximize $\RFMP(\cdot;\cdot)$ as defined in \cref{def:RFMP(d;N)}, we minimize the noise-cognizant Tikkhonov-Phillips functional in the $(N+1)$-th iteration.

Similarly, in the ROFMP, we consider the noise-cognizant Tikhonov-Phillips functional
\begin{align}
\mathcal{J}^{\mathrm{SM}}_O \left( f_{N}^{(N)} + \alpha d; \T_\daleth, \lambda, \delta \psi, \sigma\right) &\eqqcolon \left\| \frac{R^N - \alpha\proj{\mathcal{V}_N^\perp}\T_\daleth d}{\sigma} \right\|^2_{\real^\ell} + \lambda\left\|f_N^{(N)}+\alpha \left(d-b_n^{(N)}(d) \right)\right\|^2_{\Hs{1}},\\
R^{N+1} &\coloneqq R^N - \alpha_{N+1}^{(N+1)}\proj{\mathcal{V}_N^\perp} \T_\daleth d_{N+1} = \delta \psi - \T_\daleth f_{N+1}^{(N+1)},
\end{align}
for $\lambda >0$. Confer \cref{ssect:ipmps} for a full definition. Again, we first consider the minimizer with respect to $\alpha$:
\begin{align}
0 &= \frac{\partial}{\partial \alpha} \mathcal{J}^{\mathrm{SM}}_O \left( f_{N}^{(N)} + \alpha d; \T_\daleth, \lambda, \delta \psi, \sigma\right) \\
&= \frac{\partial}{\partial \alpha}\left[
\left\| \frac{R^N}{\sigma} \right\|^2_{\real^\ell} 
-2\alpha \left\langle \frac{R^N}{\sigma},\frac{\proj{\mathcal{V}_N^\perp}\T_\daleth d}{\sigma} \right\rangle^2_{\real^\ell} 
+\alpha^2\left\| \frac{\proj{\mathcal{V}_N^\perp}\T_\daleth d}{\sigma} \right\|^2_{\real^\ell} \right. \\
&\qquad\qquad \left. 
+ \lambda\left\|f_N^{(N)}\right\|^2_{\Hs{1}}
+ 2\alpha\lambda\left\langle f_N^{(N)},d-b_n^{(N)}(d) \right\rangle^2_{\Hs{1}}
+ \alpha^2\lambda\left\|d-b_n^{(N)}(d)\right\|^2_{\Hs{1}}
\right]\\
&=  
-2\left\langle \frac{R^N}{\sigma},\frac{\proj{\mathcal{V}_N^\perp}\T_\daleth d}{\sigma} \right\rangle^2_{\real^\ell} 
+2\alpha\left\| \frac{\proj{\mathcal{V}_N^\perp}\T_\daleth d}{\sigma} \right\|^2_{\real^\ell} \\
&\qquad\qquad 
+ 2\lambda\left\langle f_N^{(N)},d-b_n^{(N)}(d) \right\rangle^2_{\Hs{1}}
+ 2\alpha\lambda\left\|d-b_n^{(N)}(d)\right\|^2_{\Hs{1}},
\end{align}
which yields 
\begin{align}
\alpha_{N+1} = \alpha_{N+1}^{(N+1)} \coloneqq 
\frac{\left\langle \frac{R^N}{\sigma},\frac{\proj{\mathcal{V}_N^\perp}\T_\daleth d}{\sigma} \right\rangle^2_{\real^\ell} - \lambda\left\langle f_N^{(N)},d-b_n^{(N)}(d) \right\rangle^2_{\Hs{1}}}{\left\| \frac{\proj{\mathcal{V}_N^\perp}\T_\daleth d}{\sigma} \right\|^2_{\real^\ell} + \lambda\left\|d-b_n^{(N)}(d)\right\|^2_{\Hs{1}}}. 
\end{align}
Inserting $\alpha_{N+1}^{(N+1)}$ in $\mathcal{J}^{\mathrm{SM}}_O \left( f_{N}^{(N)} + \alpha d; \T_\daleth, \lambda, \delta \psi, \sigma\right)$, we obtain 
\begin{align}
&\mathcal{J}^{\mathrm{SM}}_O \left( f_{N}^{(N)} + \alpha_{N+1}^{(N+1)} d; \T_\daleth, \lambda, \delta \psi, \sigma\right)\\
&= 
\left\| \frac{R^N}{\sigma} \right\|^2_{\real^\ell} 
-2\alpha_{N+1}^{(N+1)} \left\langle \frac{R^N}{\sigma},\frac{\proj{\mathcal{V}_N^\perp}\T_\daleth d}{\sigma} \right\rangle^2_{\real^\ell} 
+\left(\alpha_{N+1}^{(N+1)}\right)^2\left\| \frac{\proj{\mathcal{V}_N^\perp}\T_\daleth d}{\sigma} \right\|^2_{\real^\ell} \\
&\qquad 
+ \lambda\left\|f_N^{(N)}\right\|^2_{\Hs{1}}
+ 2\alpha_{N+1}^{(N+1)}\lambda\left\langle f_N^{(N)},d-b_n^{(N)}(d) \right\rangle^2_{\Hs{1}}
+ \left(\alpha_{N+1}^{(N+1)}\right)^2\lambda\left\|d-b_n^{(N)}(d)\right\|^2_{\Hs{1}}\\
&= 
\left\| \frac{R^N}{\sigma} \right\|^2_{\real^\ell}
+ \lambda\left\|f_N^{(N)}\right\|^2_{\Hs{1}}\\ 
&\qquad -2
\frac{\left\langle \frac{R^N}{\sigma},\frac{\proj{\mathcal{V}_N^\perp}\T_\daleth d}{\sigma} \right\rangle^2_{\real^\ell} - \lambda\left\langle f_N^{(N)},d-b_n^{(N)}(d) \right\rangle^2_{\Hs{1}}}{\left\| \frac{\proj{\mathcal{V}_N^\perp}\T_\daleth d}{\sigma} \right\|^2_{\real^\ell} + \lambda\left\|d-b_n^{(N)}(d)\right\|^2_{\Hs{1}}} \\
&\qquad\qquad\qquad\qquad \times \left[\left\langle \frac{R^N}{\sigma},\frac{\proj{\mathcal{V}_N^\perp}\T_\daleth d}{\sigma} \right\rangle^2_{\real^\ell} + \lambda\left\langle f_N^{(N)},d-b_n^{(N)}(d) \right\rangle^2_{\Hs{1}} \right]\\ 
&\qquad +\left(\frac{\left\langle \frac{R^N}{\sigma},\frac{\proj{\mathcal{V}_N^\perp}\T_\daleth d}{\sigma} \right\rangle^2_{\real^\ell} - \lambda\left\langle f_N^{(N)},d-b_n^{(N)}(d) \right\rangle^2_{\Hs{1}}}{\left\| \frac{\proj{\mathcal{V}_N^\perp}\T_\daleth d}{\sigma} \right\|^2_{\real^\ell} + \lambda\left\|d-b_n^{(N)}(d)\right\|^2_{\Hs{1}}}\right)^2\\
&\qquad\qquad\qquad\qquad \times 
\left[ \left\| \frac{\proj{\mathcal{V}_N^\perp}\T_\daleth d}{\sigma} \right\|^2_{\real^\ell} +\lambda\left\|d-b_n^{(N)}(d)\right\|^2_{\Hs{1}}\right] 
\end{align}
\begin{align}
&= 
\left\| \frac{R^N}{\sigma} \right\|^2_{\real^\ell}
+ \lambda\left\|f_N^{(N)}\right\|^2_{\Hs{1}}
- \frac{\left(\left\langle \frac{R^N}{\sigma},\frac{\proj{\mathcal{V}_N^\perp}\T_\daleth d}{\sigma} \right\rangle^2_{\real^\ell} - \lambda\left\langle f_N^{(N)},d-b_n^{(N)}(d) \right\rangle^2_{\Hs{1}}\right)^2}{\left\| \frac{\proj{\mathcal{V}_N^\perp}\T_\daleth d}{\sigma} \right\|^2_{\real^\ell} + \lambda\left\|d-b_n^{(N)}(d)\right\|^2_{\Hs{1}}}\\
&= 
\mathcal{J}^{\mathrm{SM}}_O \left( f_{N-1}^{(N-1)} + \alpha_{N+1}^{(N+1)} d; \T_\daleth, \lambda, \delta \psi, \sigma\right) - \ROFMP(d;N).
\end{align}

\end{document}